\documentclass[journal]{IEEEtran}
\usepackage{cite}

\ifCLASSINFOpdf
  \usepackage[pdftex]{graphicx}
  \graphicspath{{./pictsUsedOnly/}}
\else
\fi
\usepackage[cmex10]{amsmath}
\usepackage{algorithm,algorithmic}

\usepackage{array}

\ifCLASSOPTIONcompsoc
  \usepackage[caption=false,font=normalsize,labelfont=sf,textfont=sf]{subfig}
\else
  \usepackage[caption=false,font=footnotesize]{subfig}
\fi

\usepackage{siunitx}

\newcommand{\tr}{ ^{\rm T}}
\newcommand{\inv}{ ^{-1}}

\DeclareMathOperator*{\argmin}{arg\,min}

\newcommand{\ft}[1]{ \mathcal{F}\left(#1\right)}
\newcommand{\ift}[1]{\mathcal{F}^{-1}\left(#1\right)}
\newcommand{\fts}[1]{ \mathcal{F}_{\parallel}\left(#1\right)}
\newcommand{\ifts}[1]{\mathcal{F}_{\parallel}^{-1}\left(#1\right)}
\newcommand{\ftz}[1]{ \mathcal{F}_{\perp}\left(#1\right)}
\newcommand{\iftz}[1]{\mathcal{F}_{\perp}^{-1}\left(#1\right)}
\newcommand{\ftt}[1]{ \mathcal{F}_{\rm t}\left(#1\right)}

\newcommand{\bk}{ {\bf k}}
\newcommand{\bx}{ {\bf x}}

\newcommand{\mc}[1]{ \mathcal{#1}}

\def\kwave~{\textbf{k}-\texttt{Wave}}

\newcommand{\dontshow}[1]{}
\newcommand{\query}[1]{{\color{green}#1}}

\usepackage{color}
\usepackage{amssymb}
\usepackage{amsfonts}
\usepackage{amsthm}

\newtheorem{innercustomthm}{Theorem}
\newenvironment{theorem}[1]
  {\renewcommand\theinnercustomthm{#1}\innercustomthm}
  {\endinnercustomthm}

\definecolor{darkblue}{rgb}{0,0,0.6}
\definecolor{darkgreen}{rgb}{0,0.6,0}
\definecolor{darkred}{rgb}{0.6,0,0}
\definecolor{orange}{rgb}{1,0.7,0}
\definecolor{darkorange}{rgb}{1,0.4,0.3}

\begin{document}
%
\title{Acoustic Wave Field Reconstruction from Compressed Measurements with Application in Photoacoustic Tomography}
%
%
%

\author{Marta M.~Betcke, Ben T.~Cox, Nam Huynh, Edward Z.~Zhang, Paul C.~Beard and Simon R.~Arridge
\thanks{M.~Betcke and S.~Arridge are with the Department
of  Computer Science, University College London, UK,
WC1E 6BT London, UK, e-mail: m.betcke@ucl.ac.uk.}
\thanks{B.~Cox and N.~Huynh and E.~Zhang and P.~Beard are with Department of Medical Physics, University College London, UK,
WC1E 6BT London, UK.}
}

%
%

\markboth{Journal of \LaTeX\ Class Files,~Vol.~11, No.~4, December~2012}%
{Shell \MakeLowercase{\textit{et al.}}: Bare Demo of IEEEtran.cls for Journals}
%



\maketitle

\begin{abstract}

We present a method for the recovery of compressively sensed acoustic fields using patterned, instead of point-by-point, detection.
From a limited number of such compressed measurements, we propose to reconstruct the field on the sensor plane in each time step independently assuming its sparsity in a Curvelet frame.
A modification of the Curvelet frame is proposed to account for the smoothing effects of data acquisition and motivated by a frequency domain model for photoacoustic tomography. An ADMM type algorithm, SALSA, is used to recover the pointwise data in each individual time step from the patterned measurements. For photoacoustic applications, the photoacoustic image of the initial pressure is reconstructed using time reversal in $\bk$-Wave Toolbox.  

\end{abstract}

\begin{IEEEkeywords}
ADMM methods, compressed sensing, curvelet frame, L1 minimization, photoacoustic tomography.
\end{IEEEkeywords}

%
\IEEEpeerreviewmaketitle

\section{Introduction}
%
%
%
%

\IEEEPARstart{C}{ompressed} Sensing (CS) is a new measurement paradigm, which allows for the reconstruction of sparse signals sampled at sub-Nyquist rates. Nowadays, it is common understanding that many digital signals and images admit an adequate representation with far fewer coefficients than their actual length. This phenomena is known as compressibility and it has been a driving force in many image processing applications, most notably the image compression algorithms JPEG and its successor JPEG 2000.  
CS emerged from the realization that the signals could in fact be acquired directly in their compressed form instead of sampling the signal at Nyquist rate and then compressing it and while doing so discarding most of the laboriously obtained coefficients. 

Since the seminal works by Donoho \cite{Donoho:2006cs} and C\`andes, Romberg and Tao \cite{Candes:2006rup,Candes:2006ssr}, there has been an explosion of results in the field and many applications has been suggested starting with the prototype single pixel camera \cite{Duarte:2008ricecam}.  

Here we propose a new application of CS in acoustic wave field sensing. In \cite{Candes:2005coswp} the authors proved that the acoustic field is almost optimally sparse in Curvelet frame. 
As the Curvelets essentially describe the wave front sets their propagation is well approximated through geometrical optics (high frequency asymptotic solution to the wave equation).  The wave front sets, and hence the Curvelets, propagate along the geometrical rays which are projections on $\mathbb R^d$ of phase space solutions of the corresponding Hamilton-Jacobi equation resulting from the high frequency asymptotic.
Motivated by this result we investigate the Curvelet representation of the cross-section of the wave field by the planar ultrasound sensor. While the arguments in \cite{Candes:2005coswp} do not directly apply to this situation, the planar cross-section through the acoustic wave front constitutes a singularity along a smooth curve for which Curvelets have been demonstrated to be a nearly optimal representation \cite{Candes:2004tfc}.

\subsection{Contribution}
In this paper we focus on the reconstruction problem via data recovery for a novel way of interrogating the high resolution 
ultrasound sensor using patterns instead of the more conventional sequential point-by-point interrogation. An example of such a system using a Single-Pixel Optical Camera (SPOC) was presented in \cite{Huynh:2016spoc}.
The theory of CS predicts that substantially fewer such measurements, of the order $k \log n$, need to be taken in order to capture a signal of length $n$ and sparsity $k$ leading to a substantial reduction of the acquisition time. We propose from such measurements to recover the pressure at the detector at each time step independently using CS recovery algorithms. We discuss the motivation for using a Curvelet tranform as the sparsifying transformation for the acoustic field at the detector and we propose its modification: a low-frequency Curvelet tranform tailored to the frequency range of the acoustic field on the detector. We investigate the appropriate choice of interrogation patterns and recovery algorithm for this problem. The proposed techniques can be used in various applications ranging from ultrasonic field mapping to photoacoustic tomography (PAT). In this work we focus on the latter to illustrate our method. 

For PAT applications, the immediate benefit of recovering data independently at each time step
is the decoupling of the CS reconstruction from the acoustic inversion, allowing for the recovered time series to be input into any
reconstruction algorithm for PAT for which software is readily available. 
\dontshow{This in turn allows to easily incorporate effects like acoustic absorption (nonlinear problem), or varying speed of sound in the medium (does not result in spherical means transform) which require (possibly nonlinear) numerical solve.}

\subsection{Related work}
For PAT, one and two step approaches have been investigated in the literature. In one step approaches the initial pressure is directly recovered from the compressed measurements. In \cite{Provost:2009cspat} the authors propose such an approach for a 2D problem with measurements limited in angle and frequency. It was numerically investigated using an analytic approximate inversion formula and a number of sparse representations including Wavelets and Curvelets. An approach for 3D imaging using a translatable circular detection array and Wavelets as a sparse representation was presented in \cite{Guo:2011}. A generic variational approach for 2D and 3D PAT reconstruction based on an algebraic adjoint was first discussed in \cite{Huang:2013fwi}. Two variational approaches based on the analytical adjoint were recently proposed, one in a FEM-BEM setting \cite{Belhachmi:2016adj} and the other in a $\bk$-space setting with an efficient $\bk$-Wave implementation \cite{Arridge:2016adj,Arridge:2016apat}. 

In contrast, in the context of PAT, the present paper is about a two step method. We propose to first recover the photoacoustic data in every time step independently from pattern measurements using sparsity of the data in Curvelet basis, and subsequently to reconstruct it using standard PAT reconstruction methods. 
Recently, another two stage approach has been suggested in \cite{Haltmeier:2016dsp}, where the authors explore the temporal sparsity of the data by means of a custom made transform in time. From the sparse pressure data subsequently the initial pressure is reconstructed with the Universal Back Projection formula.  

\subsection{Outline}
The remainder of this paper is organized as follows. In Section \ref{sec:pat} we introduce the forward and inverse problem for PAT and methods for their solution. Section \ref{sec:cs} briefly recapitulates the theory of compressed sensing. An appropriate multiscale representation of the time series PAT data is considered in Section \ref{sec:ms}, where we derive the frequency model of sensor data and propose a modified version of Curvelet transform tailored to the range of frequencies of the acoustic field.  
In Section \ref{sec:cspat} we discuss specific issues arising when compressively sensing photoacoustic data using patterned interrogation of optical ultrasound detector. We briefly describe the Single-Pixel Optical Camera based PAT scanner. We consider the challenges for Curvelet transforms for approximation of sensor data over the time series. We discuss choice of the interrogation patterns and the algorithm for recovery of the sensor data. In Section \ref{sec:recon} we present recovery results for the optical sensor data and the final PAT image reconstruction from both simulated and real data.

\section{Photoacoustic Tomography}\label{sec:pat}

Photoacoustic tomography (PAT) is an example of a wider range of hybrid imaging techniques, in which contrast induced by one type of wave is read out by another wave. In this way, both high contrast and high resolution can be simultaneously achieved, which is often difficult with conventional imaging techniques that usually provide either 
one or the other, but not both. 
PAT is an emerging biomedical imaging modality with both pre-clinical and clinical applications that can provide complementary information to established imaging techniques 
\cite{Wang:2009pas,Beard:2011if,Valluru:2016patonc,Zhou:2016tutpat,Xia:2015sawb}


Many PAT applications require a high resolution, three dimensional image e.g.~an image of capillaries of a few tens of microns diameter in a cm sized imaging region. Such highly resolved imaging requires an ultrasound sensor array of tens of thousands of pixels.
In one such PAT system \cite{Zhang:2008backm}, the sensor is a Fabry Perot (FP) interferometer interrogated by a laser whose focus is moved to form a raster scan of the desired resolution. For sequential sampling, such as this, the ultimate limit to the data acquisition rate (the rate at which time series are collected) is the propagation time for sound to cross the specimen, e.g.~it would take $10\mu s$ for the signal to reach the detector from 15mm depth, resulting in $100 kHz$ acquisition rate (not to be confused with the sampling rate, which might be as high as 100 MHz). No sequential scanner is close to approaching this limit, making a sequential acquisition a major practical limitation for high resolution 3D PAT.
For \textit{in vivo} applications the required acquisition times at currently achievable rates are typically a few minutes, not only resulting in motion artifacts, but limiting studies to phenomena on such long timescales.


The principle involved in PAT is to send a short (ns) pulse of near-infrared or visible light into tissue, whereupon absorption of the photons e.g.\ by haemoglobin molecules, generates a small local increase in pressure which propagates to the surface as a broadband, ultrasonic pulse. If the amplitude of this signal is recorded over an array of sensors at the tissue surface, an image reconstruction algorithm can be used to estimate the original 3D pressure increase due to optical absorption; this is the photoacoustic image, $p_0$. 

Mathematically, under an assumption of free space propagation the photoacoustic forward problem is modelled as an initial value problem for the wave equation \cite{TrJaReCo12}
\begin{subequations}\label{eq:pafp}
\begin{align}
\frac{1}{c^2(\bx)}\frac{\partial^2 p(\bx,t)}{\partial t^2}  &= \rho_0(\bx)\nabla\cdot \left(\frac{1}{\rho_0(\bx)} \nabla \right)  p(\bx,t),\, \bx \in \mathbb R^d, t \in (0,T), \\
p(\bx,0) &= p_0(\bx), \\
\frac {\partial}{\partial t} p(\bx,0) &= 0,
\end{align}
\end{subequations}
where $p(\bx, t)$ denotes the time dependent acoustic pressure in $\mathbb R^d \times (0,T), d=2,3$, $p_0(\bf x)$ its initial value and $c(\bx)$ and $\rho_0(\bx)$ are the ambient speed of sound and density, respectively. 

The photoacoustic inverse problem is to recover this initial pressure $p_0(\bx), \bx \in \Omega$ compactly supported in the region of interest  $\Omega$ from a time series measurement $g(\bx,t) = p(\bx, t), \bx \in \mc S, t \in (0,T)$  on the surface $\mc S$ (e.g.~boundary of $\Omega$) and it amounts to a solution of the following initial boundary value problem \cite{Finch:2004msph}
\begin{subequations}\label{eq:paip}
\begin{align}
\frac{1}{c_0^2(\bx)}\frac{\partial^2 p(\bx,t)}{\partial t^2}  &= \rho_0(\bx)\nabla\cdot \left(\frac{1}{\rho_0(\bx)} \nabla \right)  p(\bx,t) ,\; \bx \in \Omega , t \in (0,T),  \\
p(\bx,0) &= \mathbf{0}, \\
\frac {\partial}{\partial t} p(\bx,0) &= 0,\\
p(\bx, t) &= g(\bx, T-t), \quad \bx \in \mc S, \,  t \in (0, T),
\end{align}
\end{subequations}
with PAT data fed backwards in time as boundary values (also referred to as \textit{time reversal}). The formulation \eqref{eq:paip} holds exactly in 3D for non-trapping smooth sound speed $c(x)$ if $T$ has been chosen large enough so that $g(\bx_{\mc S},t)=0,\, t \geq T$ and the wave has left the domain $\Omega$.
Furthermore, assuming that the measurement surface $\mc S$ surrounds the region of interest $\Omega$ containing the support of the initial pressure $p_0$, problem \eqref{eq:paip} has a unique solution. 
The condition on $\mc S$ to surround $p_0$ can be relaxed under additional assumptions on $\mc S$ and $\Omega$ and smoothess of initial pressure $p_0 \in H_0^1(\Omega)$; see \cite{Kuchment:2011mtat} and the citations within. The stability of reconstruction can be obtained from microlocal analysis which provides insights into which singularities are visible and which are not \cite{Stefanov:2009tatvs,Frikel:2015mlapat}, which directly translates to stability of reconstruction of those singularities \cite{Kuchment:2011mtat}. 

In a real experiment, we can only acquire a discrete (both in time and space) subset of time series measurements, and frequently it is not possible to acquire measurements on a surface $\mc S$ surrounding the object. An example of a popular sensor violating this assumption is a planar sensor. In practice, such sensor will have finite size, resulting in invisibility of some interfaces and in turn in artefacts in the reconstructed image. 

We should mention that other approaches to image reconstruction exist. An overview of methods for the case when $\mathcal S$ is a surface surrounding $\Omega$ can be found e.g.~in \cite{Kuchment:2011mtat}.
In this paper we directly solve the time-reversal approach \eqref{eq:paip} using the pseudospectral method implemented in the $\bk$-Wave Toolbox \cite{Treeby:2010kwave}, which is an efficient numerical scheme 
for solving the wave equation in domains with heterogeneous acoustic properties, and is exact in the case of homogeneous media.
Furthermore, the methodology proposed here is tailored to high resolution detectors 
which are planar 
and so we assume $\mathcal S$ to be a finite rectangular section of the $xy$-plane
$$\mathcal S = \{ (x,y,z): \; |x| \leq x_d/2, \,  |y| \leq y_d/2, \, z =0 \}.$$
and solve the resulting initial \emph{partial} boundary value problem \eqref{eq:paip}. 
The rectangular planar detector shape allows us to use any sparsifying transform derived for natural images.

\section{Compressed sensing}\label{sec:cs}
Let $\Psi$ be a sparsifying transform 
$$\Psi: \mathbb{R}^n \rightarrow \mathbb{R}^N, $$ 
resulting in possibly an overdetermined representation $N \geq n$ and 
$\Psi\inv$ its inverse (in the frame sense). 
With $f$ we  denote the transformation of the original signal $g \in \mathbb{R}^n$
\begin{equation}
f = \Psi g. 
\end{equation}
In compressed sensing the signal $g \in \mathbb{R}^n$ is projected on a series of sensing vectors $\phi_i \in \mathbb{R}^n, i \in 1,\dots,m$, where $m \ll n$, yielding a vector $b \in \mathbb{R}^m$ of compressive measurements
\begin{equation}\label{eq:compm}
b = \Phi g + e \quad \mbox{with} \quad \Phi = [\phi_1, \phi_2, \dots, \phi_m]\tr,
\end{equation}
where $e$ is the measurement noise, $\| e\|_2 \leq \varepsilon$. 
The following recovery result guarantees that the original signal $g$, can be robustly recovered from compressed measurements \eqref{eq:compm} via solution of the 
minimization problem for $f$ \cite{Candes:2006ssr}
\begin{equation}\label{eq:l1rec}
\min_{f \in \mathbb R^n} \| f \|_1, \quad \mbox{s.t.}  \quad  \| \Phi \Psi\inv f - b \| \leq \varepsilon.
\end{equation}

\begin{theorem}{RR}[Robust recovery \cite{Candes:2008rip, Fornasier:2011cs}]\label{thm:l1rec}
Let $\delta_{k}$ be the isometry constant of $\Phi\Psi\inv$ defined as a smallest positive number such that
\begin{equation}
(1 -\delta_k) \| f \|_2^2 \leq  \| \Phi \Psi\inv f\|_2^2 \leq(1 + \delta_k) \| f \|_2^2 
\end{equation}
holds for all $k$-sparse vectors $f$. 

If $\delta_{2k} < \sqrt{2} - 1$ (relaxed to $\delta_{2k} < \frac{2}{3 + \sqrt{7/4}}$ in \cite{Fornasier:2011cs}) then the error of solution to \eqref{eq:l1rec} $f^*$ is bounded as follows
\begin{equation}\label{eq:l1err}
\| f - f^*\|_2 \leq C_1 \varepsilon + C_2 \frac{\| f - f^k \|_1}{\sqrt{k}},
\end{equation}
where $f^k$ denotes best $k$-term approximation, obtained from $f$ selecting its $k$ largest in magnitude coefficients, and $C_1, C_2$ are constants dependent only on $\delta_{2k}$. 
\end{theorem}
The robust recovery Theorem \ref{thm:l1rec} holds for any vector $f$, however the error norm $\| f - f^k\|_1$ of the $k$-term approximation is only small for $k$-sparse or compressible vector $f$  (i.e.~with fast enough decaying magnitude of the coefficients). Furthermore, for $k$-sparse vectors we have $ (\| f - f^k \|_1)/\sqrt{k} \leq  \| f - f^k \|_2$ and hence the the bound can be expressed in $L_2$ norm
\begin{equation*}
\| f - f^*\|_2 \leq C_1 \varepsilon + C_2  \| f - f^k \|_2.
\end{equation*}



\section{Multiscale representation of time series data}\label{sec:ms}
\subsection{Frequency model of sensor data}\label{sec:ms:bli} 

For constant sound speed and density $c(\bx) = c_0, \rho(\bx) = \rho_0$, the wave equation \eqref{eq:pafp} becomes homogeneous and it admits an analytic solution in Fourier domain \cite{Cox:2005fwdkspace} 
\begin{equation}\label{eq:prop}
p(\bk,t) = \cos(c_0 |\bk| t) p_0(\bk),
\end{equation}
where $\bf k$ is the frequency domain wave vector. 
Equation \eqref{eq:prop} formally allows us to calculate the PAT data time series
\begin{eqnarray}
\label{eq:pxst}
p({\bf x_{\mathcal S}},t) &=&  \left.\mathcal F^{-1} (p( \bk, t)) \right|_{\bx = \bf{x_{\mathcal S}}} \\
\nonumber &=& \left. \ift{ \cos(c_0 |\bk| t) p_0(\bk) }  \right|_{\bx = \bf{x_{\mathcal S}}}\\
\nonumber &=& 2\int_0^\infty \ifts{ \cos(c_0 |\bk| t) p_0(\bk) } d\bk_{\perp},
\end{eqnarray}
where we used the tensor decomposition $\mathcal F = \mathcal F_{\parallel}\mathcal F_{\perp}$of the 3D Fourier transform in, and orthogonal to, the sensor plane, and the fact that 
for $\bx_{\mc S}$ we have $z=0$ 
and hence $e^{i\bk_{\perp} z} = 1$.  

Since $p_0(\bx)$ is a real valued function, $p_0(\bk)$ and consequently $p(\bk,t)$ are symplectic in $\bk$ i.e.~$p_0(\bk) = p_0(-\bk)^*$.  
Furthermore, $p(\bk,t)$ as given by \eqref{eq:prop} is even in $t$. As our domain $\Omega$ is positioned in $z\geq 0$ to restore uniqueness we assume $p_0(\bx)$ to be even w.r.t.~$z=0$ plane. Consequently, the Fourier transform in the $z$ as well as the $t$ direction is equal to the cosine transform and can be calculated integrating from 0 to $\infty$; thus in what follows the two are used interchangeably. 
 
\dontshow{ 
\textbf{Alternative 1:}\\ 
This model is however an unrealistic description of the acoustic pressure recorded by the sensor. For once, it is known that the acoustic absorption in the tissue follows the frequency power law \cite{Duck1990}
\begin{equation}
\alpha = \alpha_0 \omega^\beta,
\end{equation}
where $\alpha_0$ is acoustic absorption in [$({\rm rad/s})^{-y} {\rm m}^{-1}$], $\omega$ is angular frequency in [${\rm rad/s}$] and $\beta \in [1,2]$. \query{If I write it like that, why do I not solve the problem with absorption?} Further deterioration can occur through the interaction of the wave with the sensor. Finally, the pressure is read by an optical readout system, introducing spatial sampling frequency \query{and possibly a PSF?}. 

On the other hand the frequency roll off, as in Blackman window, is necessary to avoid Gibbs phenomena in the numerical methods \cite{Treeby:2011gf}, however the model tests performed in \query{citation?} suggest that it in fact resembles the behavior of the physical measurements. Hence, in what follows we adapt this model.

\textbf{Alternative 2:}\\
}

In experiments however, $p(\bx_{\mathcal S},t)$ is available to us only indirectly, through measurements. To account for limitations of physical equipment, we introduce a degradation operator, $\mathcal D$. Note, that here $\mathcal D$ is not a measurement operator, but it models the effect of the finite size and temporal response of the
measurement system which reads the acoustic time series.
In other words, $g(\bx_{\mathcal S},t) = \mathcal D p(\bx_{\mathcal S},t)$ describes the 
filtering 
of the acoustic pressure $p(\bx_{\mathcal S},t)$ 
by the physical presence of the sensor
before the sensing vector, $\phi_i$ is applied to it to collect the compressive measurement, $b_i$. 

In this work we assume that $\mathcal D$ is a band limiting spatial and temporal filtering operator acting on the time dependent pressure on optical sensor, which  for simplicity we describe in the frequency domain
\begin{eqnarray}
&\mc D[p](\bk_{\mathcal S}, \omega) = w_t(\omega) w_{\parallel}(\bk_{\mc S}) p(\bk_{\mc S},\omega),
\end{eqnarray}
where $w_{\parallel}(\bk_{\mc S})$ and $w_t(\omega)$ are 
some frequency window functions on the sensor and in time.

Taking Fourier transform of \eqref{eq:pxst} in both variables we obtain
\begin{eqnarray}
\nonumber &p(\bk_{\mc S}, \omega) = \ftt{\fts{p(\bx_{\mc S},t)}}\\ 
\nonumber &= \ftt{ 2\int_0^\infty \fts{  \ifts{ \cos(c_0|\bk|t) p_0(\bk) } } d\bk_{\perp} }\\
\nonumber &= \ftt{ 2\int_0^\infty  \frac{\omega/c_0^2 \cos(\omega t)}{\sqrt{(\omega/c_0)^2-|\bk_{\mc S}|^2}} p_0(\bk_{\mc S},\sqrt{(\omega/c_0)^2-|\bk_{\mc S}|^2}) d\omega }\\
\label{eq:pkso} &= \frac{\omega/c_0^2}{\sqrt{(\omega/c_0)^2-|\bk_{\mc S}|^2}} p_0(\bk_{\mc S},\sqrt{(\omega/c_0)^2-|\bk_{\mc S}^2|}),
\end{eqnarray}
where changing the integration variable to $\omega$, $\omega/c > |\bk_{\mc S}|$, in the third line allowed us to interpret the integral as inverse cosine transform in $\omega$. Equation \eqref{eq:pkso} connects the Fourier transform of the pressure time series on the detector with the Fourier transform of initial pressure and is the basis of the reconstruction formula derived in \cite{Koestli:2001f,Xu:2002efdr}.
From  \eqref{eq:pkso} it is obvious that application of $\mc D[p](\bk_{\mc S},\omega)$ corresponds to the application of a filter window $w_t(c_0|\bk|) w_{\parallel}(\bk_{\mc S})$ to $p_0(\bk)$ 
\begin{eqnarray*}
\mc D[p](\bk_{\mc S}, \omega) =& \frac{w_t(\omega) w_{\parallel}(\bk_{\mc S}) \omega/c_0^2}{\sqrt{(\omega/c_0)^2-|\bk_{\mc S}|^2}} p_0(\bk_{\mc S},\sqrt{(\omega/c_0)^2-|\bk_{\mc S}|^2})\\
=& w_t(c_0|\bk|) w_{\parallel}(\bk_{\mc S}) \frac{|\bk|}{|\bk_{\perp}|c_0} p_0(\bk) =: \mc D_\Omega[p_0](\bk)
\end{eqnarray*}
for scalar $\bk_{\perp} > 0$. Thus only a smoothed initial pressure can be recovered from the data. Conversely, when simulating PAT data, the initial pressure can be smoothed with $w_t(c_0|\bk|) w_{\parallel}(\bk_{\mc S})$ before forward propagation instead of smoothing the sensor data with $\mc D$, which can be useful to eliminate Gibbs phenomena in $\bk$-space methods; see also \cite{Tabei:2002kspace,Treeby:2011gf} where a Blackman window was applied to $p_0$.


\subsection{Curvelets}\label{sec:ms:c}
In this work we suggest using Curvelets to represent each of the measured time series data $g(\bx_{\mathcal S}, t)$. 
For the planar sensor $\mathcal S$ this corresponds to a planar cross-section of the wave field $p(\bx, t)$ at a given time $t$ (which we observe over a finite section of the plane, $\mc S$). 
Due to PAT forward problem being an initial value problem, the corresponding wave field $p(\bx, t)$ is essentially smooth away from the $p_0$-shaped (with corners smoothed) wave front and the same holds for their planar cross-sections.
As the time evolves $p(\bx,t)$ develops overlapping wave fronts. Those however can be treated as superposition and hence we only need to be able to represent an individual wave front.    

In what follows we assume $g(\bx_{\mathcal S}, t)$ to be sampled on an $n_1 \times n_2$ grid on the 
sensor, resulting in an image 
$g^t[i_1,i_2], \, i_1 = 1,\dots,n_1 \, i_2 = 1,\dots, n_2$,
where $(\cdot)$ denotes function evaluation and $[\cdot]$ array indexing (following the original  notation in \cite{Candes:2006fdct}).

Curvelets \cite{Candes:2004ntfc, Candes:2006fdct} are a multiscale pyramid with many directions and positions at each scale. Curvelets obey parabolic scaling, meaning that at scale $2^{-j}$ Curvelet has an envelope which aligns along a ridge of length $2^{-j/2}$ and width $2^{-j}$, resulting in a location, direction and scale dependent frequency plain tiling. 
Curvelets have been shown to be nearly optimal representation  of objects smooth away from (piecewise) $\mathcal C^2$ singularities \cite{Donoho:2001oar, Kashin:1987acos} and hence the error of the $k$-term curvelet approximation (corresponding to taking the $k$ largest magnitude coefficients) can be bounded as \cite{Candes:2004tfc} 
\begin{equation}
\|g - g^k\|^2_2 \leq C \cdot (\log k)^3 \cdot k^{-2}.
\end{equation}

\dontshow{
\begin{figure}[h]
\begin{center}
\includegraphics[width=0.3\textwidth]{corona}
\caption{Tiling of the frequency space induced by curvelets, courtesy \cite{Candes:2006fdct}.}
\label{fig:corona}
\end{center}
\end{figure}
}

Computation of Curvelets at the finest scale, $2^{-J}$, is not straight forward because in order to capture the direction of the wave it is necessary that the Curvelet is sampled more finely than the maximal frequency in the image. As a result the frequency domain support of a Curvelet at the finest scale $J$ with the direction $\theta_{\ell}$,  $\tilde U_{J,\ell}$, does not fit into the fundamental cell $[-n_1/2, n_1/2-1] \times [-n_2/2, n_2/2-1]$. 
One solution given in \cite{Candes:2006fdct} is to wrap the Fourier transform back onto the fundamental cell which effectively periodizes it 
\begin{eqnarray*}
& \tilde U_{J,\ell} [\left(i_1 + \frac{n_1}{2}\right) \bmod n_1 - \frac{n_1}{2} , \left( i_2 + \frac{n_2}{2} \right) \bmod n_2- \frac{n_2}{2}] \\
&= \tilde U_{J, \ell} (2\pi i_1, 2\pi i_2), 
\end{eqnarray*}
where the indices $i_{1,2}$ are chosen in the support of $\tilde U_{J, \ell}(\cdot,\cdot)$.
The periodization in Fourier domain corresponds to undersampling in space, which in turn causes aliasing. In \cite{Candes:2006fdct} this effect was shown to account for less than 10\% of the squared norm of the coefficients. 

After periodization, the Fourier transform is multiplied by a $\mathcal C^{\infty}$ partition of unity  window i.e.~a window which weights the original frequency and its periodic extension such that the sum of their squares is equal to 1. This acts to preserve the norm of the signal throughout periodization thereby guaranteeing that the computed transform is a numerical isometry. The $\mathcal C^{\infty}$ window used in \cite{Candes:2006fdct} is a tensor product of the following one dimensional $\mathcal C^{\infty}$ windows 
\begin{equation}\label{eq:cinfw}
w[i] = \left\{\begin{array}{ll}
a(x_i) h_{\infty}(1-x_i),  & i = 1,\dots, m_1,\\
1, &  i = m_1+1,\dots, 3m_1+1,\\
a(x_i) h_{\infty}(x_{i-3m-1}), & i = 3m_1+2,\dots, 4m_1+1,
\end{array}\right. 
\end{equation}
where $x_i = (i-1)/(m_1-1), \; i =1, \dots m_1$, $m_1 = \lfloor n_1/3 \rfloor$,
$h_{\infty}$ is a  $\mathcal C^{\infty}$ monotonically decreasing function 
and 
\begin{equation}\label{eq:afact}
a(x) = (h_{\infty}(1-x)^2 + h_{\infty}(x)^2)^{-1/2}
\end{equation}
is a normalizing factor. Expression \eqref{eq:cinfw} assumes $n_1 \bmod 3 =0$, and hence $m_1 = n_1/3$ resulting in a window of length $4m_1+1$, while corresponding windows can be derived in the cases $n_1 \bmod 3 = 1,2$.

\subsection{Low-frequency Curvelets}\label{sec:ms:lc}

The sampled time series data, $g^t$, 
is band limited with the high frequencies roll off  (see  Section \ref{sec:ms:bli}) in contrast to standard images which are ``sharp'' i.e.~high frequencies are not damped. 
This raises the question, if we should represent those rolling off frequencies in the same way as the undamped frequencies, as the standard Curvelets would do at a cost which scales with $n = n_1 n_2$ as $n \log n$. 

As a compromise, we propose here to use a \emph{low-frequency Curvelet} representation for the high resolution image 
$g^t$, i.e.~to compute Curvelets corresponding to a subset of the frequencies excluding the Cartesian annulus of the highest frequencies, $[-n^{LF}_1/2, n^{LF}_1 /2-1] \times [-n^{LF}_2/2,  n^{LF}_2 /2-1],  n^{LF}_l < n_l, l=1,2$.
The high frequencies in the annulus $[-n_1/2, n^{LF}_1 /2-1] \cup [n^{LF}_1 /2, n_1/2-1] \times [-n_2/2, n^{LF}_2 /2-1] \cup [n^{LF}_2 /2, n_2/2-1]$
present in the high resolution image $g^t$ can then be used to fill up the frequency range needed for computation of the low-frequency Curvelets at the finest scale, instead of having to periodize. In this way, we avoid undersampling of the low-frequency Curvelets at the finest scale which now however corresponds to frequencies lower than when the Curvelets are computed for the full set of frequencies. 
We are going to concentrate on the implementation via wrapping because in this case the original Curvelet transform for the full set of frequencies results in a numerical isometry.
\dontshow{
\begin{figure}[h]
\begin{center}
\includegraphics[width=0.3\textwidth]{lowpass_corona}
\caption{Tiling of the frequency space induced by low pass curvelets. The shadowed outmost Cartesian annulus corresponds to frequencies discarded while computing the low-pass Curvelet transform. In fact, the frequencies at the lower end of the grayed out window are used while computing the Curvelets at the finest scale to avoid undersampling.}
\label{fig:lowpass_corona}
\end{center}
\end{figure}
}

In what follows,  we are going to denote the standard Curvelet transform for an $n_1 \times n_2$ image with $J$ scales as $\mc C_{J}^{n_1,n_2}$ and the low-frequency Curvelet transform corresponding to frequencies up to $n^{LF}_1/2 \times n^{LF}_2/2$ (we will use such simplified notation to denote both the negative and positive frequencies), $n^{LF}_l < n_l,\,l=1,2$ for the same image with $\mc C_{J,n^{LF}_1, n^{LF}_2}^{n_1,n_2}$. 

Let $\tilde U_{J,\ell}^{LF} = [-N_1^{LF}/2, N_1^{LF}/2-1] \times [-N_2^{LF}/2, N_2^{LF}/2-1] $ with $N_l^{LF}=2\lfloor 2m^{LF}_l \rfloor +1,\, m^{LF}_l =  n_l^{LF}/3,\, l=1,2$ denote the frequency domain support of the finest scale low-frequency Curvelet. In general we have two cases: 
\begin{itemize}
\item[i)] {\bf $N_l^{LF}/2 \leq n_l/2,\, l=1,2$:} 
the frequency support of the finest scale low-frequency Curvelet $\tilde U_{J,\ell}^{LF}$ is contained in the fundamental cell $[-n_1/2, n_1/2-1] \times [-n_2/2, n_2/2-1]$; 
\item[ii)] $N_l^{LF}/2 > n_l/2,\, l=1,2$: the finest scale low-frequency Curvelet support $\tilde U^{LF}_{J,\ell}$ extends to frequencies outside of the fundamental cell. 
\end{itemize}

In the first case, to compute the standard Curvelets corresponding to image with maximal frequency $n_1^{LF}/2 \times n_2^{LF}/2$, $\mc C_{J}^{n_1^{LF},n_2^{LF}}$, 
one would compute the Fourier transform of $g^t$, restrict it to frequencies up to $n_1^{LF}/2 \times n_2^{LF}/2$ and subsequently periodize it while weighting with the $C^{\infty}$ window \eqref{eq:cinfw}. Instead, to compute the low-frequency Curvelets, $\mc C_{J,n_1^{LF},n_2^{LF}}^{n_1,n_2}$,
we apply to the full Fourier transform of $g^t$ (with frequencies  up to $n_1/2 \times n_2/2$) a rectangular low-pass window with a cut-off frequency $N_l^{LF}/2,\, l=1,2$ being the highest frequency in the support $\tilde U_{J,\ell}^{LF}$. Then we proceed as for standard Curvelets 
$\mc C_{J}^{n_1^{LF},n_2^{LF}}$ to compute the finest scale coefficients. 

In the second case, to compute the low-frequency Curvelets, $\mc C_{J,n_1^{LF},n_2^{LF}}^{n_1,n_2}$, we still need to use priodization to fill up the support $\tilde U_{J,\ell}^{LF}$. However, we only periodically extend the Fourier transform of $g^t$ beyond the range of available frequencies, here $n_1/2 \times n_2/2$, rather then $n_1^{LF}/2 \times n_2^{LF}/2$ as would be done using standard Curvelets, $\mc C_{J}^{n_1^{LF},n_2^{LF}}$. This results in a $C^{\infty}$ window with shorter (and steeper) flanks as we 
only need to fill in the frequencies in the range $[-N_1^{LF}/2, -n_1/2-1] \cup [n_1/2, N_1^{LF}/2-1] \times [-N_2^{LF}/2, -n_2/2-1] \cup [n_2/2, N_2^{LF}/2-1]$. The corresponding window is a tensor product of the following one dimensional $\mc C^\infty$ windows
\begin{equation}\label{eq:cinfw_lf}
w^{LF}[i] = \left\{\begin{array}{ll}
a(x_i) h_{\infty}(1-x_i),  & i = 1,\dots, N_1^{LF} - n_1,\\
1, &  i = N_1^{LF} - n_1+1,\dots, n_1, \\
a(x_i) h_{\infty}(x_{i-n_1}), & i =     n_1+1, \dots, N^{LF},
\end{array}\right. 
\end{equation}
where $x_i = (i-1)/(N_1^{LF} - n_1-1), \; i =1, \dots, N_1^{LF} - n_1$, and $a$ and $h_\infty$ 
are as before. 
 
\begin{figure}
\centering
\subfloat[][$\psi \in \mc C_3^{256,256}$]{
\includegraphics[width=0.24\textwidth,trim={0 0 0 2.5em},clip]{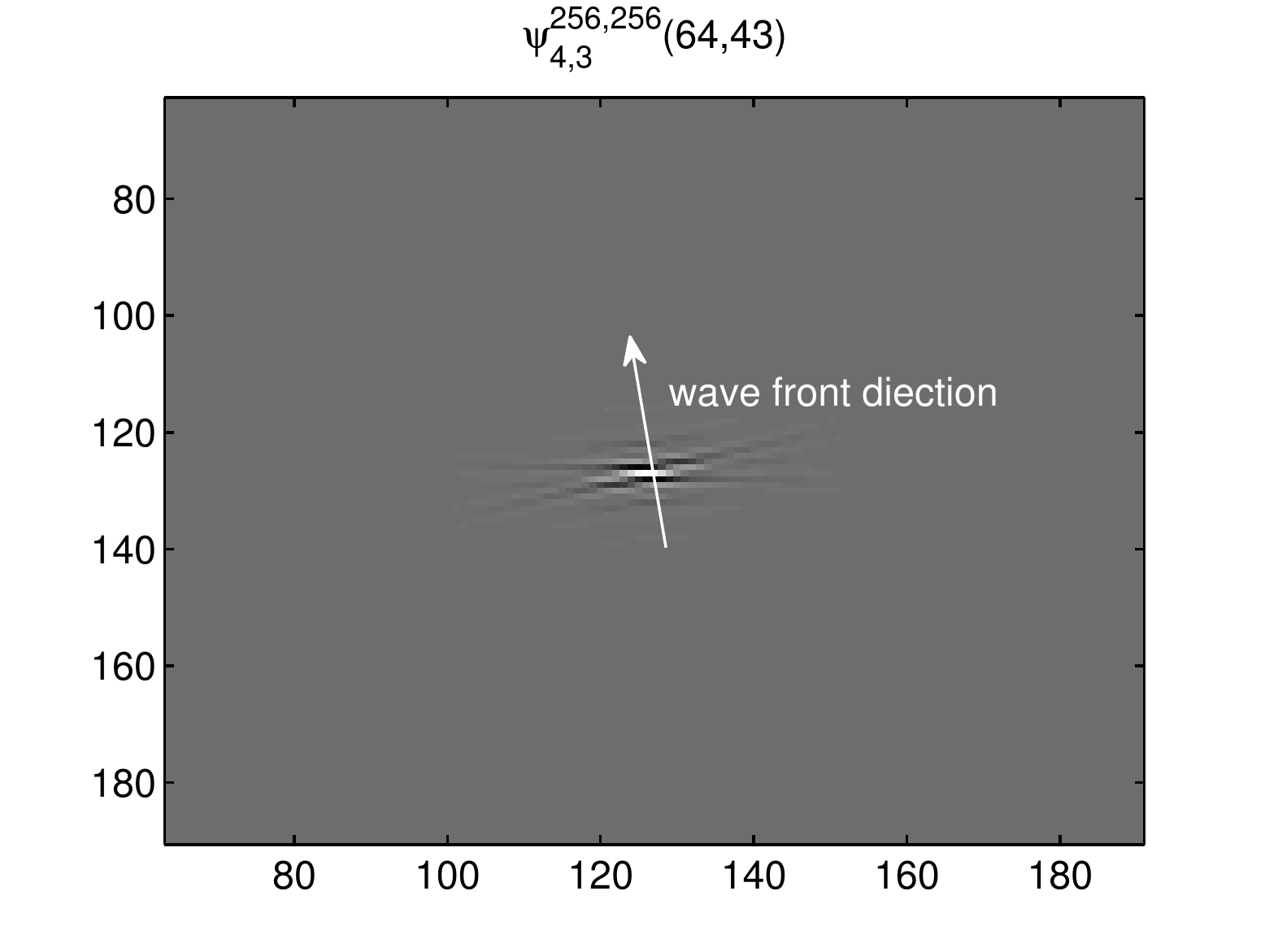}}
\subfloat[][$\psi \in \mc C_{3,192,192}^{256,256}$]{
\includegraphics[width=0.24\textwidth,trim={0 0 0 2.5em},clip]{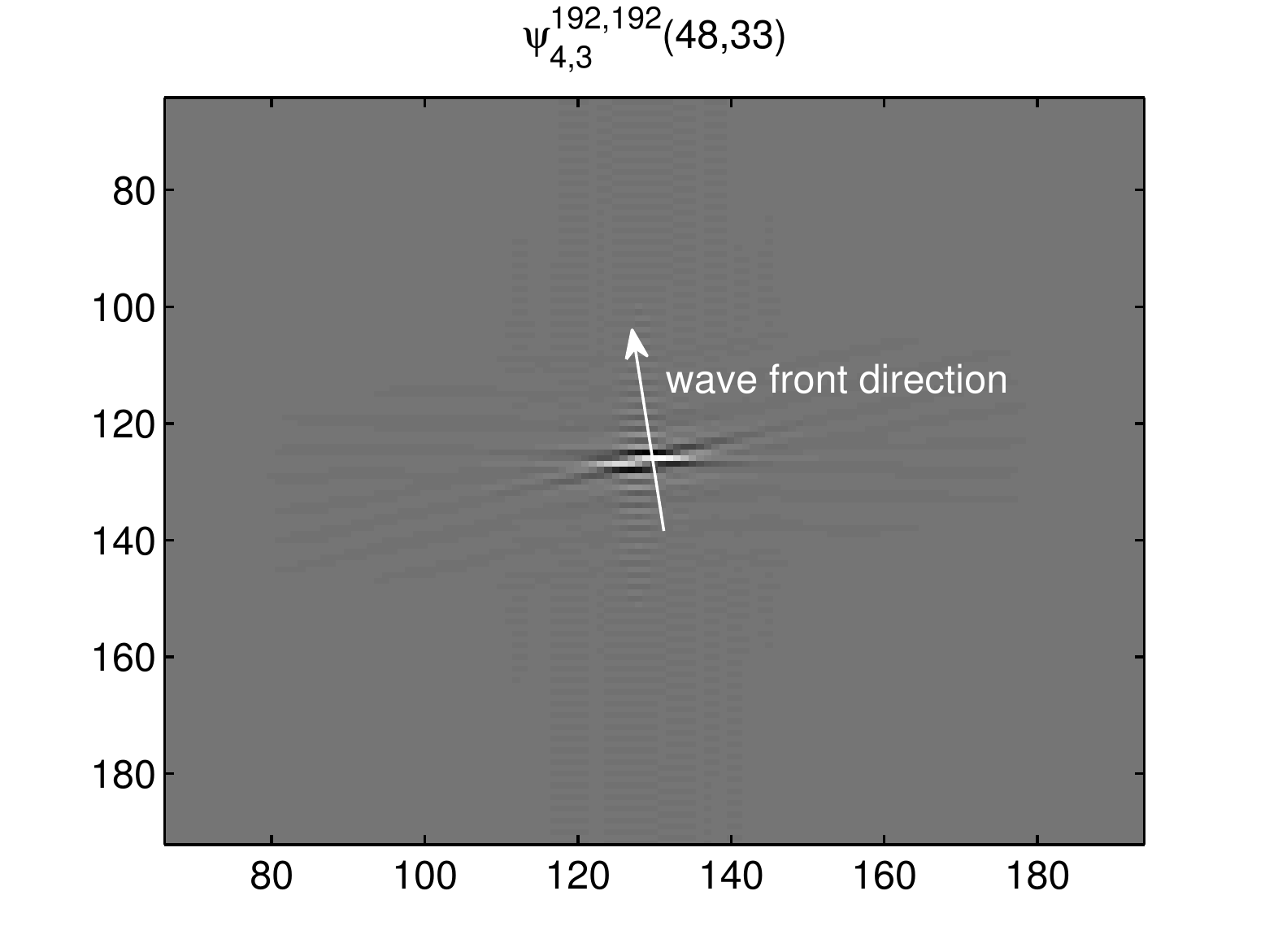}}\\
\subfloat[][$\log|\psi|, \; \psi \in \mc C_3^{256,256}$]{
\includegraphics[width=0.24\textwidth,trim={0 0 0 3em},clip]{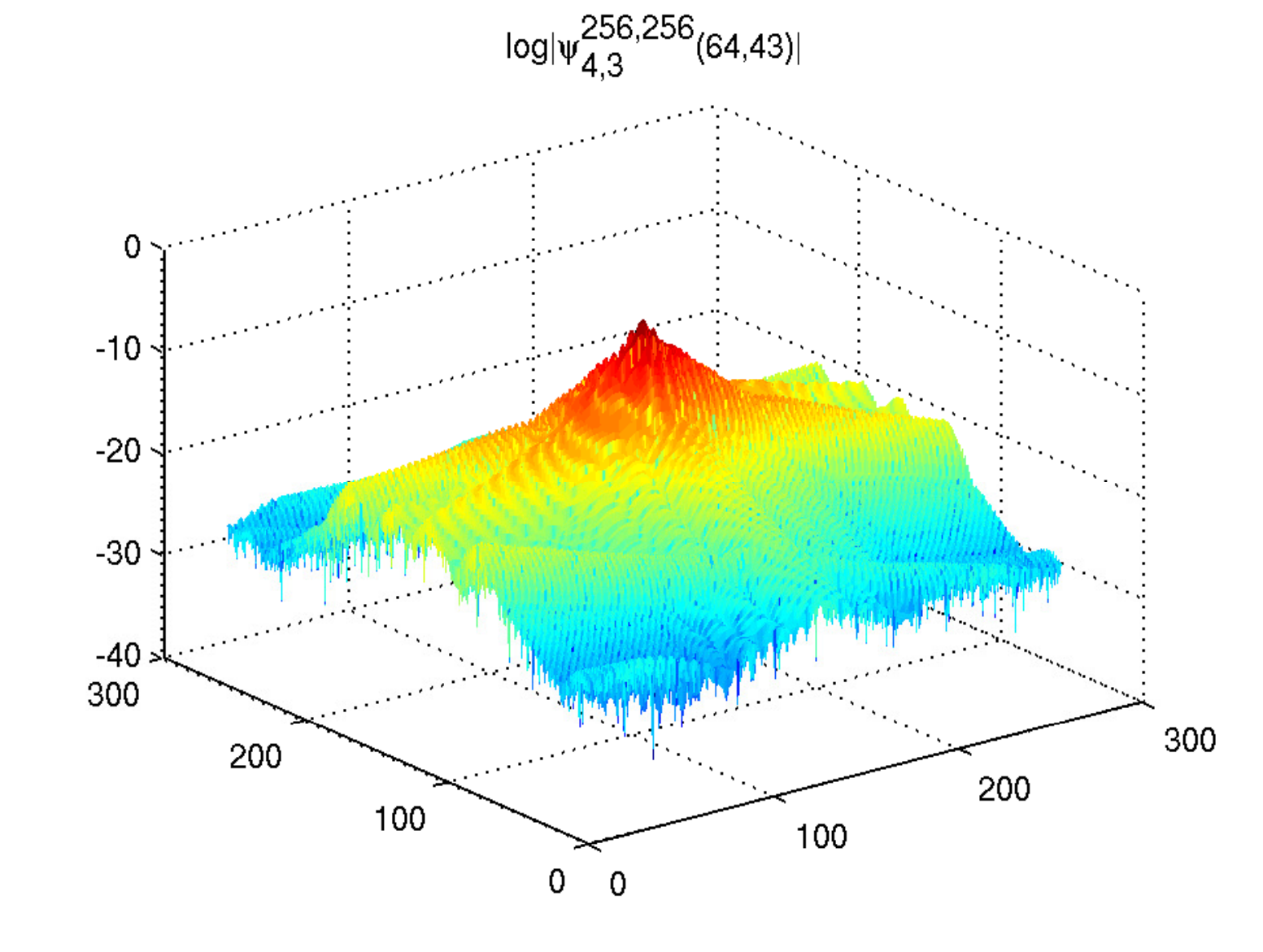}}
\subfloat[][$\log|\psi|, \; \psi \in \mc C_{3,192,192}^{256,256}$]{
\includegraphics[width=0.24\textwidth,trim={0 0 0 3em},clip]{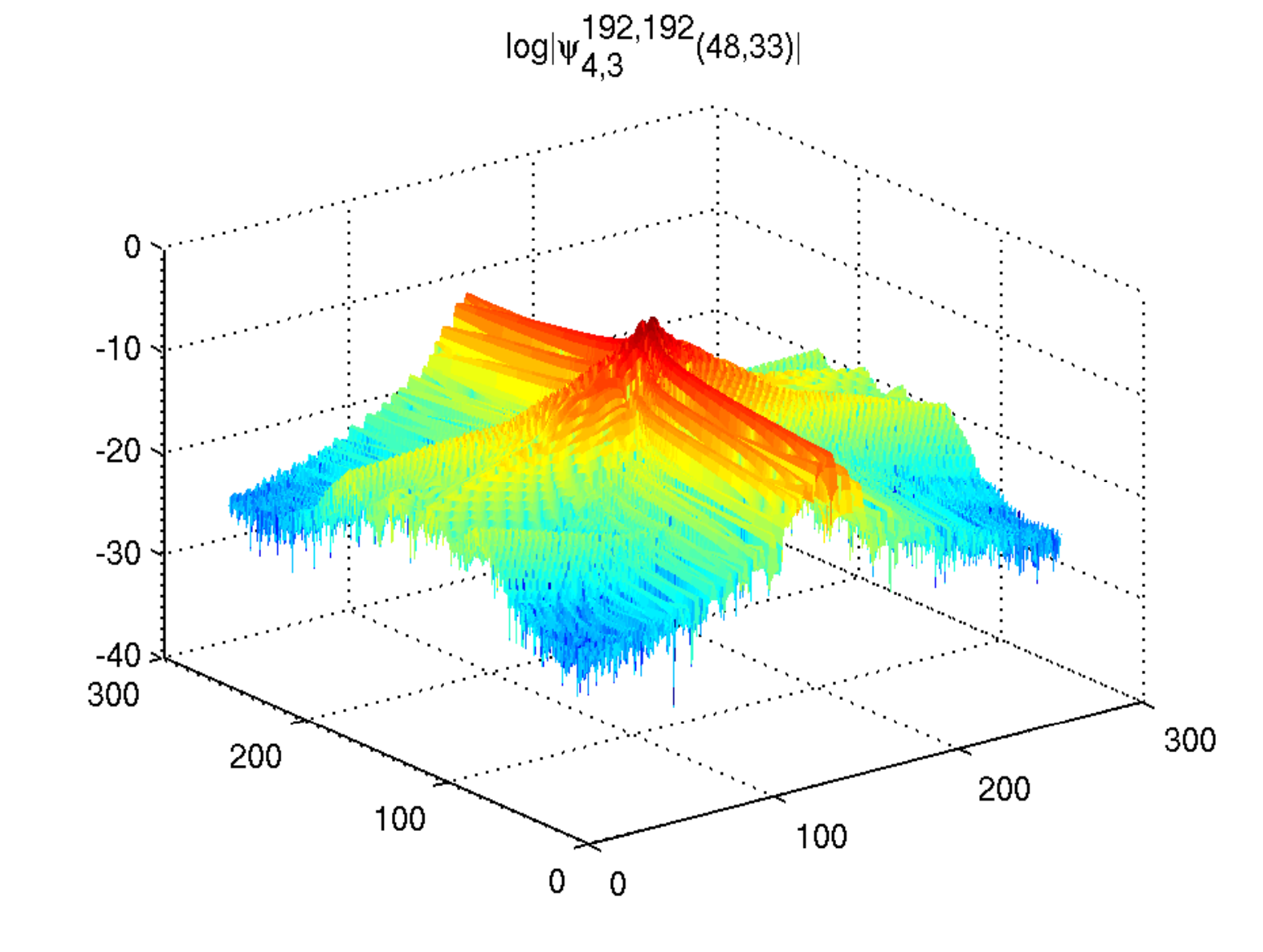}}
\caption{Zoomed in image of (a) standard Curvelet $\mc C_3^{256,256}$,  (b) low-frequency Curvelet $\mc C_{3,192,192}^{256,256}$,  at the finest level 3. 
(c), (d) Surface plot of log amplitude of (a), (b).}
\label{fig:logCCLF}
\end{figure}

Subsequently, all the lower scales of low-frequency Curvelets can be computed exactly as for the standard Curvelets $\mc C_{J}^{n_1^{LF},n_2^{LF}}$. The construction of low-frequency Curvelets above results in a numerical isometry on the restriction to frequency range up to $\min\{N_1^{LF}/2,n_1/2\} \times \min\{N_2^{LF}/2,n_2/2\}$. Thus low-frequency Curvelet transform is a different transformation to the original Curvelet transform which however for a choice of $n_l^{LF}/2$ such that $N_l^{LF}/2 \geq n_l/2,\, l=1,2$, is also an isometry on $[-n_1/2, n_1/2-1] \times [-n_2/2, n_2/2-1]$ as the standard Curvelets $\mc C_{J}^{n_1,n_2}$. Figures \ref{fig:logCCLF}(a),(b) show the finest scale Curvelet and low-frequency Curvelet, respectively. The plots of the logarithm of the amplitude in Figures \ref{fig:logCCLF}(c),(d) illustrate that the decay of the low-frequency Curvelets along the wavefront direction is slower than that of the standard Curvelets, which is due to the frequencies higher then $n^{LF}_l/2, \, l=1,2$ used in computation of Curvelets at the finest scale.

Recapitulating, there are two major benefits of such low-frequency Curvelet transform. First, is the super linear reduction in computation cost of the transform, which is particularly beneficial for solution of the CS recovery problem \eqref{eq:l1rec} which involves repeated application of $\Phi \Psi\tr$ and its adjoint $\Psi\Phi\tr$. Second, is the ability to effectively represent realistic PAT data, in which due to measurement process the high frequencies are damped as in the model derived in Section \ref{sec:ms:bli}. This results in amplitudes of low-frequency Curvelet coefficients being higher and exhibiting a quicker decay than those of the standard Curvelet coefficients (see Figure \ref{fig:curv} and the accompanying discussion in Section \ref{sec:recon:clock}), and consequently in higher robustness to noise and imperfect compressibility. 



\section{Compressed sensing of optical ultrasound detector}\label{sec:cspat}
\subsection{Single-pixel optical camera}\label{sec:cspat:spoc}
In a series of publications \cite{Huynh:2014pi, Huynh:2015rtusm, Huynh:2016spoc} we introduced a single-pixel optical camera (SPOC) for ultrasonic and photoacoustic imaging. With the SPOC, instead of recording the pressure on the detector point-by-point, the entire active area of the optical ultrasound sensor is illuminated, and using a digital micro-mirror device (DMD) a pattern $\phi_j$ is applied to the wide field light reflected from the sensor resulting in a compressed measurement  
\begin{equation}
b_j = \phi_j\tr g(\bx_{\mc S}, t_i), \quad t_i \in (0,T).
\end{equation}
Figure \ref{fig:spoc} shows a sketch of the operational principle of SPOC while for technical details of the system we refer to \cite{Huynh:2016spoc}.    

\subsection{Sparse representation of the sensor data}
In PAT the entire time series for one point or pattern is acquired with one excitation, $\phi_j\tr g(\bx, t_i)$, $\bx \in \mathcal S, t_i \in (0, T)$. This has the consequence that we acquire the same number of compressed/point measurements of the wave field at each time step $t_i$. Furthermore, due to the rate at which the DMD can change patterns, at least at present, we are limited to use of only one pattern throughout the acoustic propagation. 

However, as the wave propagates, the complexity of the sensor data varies with time. After a wave front reaches the detector, its cross-section with the detector plane expands and more wave fronts, corresponding to features farther away, reach the detector. As a result, in general the complexity of the wave field at the sensor over time first increases and at some point it starts to decrease again corresponding to the tail of the wave field. If the complexity is reflected by the sparsity, the error of the best $k$-term approximation for the sensor data varies throughout the time series and with it the recovery error bound in the robust recovery Theorem \ref{thm:l1rec}. 
Furthermore, due the sound intensity, which is $\propto$ $p(\bx,t)^2$, obeying the inverse square law, 
$p(\bx,t)$ decays as inverse distance of $\bx$ to the source. This means that as with increasing $t$, $p(\bx_{\mathcal S}, t)$ encodes information about $p_0$ further away from the sensor, the recorded pressure amplitudes decrease linearly for the same initial pressure amplitude value, ultimately resulting in a lower point wise signal to noise ratio for longer propagation times.

\dontshow{[
Smoothing of the wavefront vs construction of sparse functions - how to do this optimally. This is both simulation and real data. Curvelet: How to choose the number of levels for PA signal, wavelet on the finest (directionality of the wave field at the finest level). The effect of the smoothing on scales (quicker decay of the coefficients at all scales, but server damping of those at the highest scales, window used for smoothing eliminates the corner curvelets)

Discuss the compression, and compression at each scale. This will tell is, which features get lost if the applied compression is to high. For the reconstruction this means, that if the algorithms reconstructs the largest coefficients, those features will be reconstructed while the other lost (how does the noise) - is this true? Do the recovery results guarantee that the largest coefficients will be recovered ? 

The wave field is differently compressible at different time steps. Thus, we will automatically under sample some steps while oversample others. 
If we can guarantee that e.g. the largest Curvelet coefficients are recovered, can we perform some kind of analysis what happens during Fourier PAT reconstruction? 
Idea: Assume we recover the lower frequencies in the corona (Curvelets partitioning Fourier space). The Curvelts of images are trees. So I can guess which locations (approx) would there be new Curvelets. Also 

The planar crossection of the wave field at given time does not result in C2 edges. Should we discuss this?  
]}

\subsection{Sensing patterns}
In our experiments we used scrambled Hadamard patterns
\begin{equation}\label{eq:sHad}
H_j^s = P_r H_j P_c,
\end{equation}
where $H_j$ is the $2^j\times 2^j$ Hadamard matrix and $P_c, P_r \in \{0,1\}^{2^j\times 2^j}$ denote permutation matrices for columns and rows, respectively. For compressed sensing we only select first $m \ll n=2^j$ rows of $H_j^s$. The application of $H_j^s$ to a vector $H_j^s v = P_r H_j P_c v$ amounts to application of the permutation matrix $P_c$ to $v$, performing Hadamard transform on the permuted vector $P_c v$, and subsequently permuting the rows. Thus the scrambled Hadamard transform can be computed at essentially the same cost as the fast Hadamard transform, while scrambled Hadamard matrices have recovery properties similar to those of random Bernoulli matrices see e.g.~\cite{Foucart:2010}. 

In practice using the DMD, it was only possible to apply binary patterns. In order to make use of properties of Hadamard transform such as othogonality and self inversion, the experimental Hadamard matrix $H_j^{(0,1)}$ needs to be transformed into the Hadamard matrix $H_j$ using the simple relation
\begin{equation}\label{eq:tBinH}
H_j = \left(2H_j^{(0,1)} -\mathbf{11^T}\right)/\sqrt{2^j},
\end{equation}
and correspondingly the measured data $w = H_j^{(0,1)} g$ into
\begin{equation}\label{eq:tBinData}
H_j g = \frac{1}{\sqrt{2^j}}\left(2H_j^{(0,1)} -\mathbf{11^T}\right) g = \frac{1}{\sqrt{2^j}}\left(2 w - \mathbf{1} w(1)\right).
\end{equation}
Here, we used that $w(1) = \mathbf{1^T} g$ corresponds to the measurement acquired with `all-1s' pattern, which is the first row of $H^{(0,1)}$. It is immediately clear that the same transformation \eqref{eq:tBinH}, \eqref{eq:tBinData} holds for scrambled Hadamard matrices \eqref{eq:sHad} (with 1 in $w(1)$ replaced by the `all-1s' row number after row permutation).      

The light reflected from the DMD is integrated by a photodiode. In order to best utilize the dynamic range of the photo diode, it is necessary to keep the optical power incident on the photodiode in the same range for each pattern. All but the `all-1s' pattern are composed of an equal number of 0 and 1s. Therefore, the `all-1s' pattern was replaced with a vector which first half entries are 0 and the second half 1. As the negative of this vector (0 becomes 1 and vice versa) is exactly the $2^j/2+1$ row of $H_j^{(0,1)}$, the data corresponding to the `all-1s' pattern can be constructed by adding the data from the modified first row pattern and $(2^j/2+1)$th row pattern. Again, this construction is not affected by  scrambling.

\subsection{Recovery of the sensor data}
The sampled PAT data at each time step, $g(\bx_{\mc S}, t)$, can be recovered by solving the optimization problem \eqref{eq:l1rec}. The sensing matrix $\Phi$ in \eqref{eq:l1rec} is set to be the first $m$ rows of the scrambled Hadamard matrix $H^{s}_{\log_2(n)}$, $n = n_1 n_2$, 
\begin{equation}
\Phi = S H^s_{\log_2(n)} = S P_r H_{\log_2(n)} P_c,
\end{equation}
where $S \in \{ 0,1 \}^{m \times n}$
is a binary subsamplig matrix such that $S\tr S$ is a $n \times n$ diagonal matrix with ones at positions corresponding to the chosen and zeros to the ommitted rows, respectively, while
$S S\tr = I^{m\times m}$ is an $m \times m$ identity matrix. Consequently, we have $\Phi\Phi\tr = I^{m\times m}$.

The vector of measurements, $b$, is computed from experimantal measurements using equation \eqref{eq:tBinData}.
The sparsifying transform $\Psi$ is chosen as an orthonormal low-pass or standard Curvelet transform on the $n_1 \times n_2$ optical sensor $\mc S$ and $f^t$ is the (sparse) vector of sought for coefficients, $f^t = \Psi g^t$.

To take advantage of the structure of the problem we solve \eqref{eq:l1rec} using the Split Augmented Lagrangian Shrinkage Algorithm (SALSA) proposed in \cite{Afonso:2010salsa}. SALSA is an ADMM scheme which solves the unconstrained problem
\begin{equation}\label{eq:l1rec_reg}
\min_f \zeta(f) := \frac{1}{2}\| \Phi\Psi\tr f - b \|_2^2 + \tau \| f \|_1.
\end{equation}
A version used in this work is summarized in Algorithm \ref{alg:salsa}. 

\begin{algorithm}[h]
\caption{Split Augmented Lagrangian Shrinkage Algorithm (SALSA), \cite{Afonso:2010salsa}.}\label{alg:salsa}
\begin{algorithmic}[1]
\STATE Choose $\mu >0$, $v_0$ and $d_0$
\STATE $i := 0$
\REPEAT
\STATE 
$f_{i+1} = \argmin_f \| \Phi\Psi\tr f - b \|_2^2 + \mu \| f - v_i - d_i \|_2^2$ 
\STATE 
$v_{i+1} = \argmin_v \tau \| v \|_1 + \mu/2\| f_{i+1} - v - d_i\|_2^2$
\STATE $d_{i+1} = d_i - (f_{i+1} - v_{i+1})$
\STATE $i = i+1$
\UNTIL{$ |\zeta(f_{i+1}) - \zeta(f_{i})|/\zeta(f_{i}) < \varepsilon $}
\end{algorithmic}
\end{algorithm}
The quadratic minimization problem in line 4 leads to the linear system
\begin{equation}
f_{i+1} = (\Psi\Phi^{\rm T}\Phi\Psi\tr + \mu I)\inv (\Psi \Phi\tr b + \mu (v_i + d_i) ),
\end{equation}
where using Sherman Morrison Woodbury formula the inverse can be expressed in terms of application of $\Phi\Psi\tr$ and its adjoint
\begin{equation}
(\Psi\Phi\tr\Phi\Psi\tr + \mu I)\inv = \frac{1}{\mu}\left(I - \frac{1}{\mu+1} \Psi\Phi\tr \Phi\Psi\tr\right).
\end{equation}
The proximal operator in line 5 amounts to point wise soft thresholding 
\begin{equation}
\mc T_{\tau/\mu} (x) = \mbox{sign}(x)(|x| - \tau/\mu)_+. 
\end{equation}

\section{PAT reconstruction with the recovered sensor data}\label{sec:recon}
We present reconstruction results for three photoacoustic image reconstruction problems, reconstruction from: simulated pattern data, synthesized pattern data from a point-by-point scanner data, and data acquired with the SPOC (Section \ref{sec:cspat:spoc}). After the PAT data has been recovered using the proposed method for acoustic field reconstruction, the PAT images are reconstructed using time reversal via first order method in $\bk$-Wave Toolbox.   

\subsection{Simulated data: clock phantom}\label{sec:recon:clock}
\begin{figure}[ht]
\centering
\begin{minipage}{0.29\textwidth}
\centering
\includegraphics[width=1\textwidth,trim={0 4em 0 8em},clip]{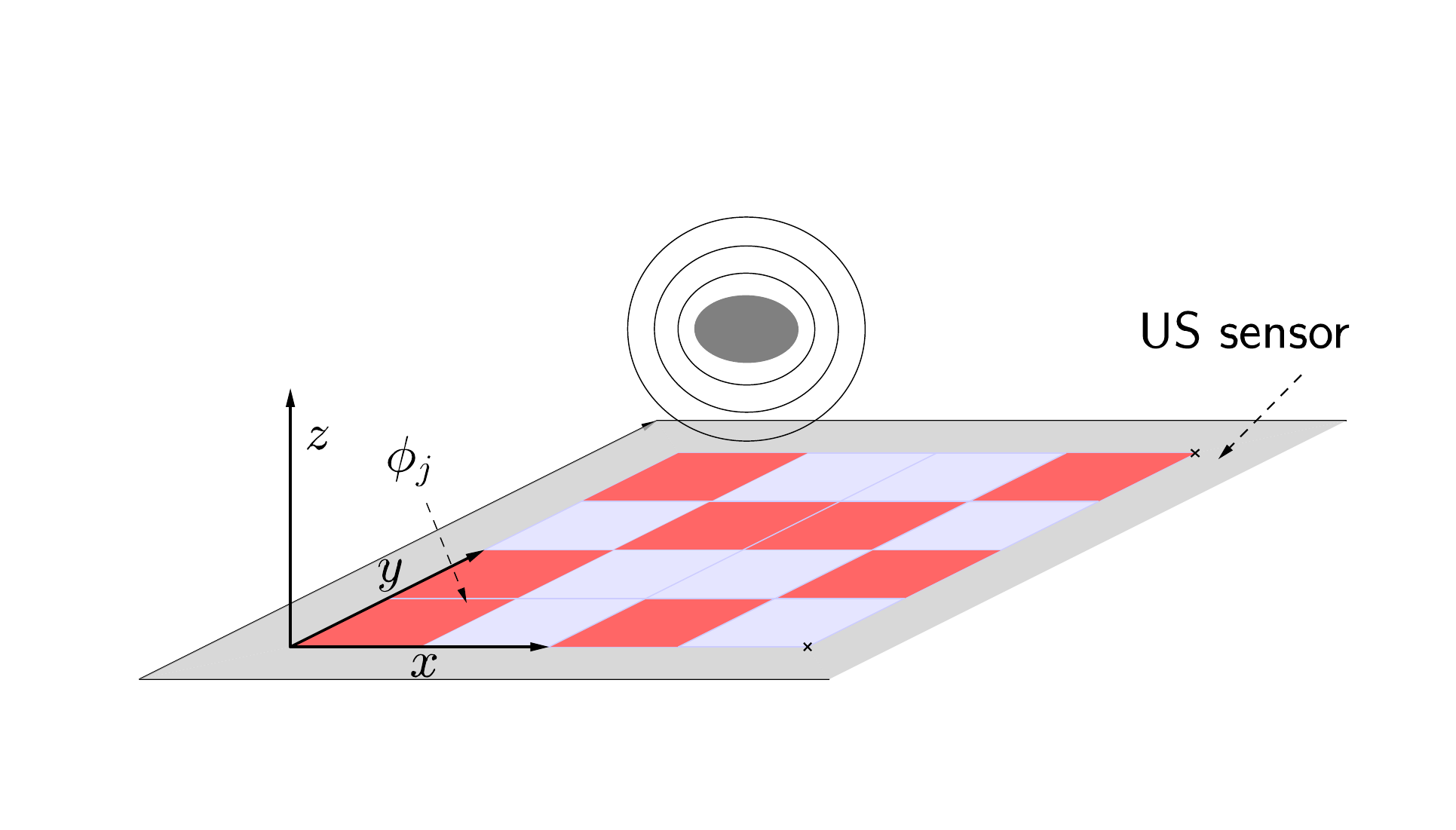}
\caption{Single-pixel optical camera.}
\label{fig:spoc}
\end{minipage}\hfill
\begin{minipage}{0.19\textwidth}
\centering
\includegraphics[width=\textwidth,trim={0 10em 0 15em},clip]{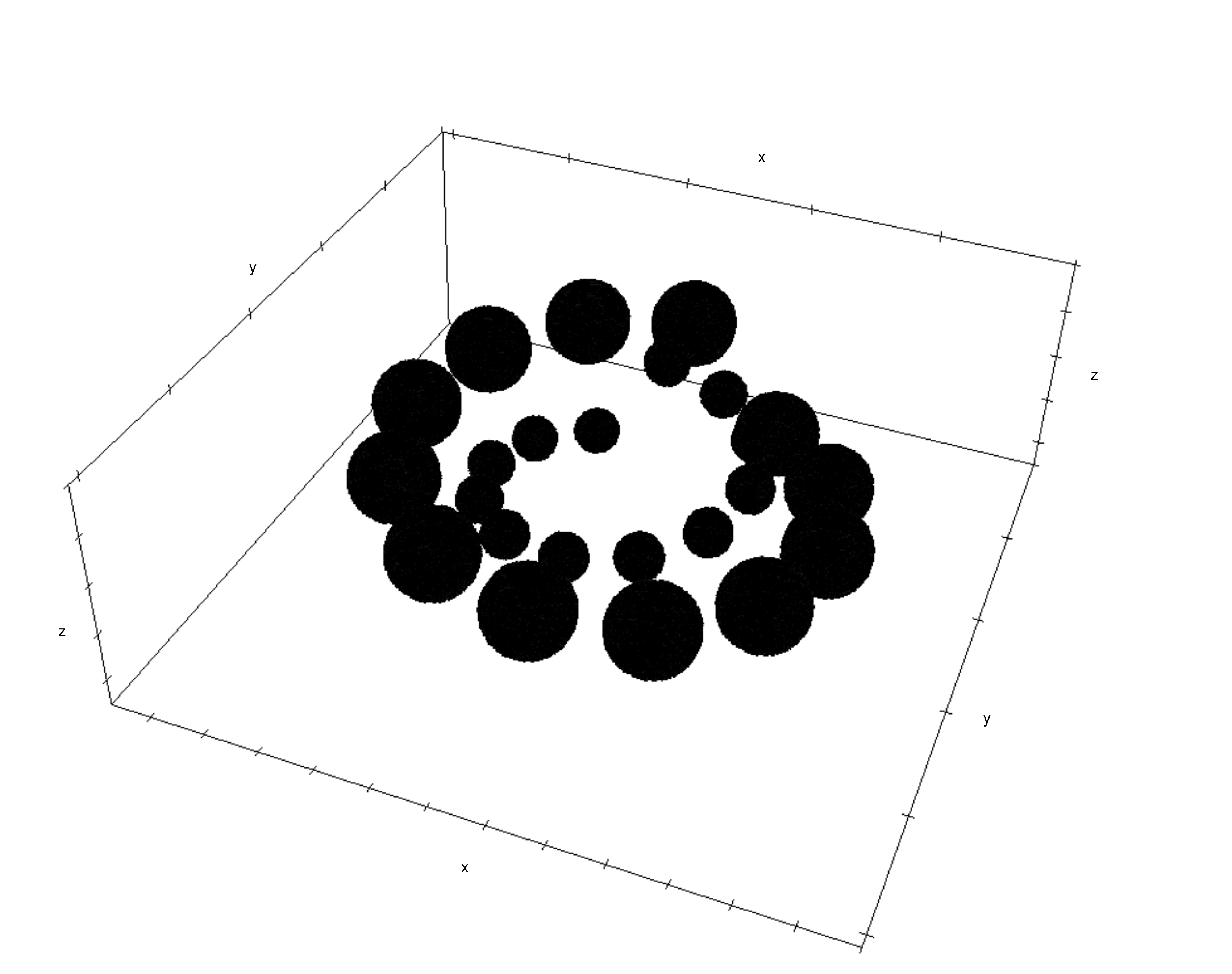}
\caption{Clock phantom.}
\label{fig:clock_phantom}
\end{minipage}
\end{figure}

In our first example we simulate the PAT data for initial pressure distribution $p_0$ depicted in Figure \ref{fig:clock_phantom}.
The ``clock phantom'' is a collection of balls which give rise to spherical wave forms expanding uniformly in all directions, which is a difficult case for directional basis like Curvelets.  

We consider a volume of $ 256\times 256\times 92$ voxels of size $0.05 \si{\milli\meter}^3$ and the sensor of matching resolution, $256\times 256$, placed at $z=0$. We assume homogeneous ambient speed of sound and density of $1500 \si{\meter/\second}$ and $1000 \si{\kg/\meter^3}$, respectively. The pressure at the sensor is sampled every $20 \si{\nano\second}$. Both, the forward and inverse PAT problems are solved with first order method from \textbf{k}-Wave Toolbox. Before the PAT time series is simulated, $p_0$ is smoothed with the 3D Blackman window. 

We attempt to recover the PAT data from 18\% of noiseless compressed measurements obtained with binary scrambled Hadamard patterns. We solve the recovery problem in each time step using SALSA with the data dependent regularization parameter $\tau =0.01 \max(|\Psi\Phi\tr b^t|)$ and $\mu = 5 \max(|\Psi\Phi\tr b^t|)/\|b^t\|_2$. We stop the algorithm if the relative change in the objective function $\zeta(f)$ \eqref{eq:l1rec_reg} drops below $5\cdot 10^{-4}$ or after $100$ iterations. 

We start by examining the utility of the Curvelet transform as a sparsifying transform for PAT data recovery problem. In particular we compare the standard Curvelet transform $\mc C^{256,256}_{3}$ with the low-frequency Curvelet transform $\mc C^{256,256}_{3,192,192}$, the lowest resolution transform which is still an isometry on the full frequency range (up to $256/2\times 256/2$). 

\begin{figure}[ht]
\centering
\subfloat[][$\log |\Psi g(\bx_{\mc S}, t_{100})|$]{
\includegraphics[width=0.2\textwidth,trim={0em 0em 0em 2.5em},clip]{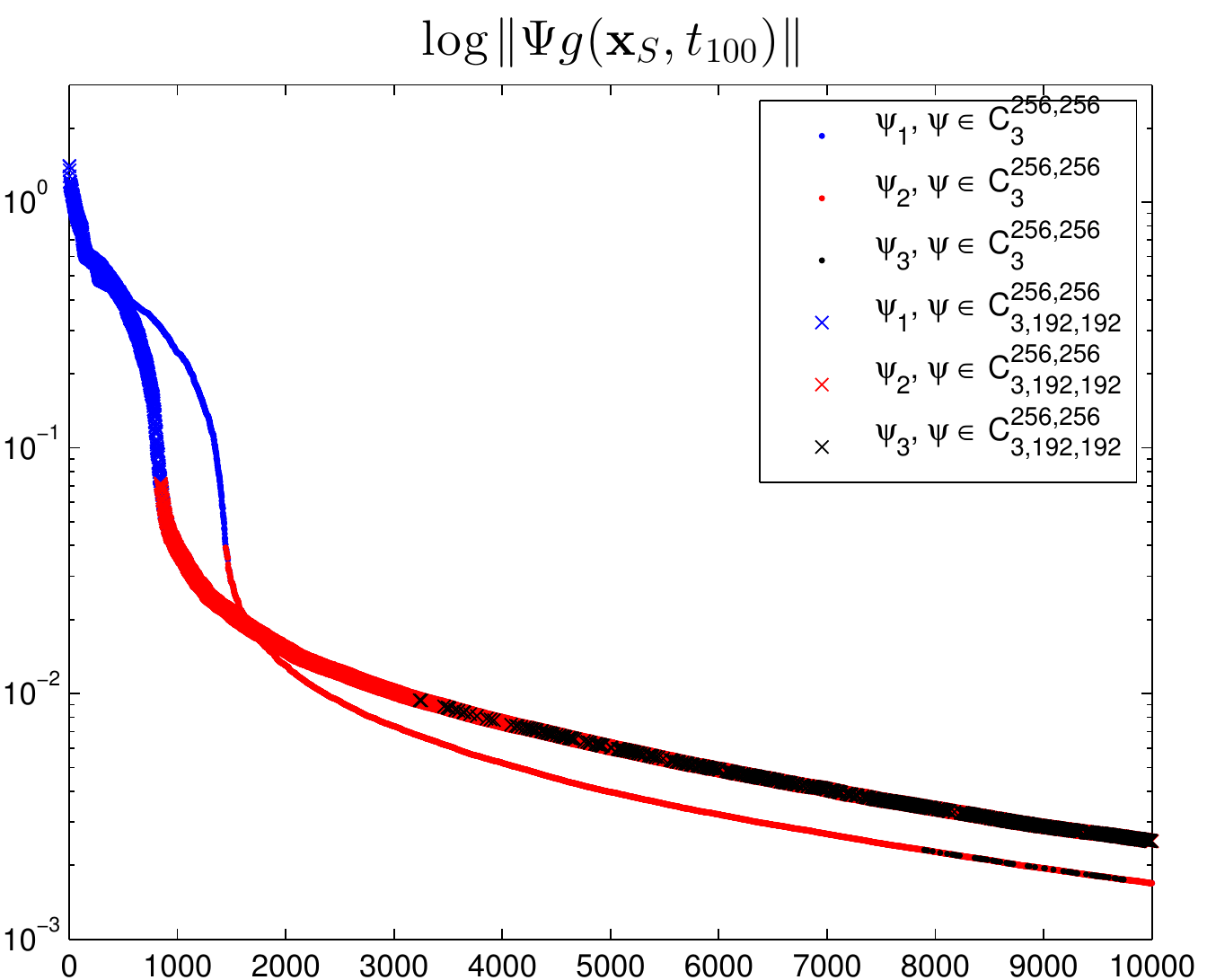}}\hspace{0.03\textwidth}
\subfloat[][$\log |\Psi g(\bx_{\mc S}, t_{230})|$]{
\includegraphics[width=0.2\textwidth,trim={0em 0em 0em 2.5em},clip]{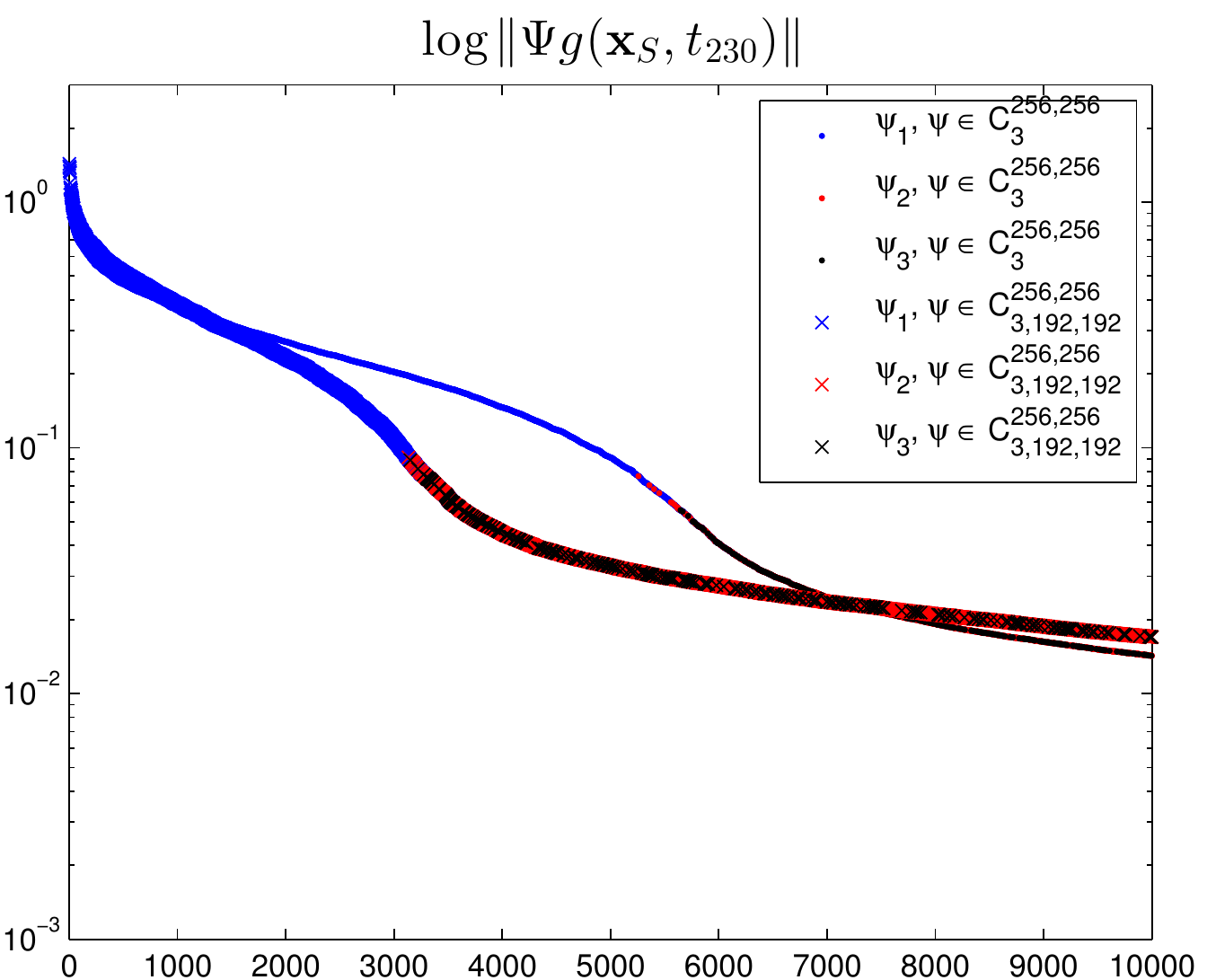}}\\
\subfloat[][Error of $\hat g(\bx_{\mc S}, t_{100})$]{
\includegraphics[width=0.2\textwidth,trim={0em 0em 0em 2.5em},clip]{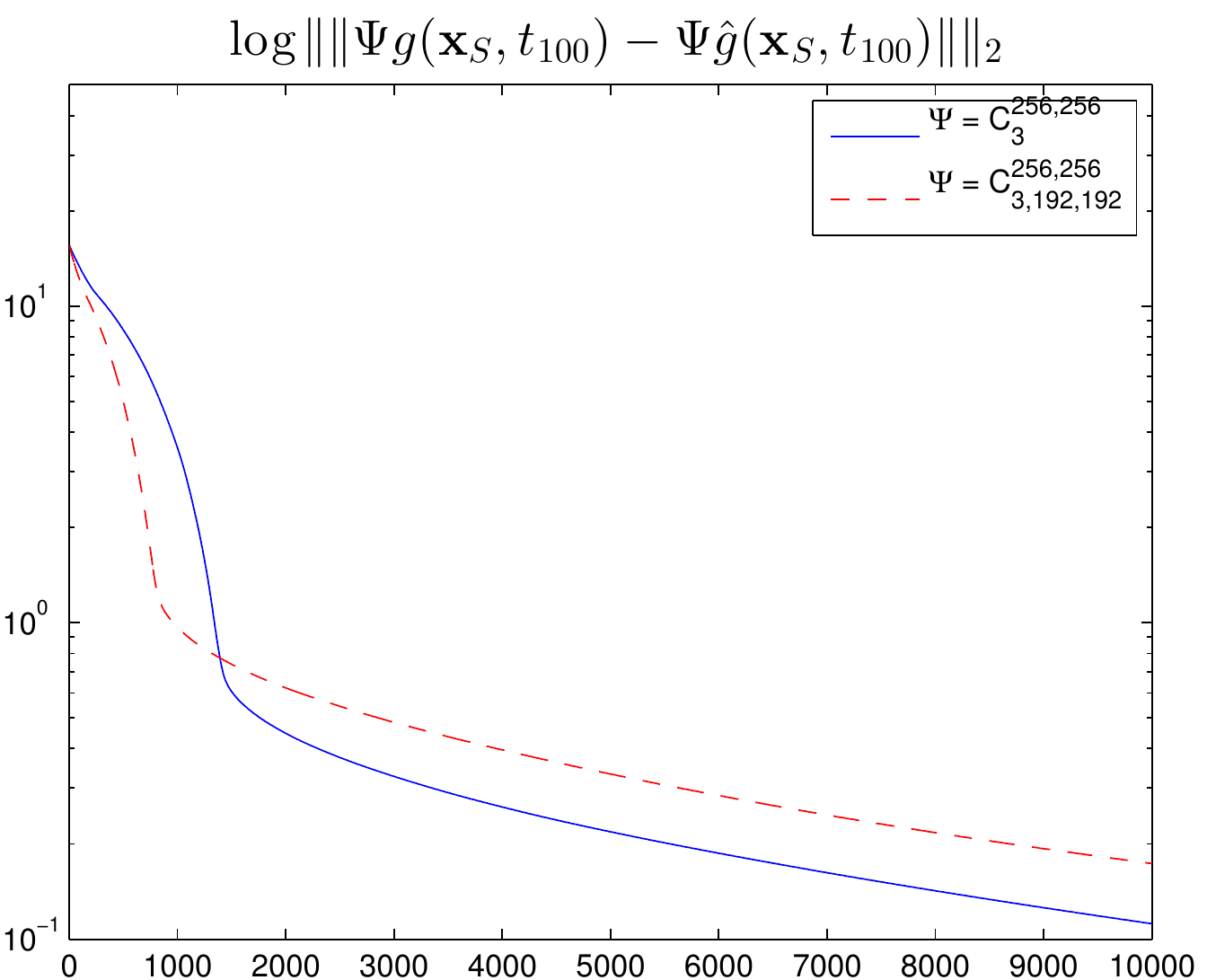}}\hspace{0.03\textwidth}
\subfloat[][Error of $\hat g(\bx_{\mc S}, t_{230})$]{
\includegraphics[width=0.2\textwidth,trim={0em 0em 0em 2.5em},clip]{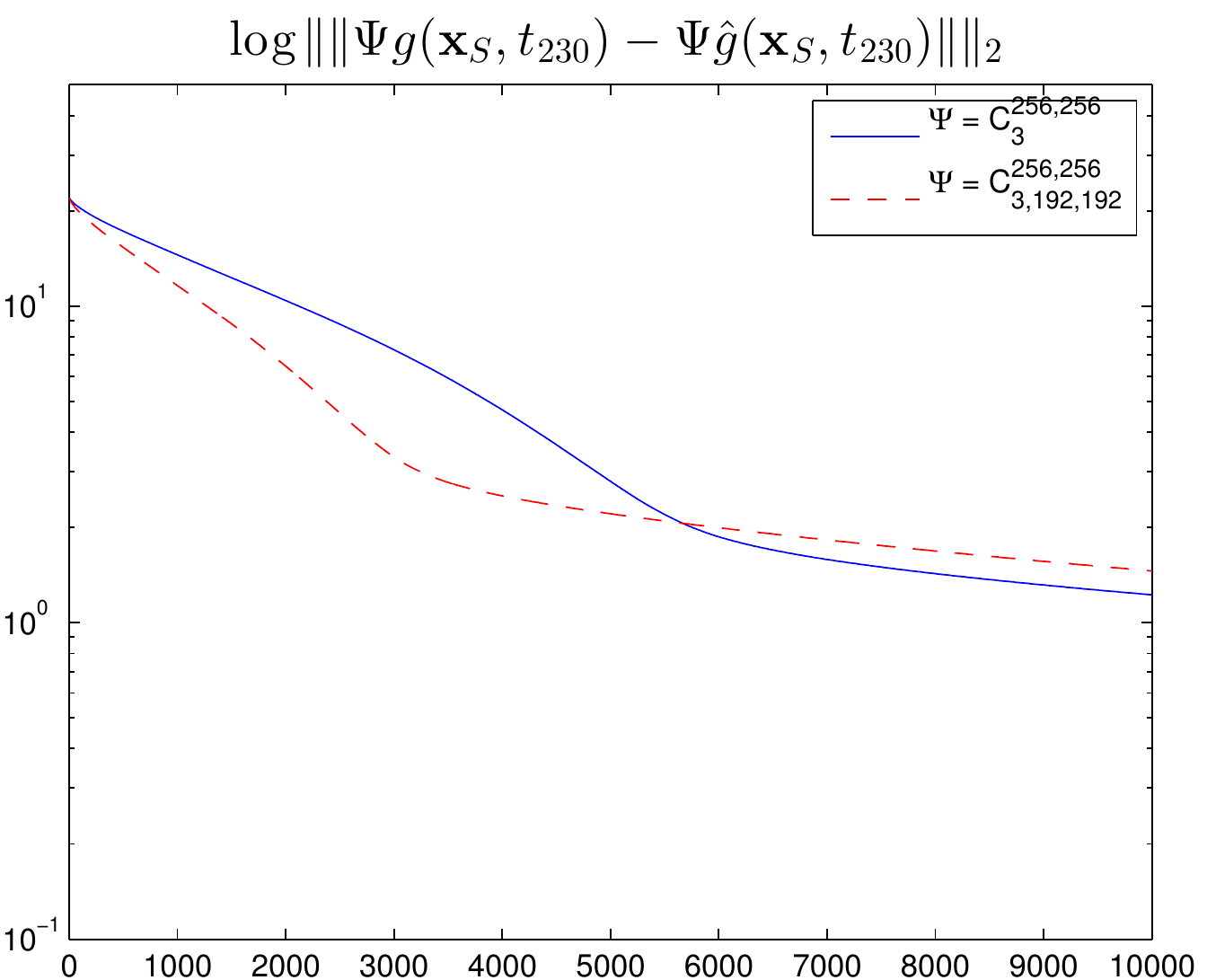}}
\caption{Clock phantom. The decay of log amplitude of Curvelet coefficients of $g(\bx_{\mc S},t_i)$ at time steps  (a) $t_{100}$, (b) $t_{230}$; the colors correspond to coefficients at different scales. (c), (d) The corresponding compression error in log scale.}
\label{fig:curv}
\end{figure}

\begin{figure}[ht]
\subfloat[][$g(\bx_{\mc S},t_{100})$]{ 
\includegraphics[width=0.160\textwidth,trim={0 0 1em 2.5em},clip]{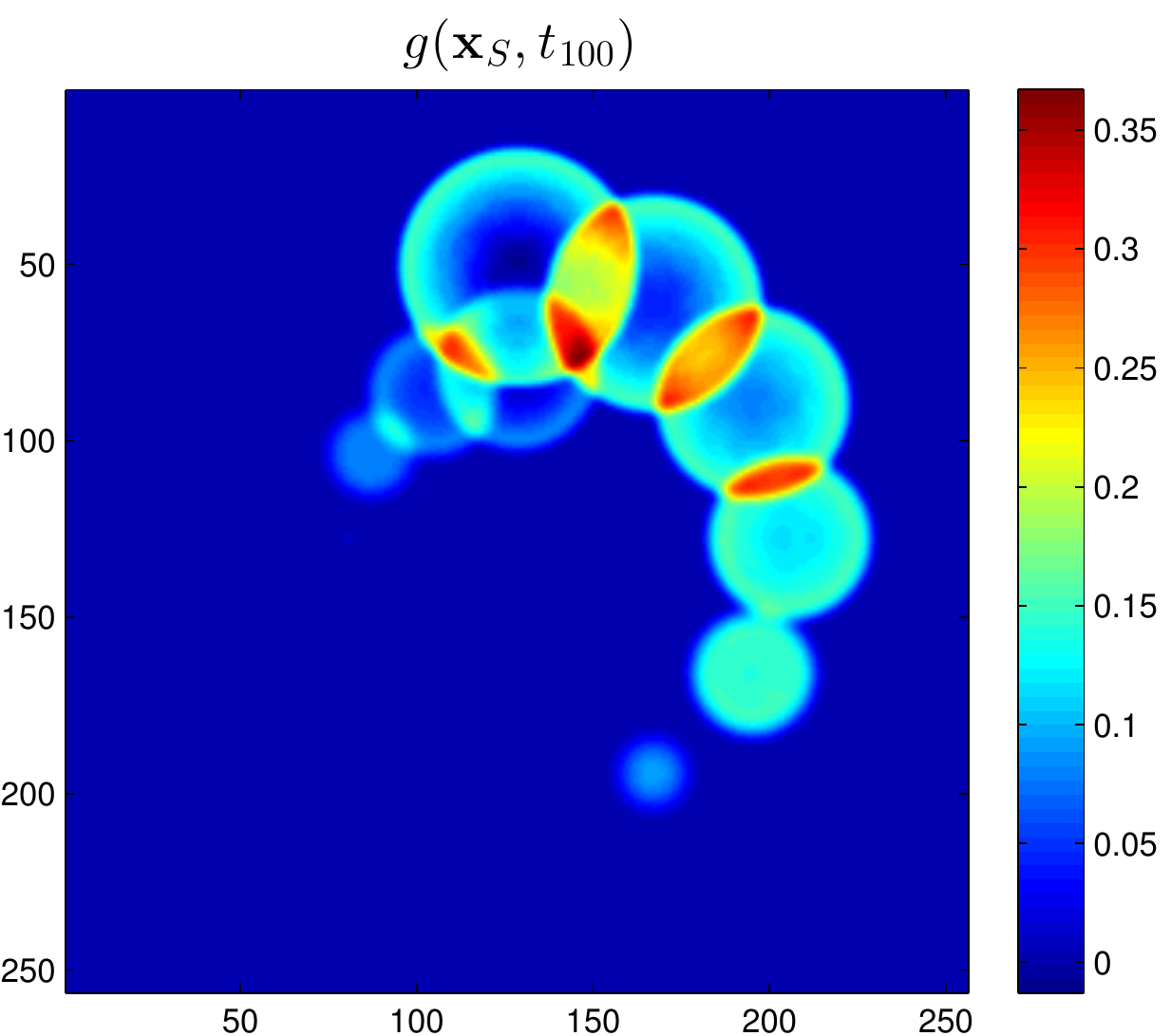}}
\subfloat[][$\hat g(\bx_{\mc S},t_{100})$]{ 
\includegraphics[width=0.137\textwidth,trim={0 0 0em 2.5em},clip]{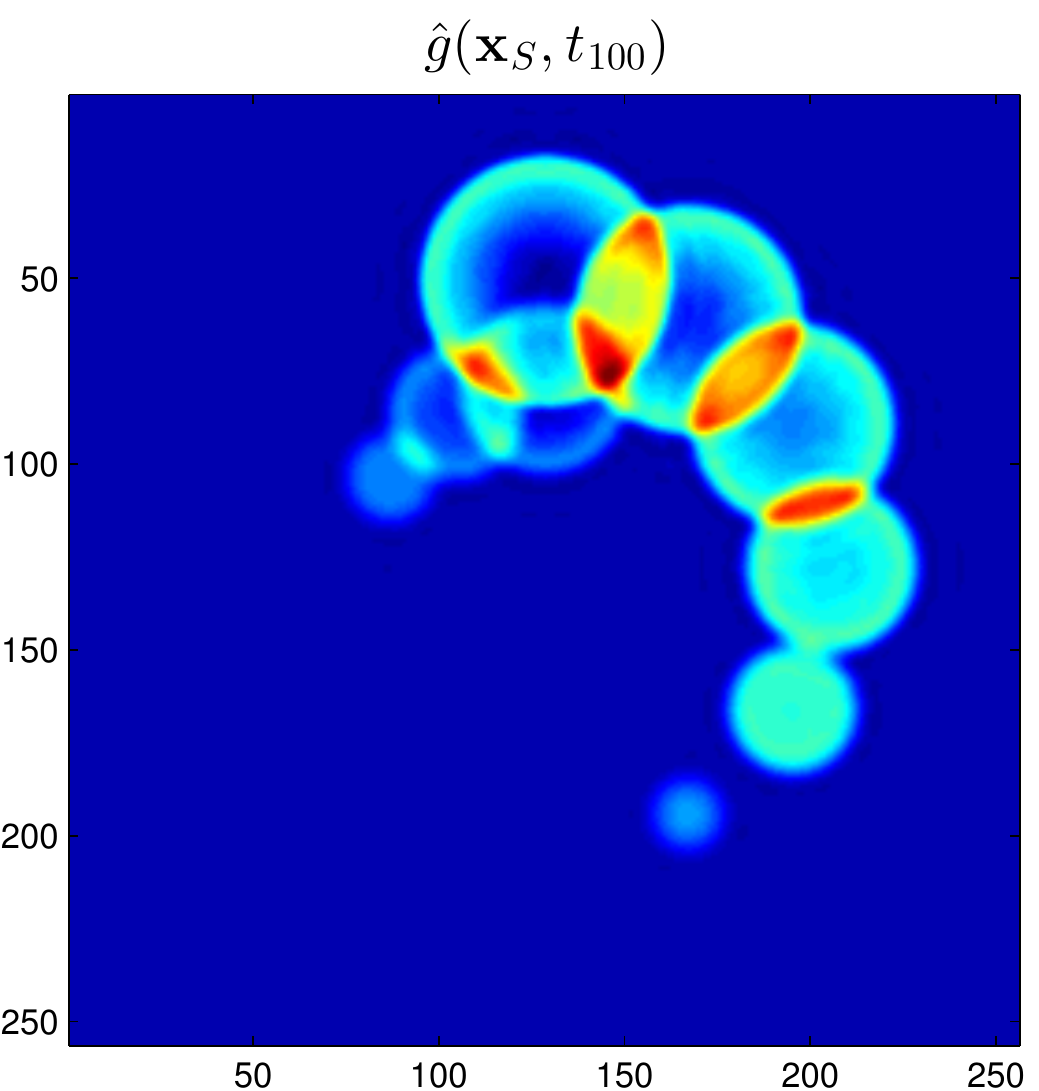}}
\hspace{0.008\textwidth}
\subfloat[][$\tilde g(\bx_{\mc S},t_{100})$]{ 
\includegraphics[width=0.137\textwidth,trim={0 0 0em 2.5em},clip]{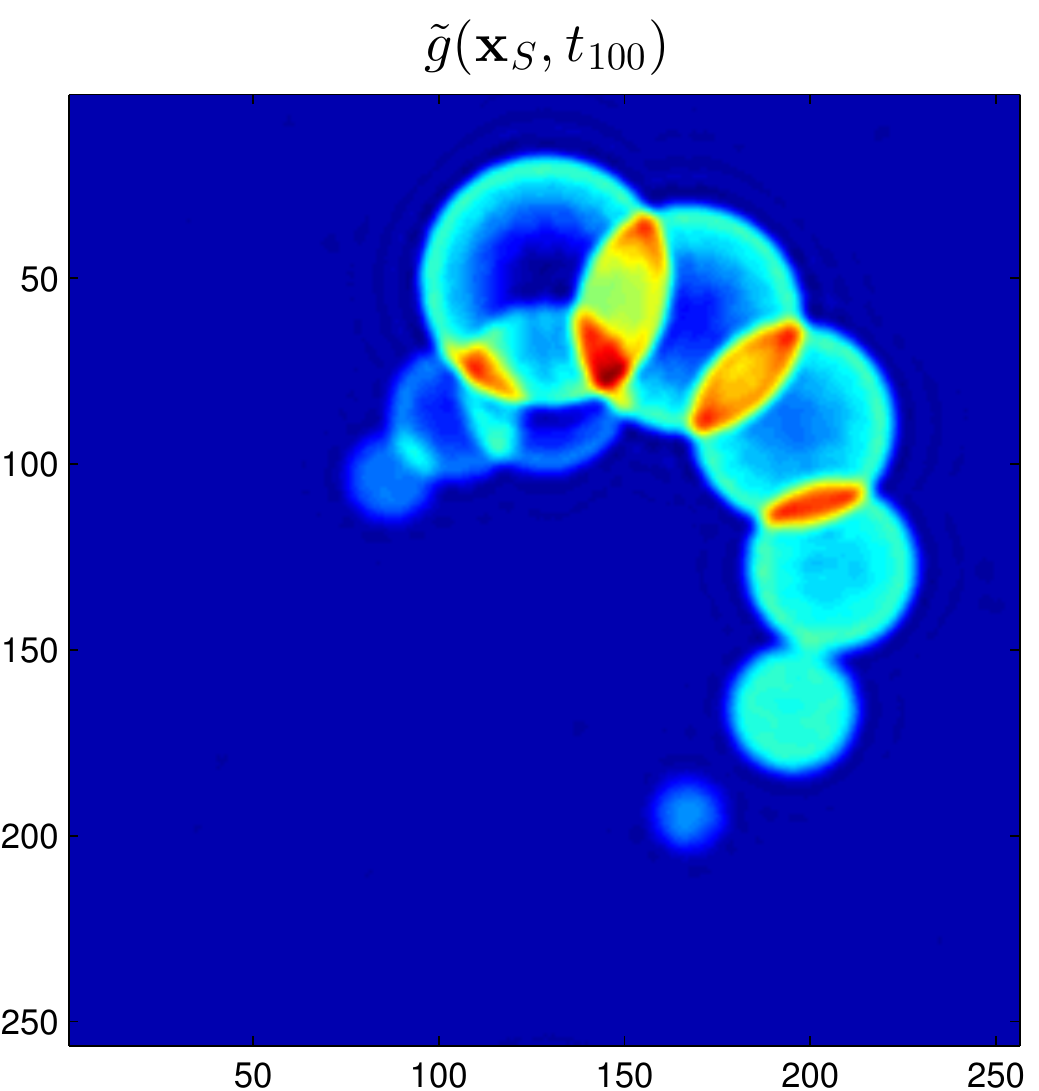}}
\\
\phantom{\hspace{0.001\textwidth}}
\subfloat[][$\Psi g(\bx_{\mc S},t_{100})$]{
\includegraphics[width=0.13\textwidth,trim={0 0 0em 2.5em},clip]{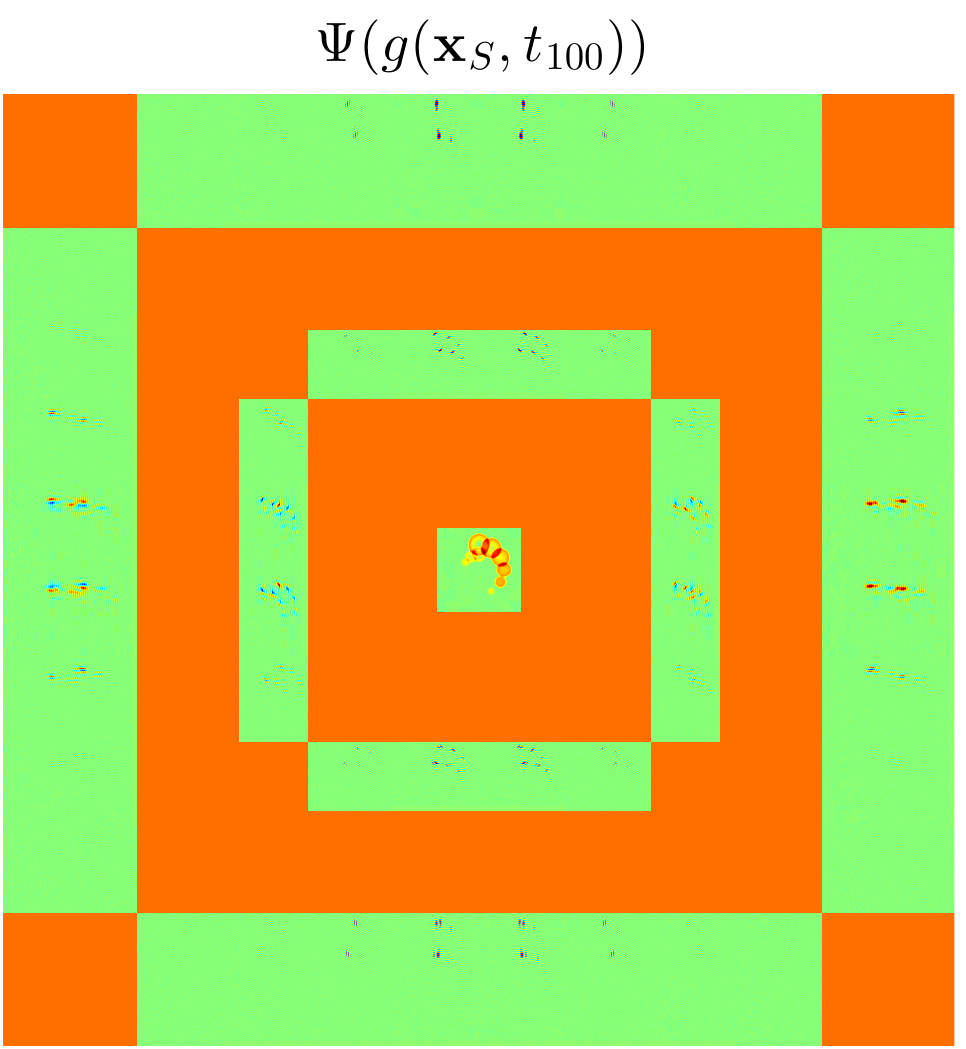}}
\hspace{0.02\textwidth}
\subfloat[][$\Psi \hat g(\bx_{\mc S},t_{100})$]{
\includegraphics[width=0.13\textwidth,trim={0 0 0em 2.5em},clip]{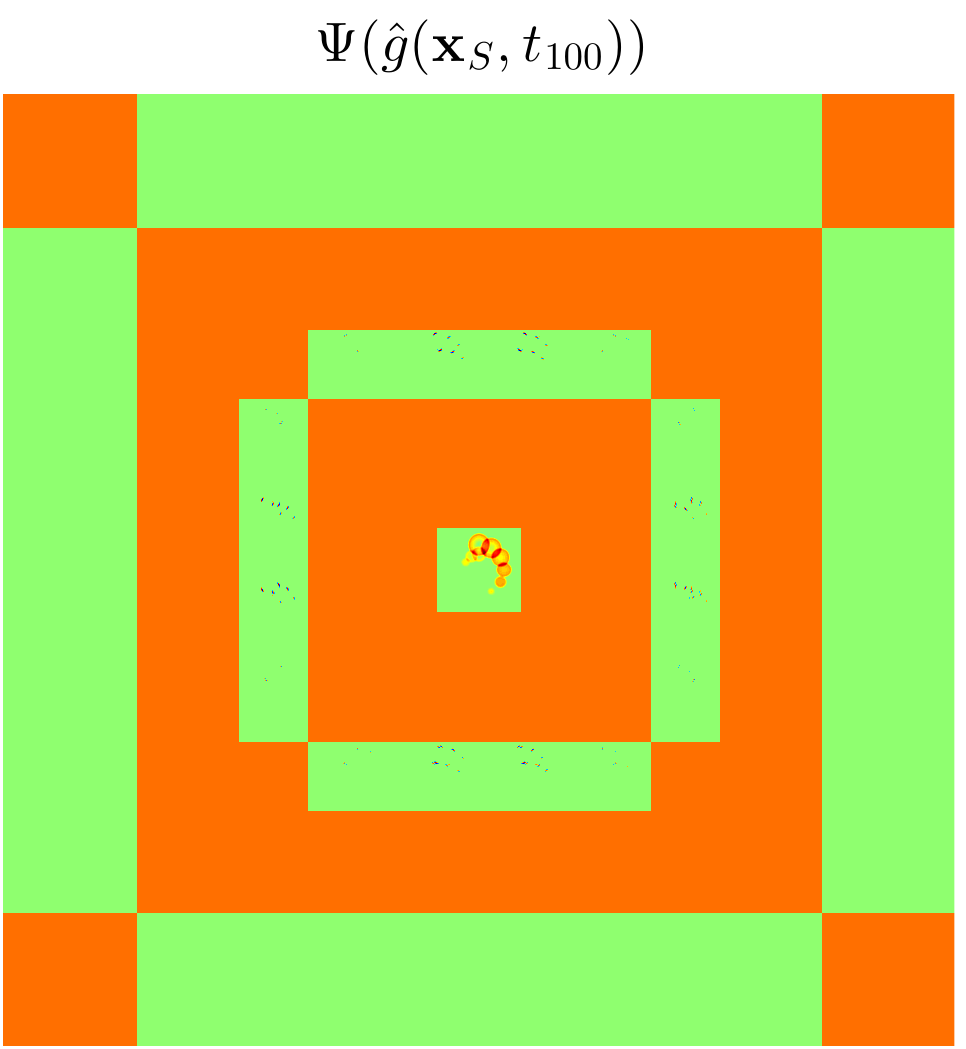}}
\hspace{0.018\textwidth}
\subfloat[][$\Psi \tilde g(\bx_{\mc S},t_{100})$]{
\includegraphics[width=0.13\textwidth,trim={0 0 0em 2.5em},clip]{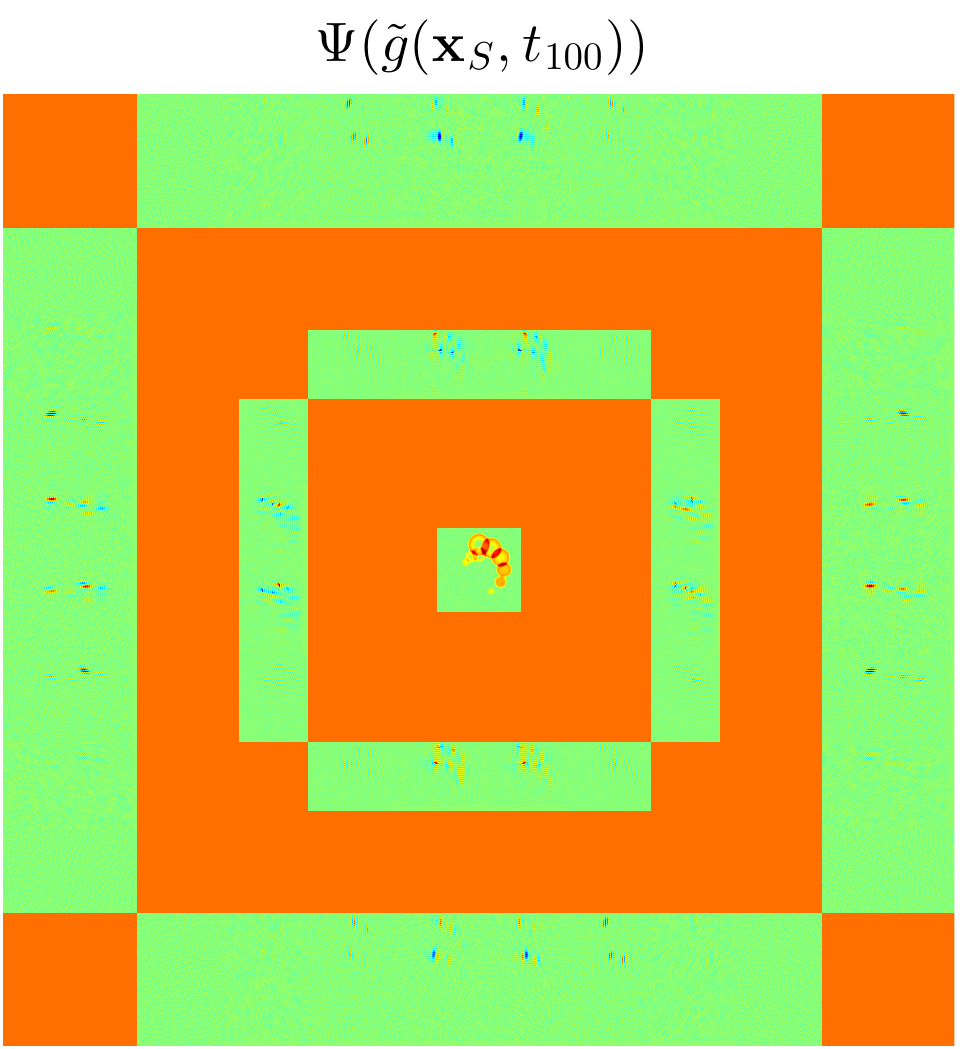}}
\\
\phantom{\hspace{0.16\textwidth}}
\subfloat[][Error of $\hat g(\bx_{\mc S},t_{100})$]{ 
\includegraphics[width=0.155\textwidth,trim={0.5em 0 1.5em 2.5em},clip]{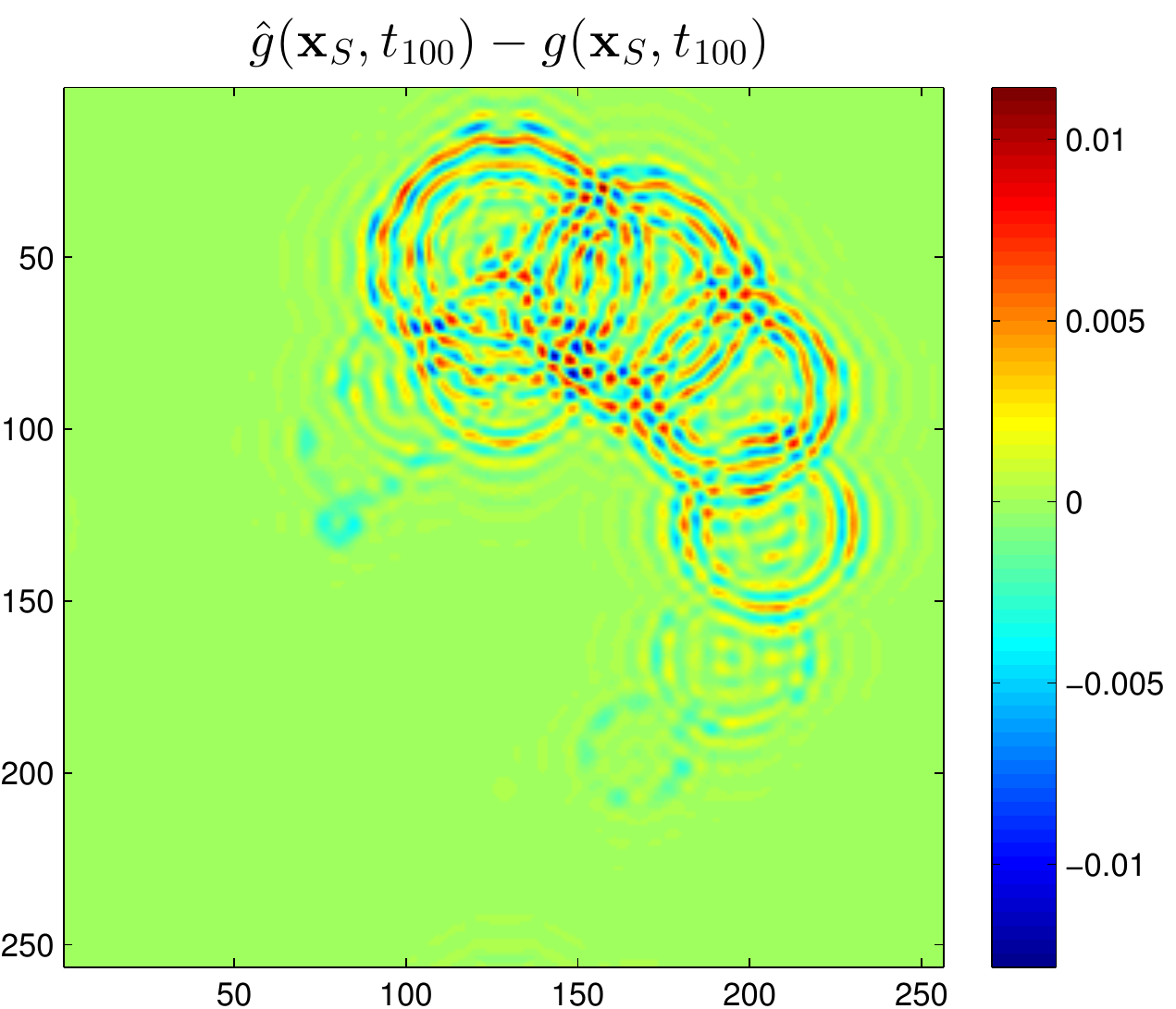}}
\subfloat[][Error of $\tilde g(\bx_{\mc S},t_{100})$]{ 
\includegraphics[width=0.155\textwidth,trim={0.5em 0 1em 2.5em},clip]{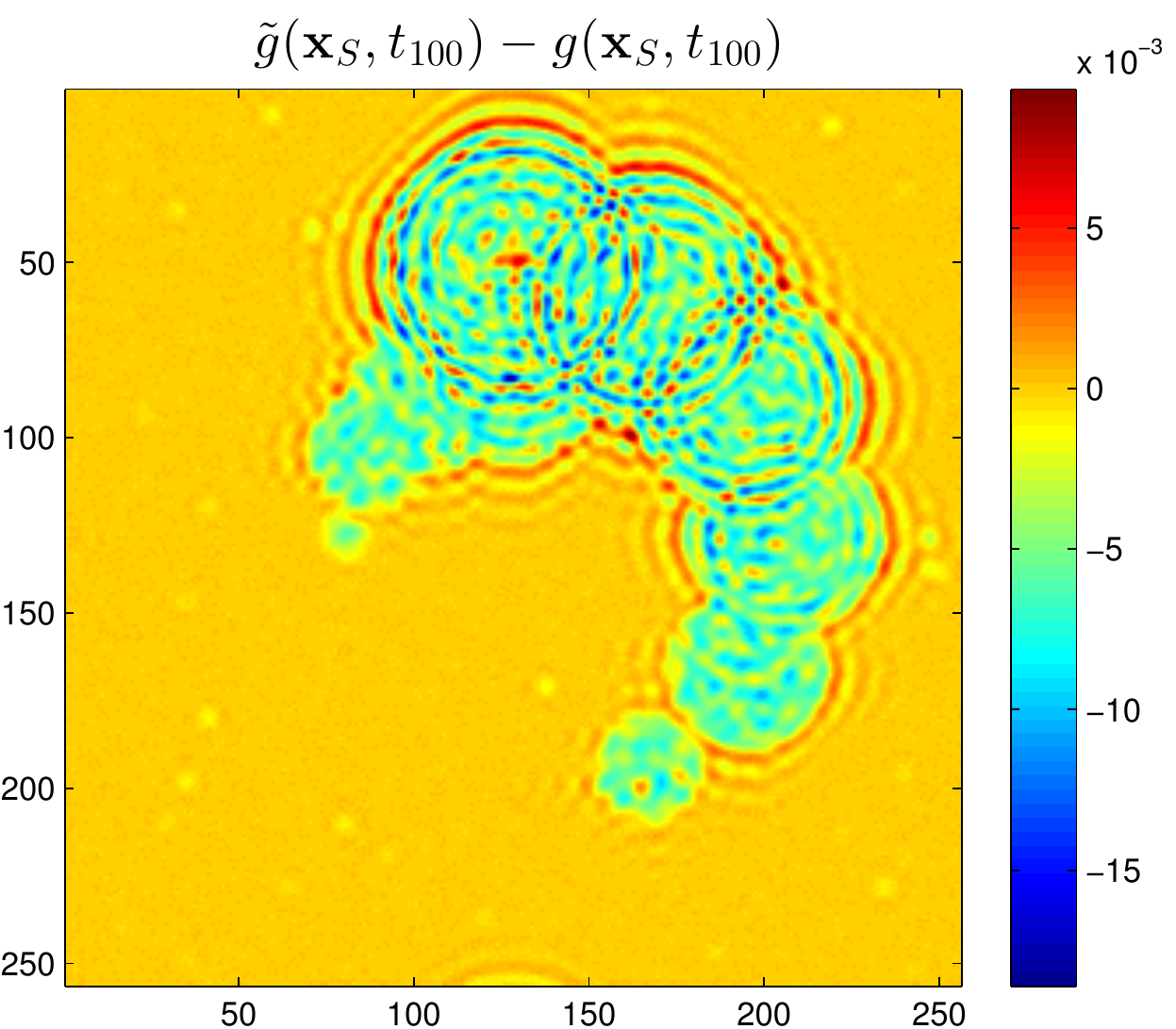}}
\caption{Clock phantom. PAT data at time step $t_{100}$, (a) simulated full data $g(\bx_{\mc S},t)$, (b) compressed $\hat g(\bx_{\mc S},t)$ and (c) reconstructed $\tilde g(\bx_{\mc S},t)$ data with standard Curvelet transform, $\mc C_3^{256,256}$.
The corresponding Curvelet coefficients (d-f) and (g) compression error $\hat g(\bx_{\mc S},t) - g(\bx_{\mc S},t)$, (h) recovery error $\tilde g(\bx_{\mc S},t) - g(\bx_{\mc S},t)$.}
\label{fig:gt100C}
\end{figure}

\begin{figure}[ht]
\subfloat[][$g(\bx_{\mc S},t_{230})$]{ 
\includegraphics[width=0.160\textwidth,trim={0 0 1em 2.5em},clip]{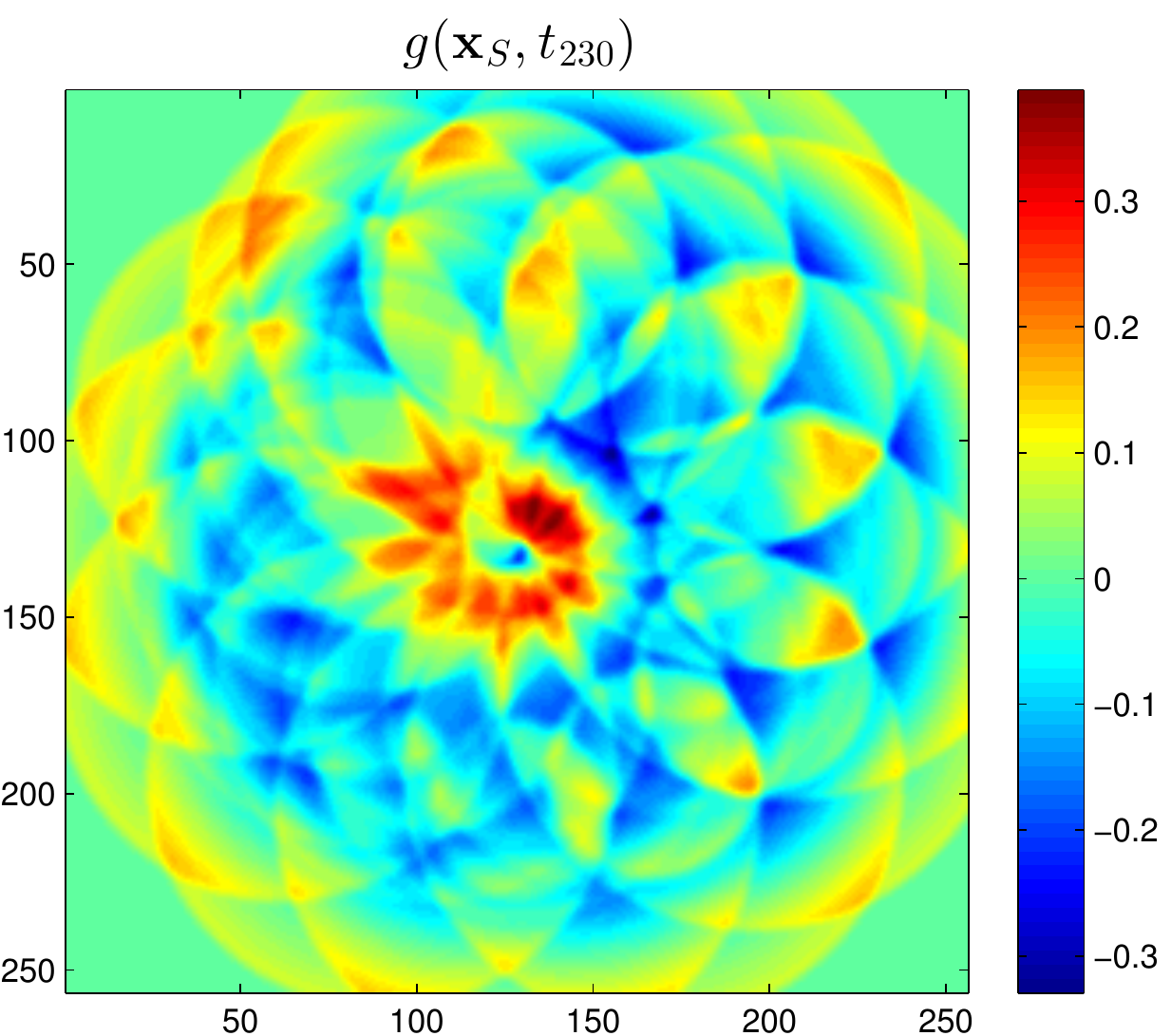}}
\subfloat[][$\hat g(\bx_{\mc S},t_{230})$]{ 
\includegraphics[width=0.137\textwidth,trim={0 0 0em 2.5em},clip]{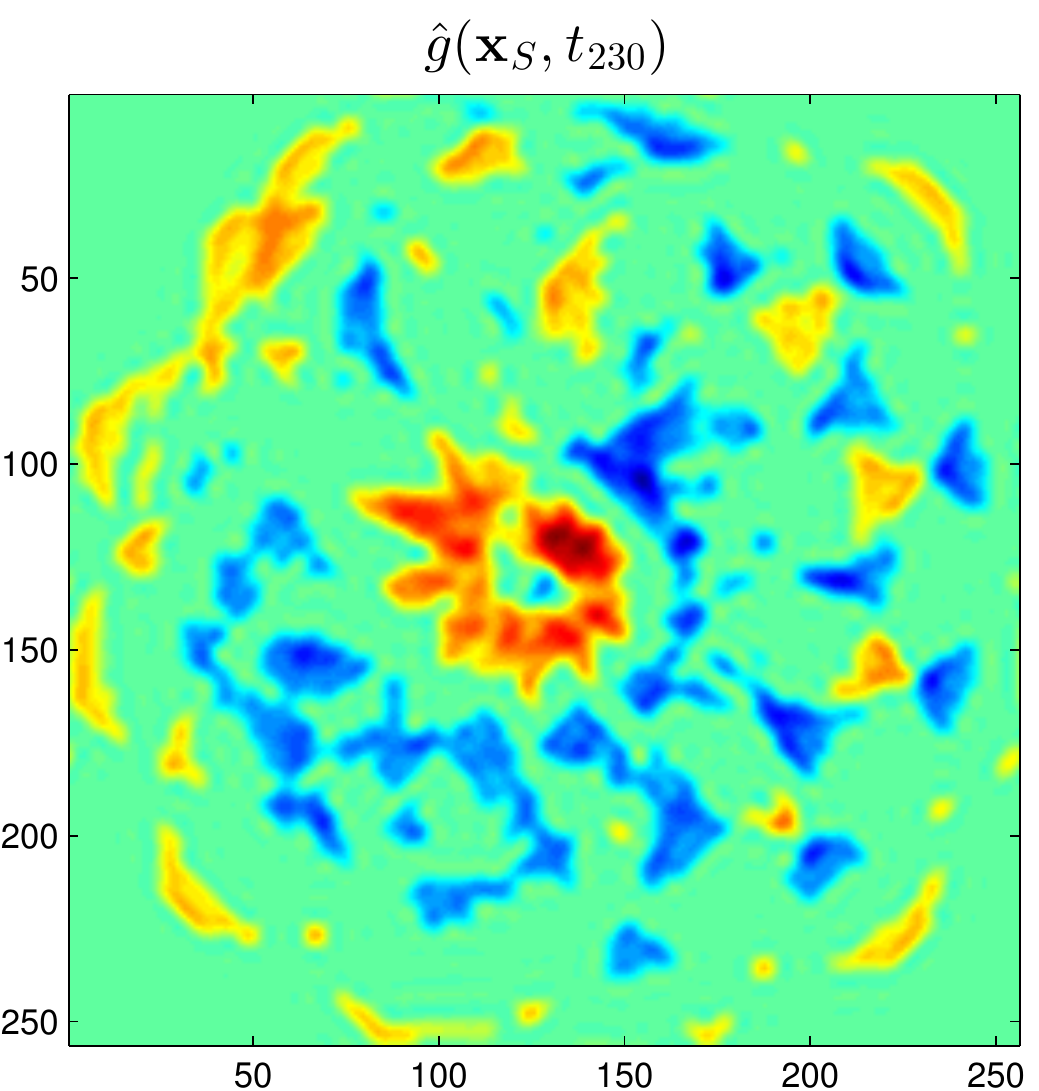}}
\hspace{0.008\textwidth}
\subfloat[][$\tilde g(\bx_{\mc S},t_{230})$]{ 
\includegraphics[width=0.137\textwidth,trim={0 0 0em 2.5em},clip]{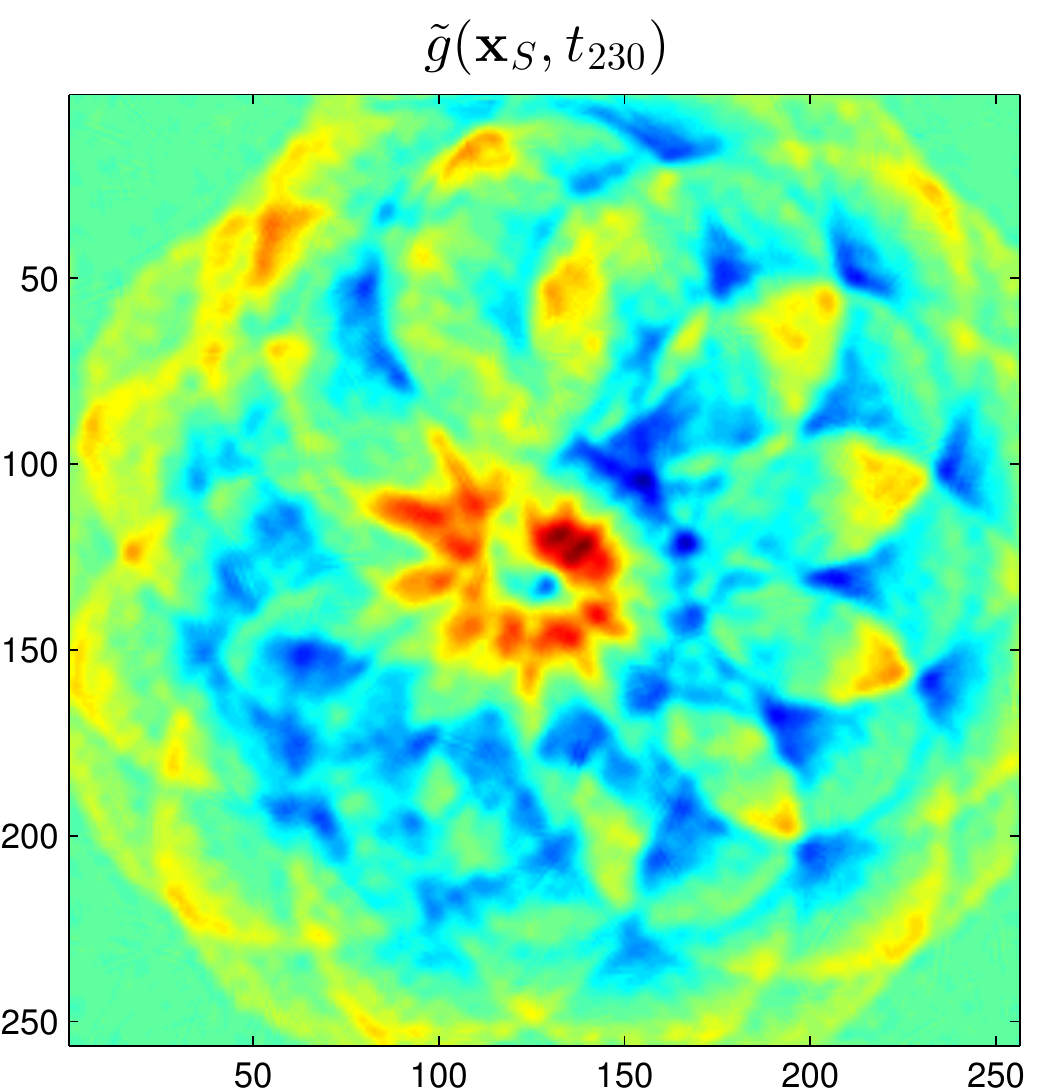}}
\\
\dontshow{
\phantom{\hspace{0.001\textwidth}}
\subfloat[][$\Psi g(\bx_{\mc S},t_{230})$]{
\includegraphics[width=0.13\textwidth,trim={0 0 0em 2.5em},clip]{clock_93x256x256_im_psi_sd_t230_fast_hadamard_11802x65536_curvelet_l3_real_fcurv_L1_tau0p01_SALSA_v2_stop1_tol0p0005_maxit100_mu500}}
\hspace{0.02\textwidth}
\subfloat[][$\Psi \hat g(\bx_{\mc S},t_{230})$]{
\includegraphics[width=0.13\textwidth,trim={0 0 0em 2.5em},clip]{clock_93x256x256_im_psi_csd_t230_fast_hadamard_11802x65536_curvelet_l3_real_fcurv_L1_tau0p01_SALSA_v2_stop1_tol0p0005_maxit100_mu500}}
\hspace{0.018\textwidth}
\subfloat[][$\Psi \tilde g(\bx_{\mc S},t_{230})$]{
\includegraphics[width=0.13\textwidth,trim={0 0 0em 2.5em},clip]{clock_93x256x256_im_psi_rsd_t230_fast_hadamard_11802x65536_curvelet_l3_real_fcurv_L1_tau0p01_SALSA_v2_stop1_tol0p0005_maxit100_mu500}}
\\
}
\phantom{\hspace{0.16\textwidth}}
\subfloat[][Error of $\hat g(\bx_{\mc S},t_{230})$]{ 
\includegraphics[width=0.155\textwidth,trim={0.5em 0 1.5em 2.5em},clip]{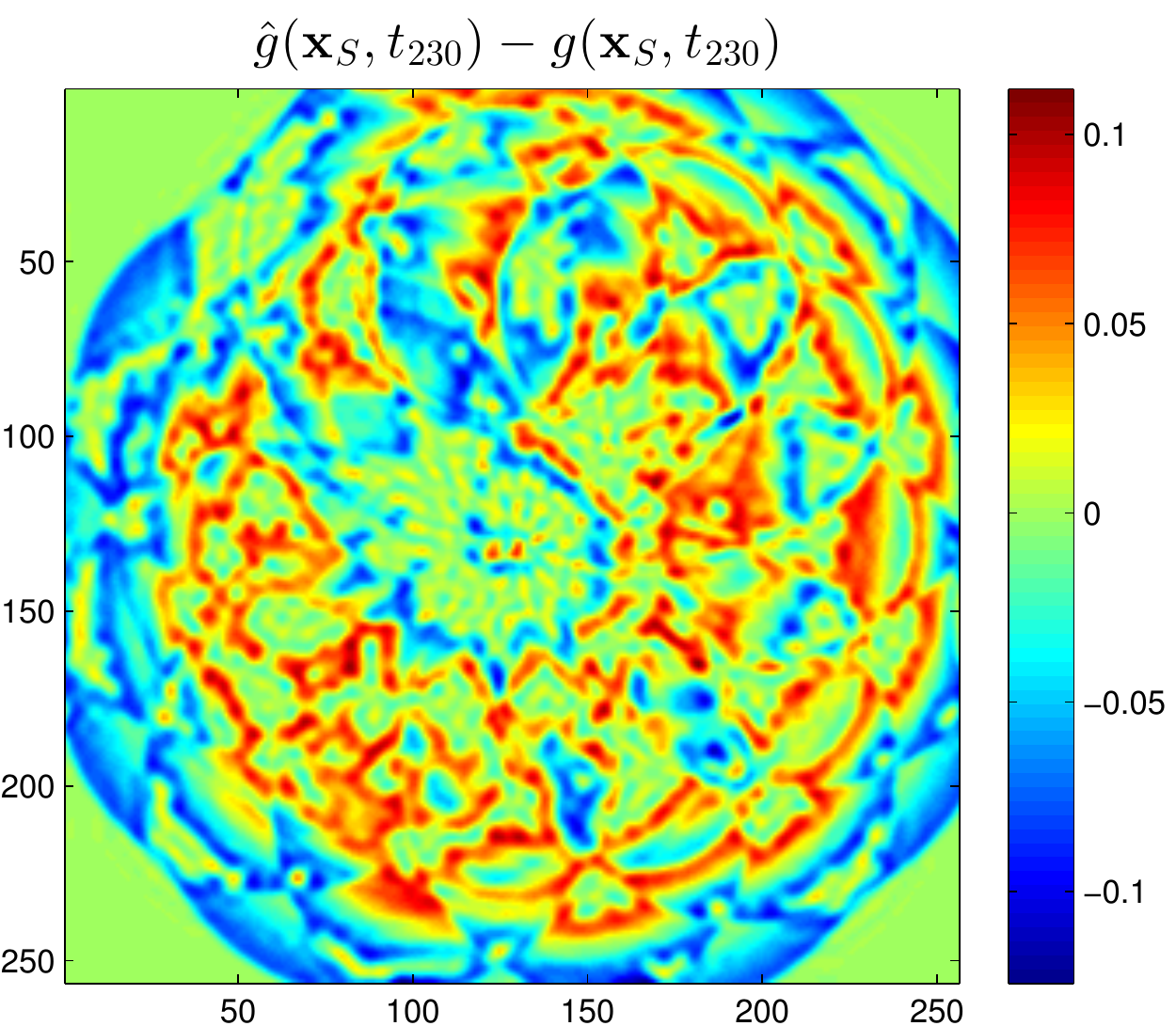}}
\subfloat[][Error of $\tilde g(\bx_{\mc S},t_{230})$]{ 
\includegraphics[width=0.157\textwidth,trim={0.5em 0 1em 2.5em},clip]{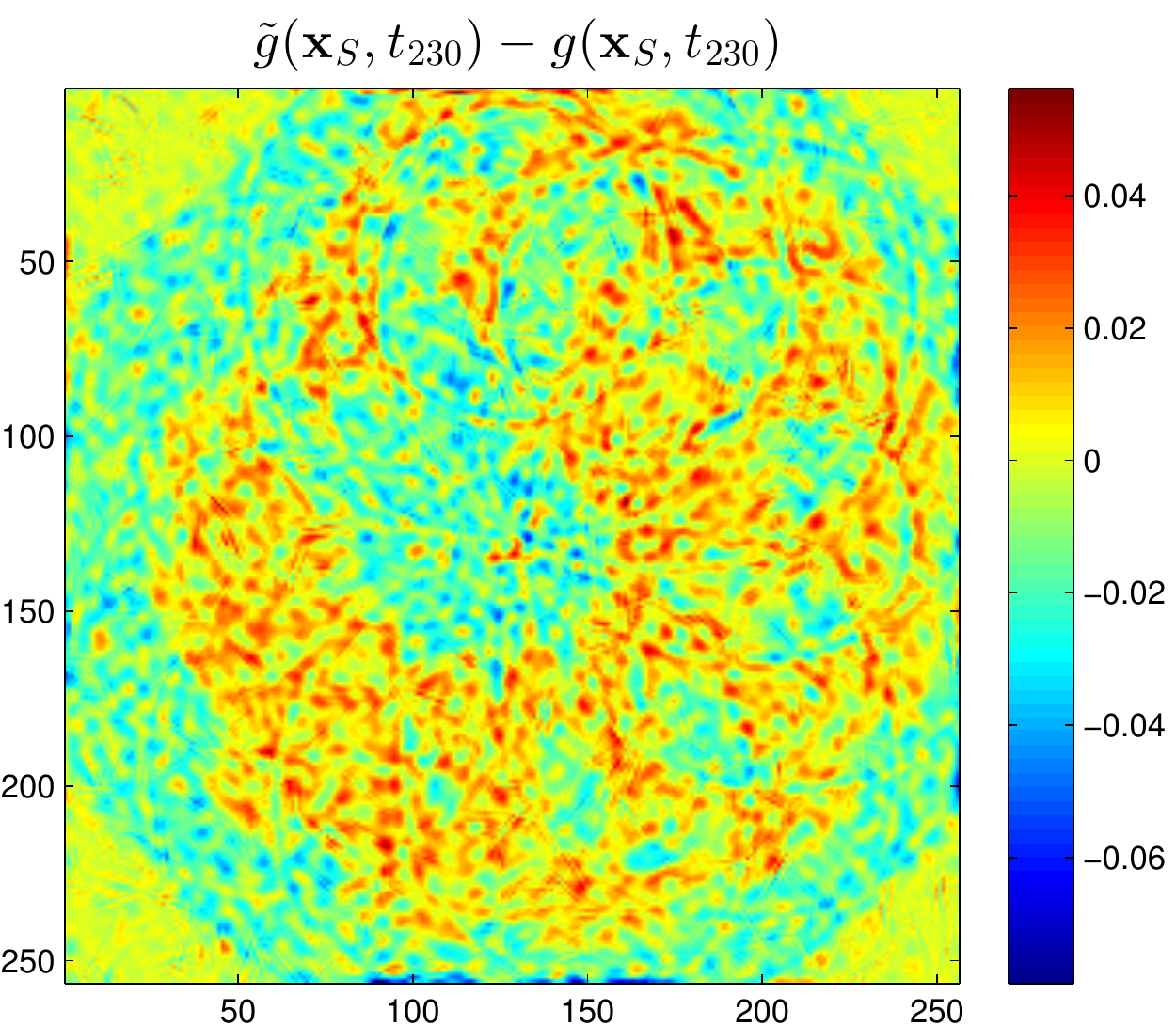}}
\caption{Clock phantom. PAT data at time step $t_{230}$, (a) simulated full data $g(\bx_{\mc S},t)$, (b) compressed $\hat g(\bx_{\mc S},t)$ and (c) reconstructed $\tilde g(\bx_{\mc S},t)$ data with standard Curvelet transform, $\mc C_3^{256,256}$. (d) Compression error $\hat g(\bx_{\mc S},t) - g(\bx_{\mc S},t)$, (e) recovery error $\tilde g(\bx_{\mc S},t) - g(\bx_{\mc S},t)$.}
\label{fig:gt230C}
\end{figure}

First, we consider the approximation properties of the standard Curvelet and low-frequency Curvelet transforms.
Figure \ref{fig:curv}(a-b) shows the decay of the amplitudes of Curvelet coefficients for standard Curvelets $\mc C^{256,256}_{3}$ and low-frequency Curvelets $\mc C^{256,256}_{3,192,192}$ of the photoacoustic field at time steps $t_i=100, 230$ corresponding to the fields depicted in Figures \ref{fig:gt100C}(a), \ref{fig:gt230C}(a). For the largest amplitude coefficients the low-frequency Curvelets consistently exhibit a quicker decay and correspondingly lower approximation error shown in Figure \ref{fig:curv}(c-d). While in the early time steps $t_i=100$ the decay rate difference is most pronounced for later time steps the difference gets smaller but is distributed over more coefficients. In all cases eventually the approximation error of standard Curvelet transform falls below that of the low-frequency Curvelet transform but this is only after all the significant coefficients has been captured.  

\begin{figure}[ht]
\centering
\subfloat[][$\Psi = \mc C_3^{256,256}$]{
\includegraphics[width=0.225\textwidth,trim={0 0 0 1em},clip]{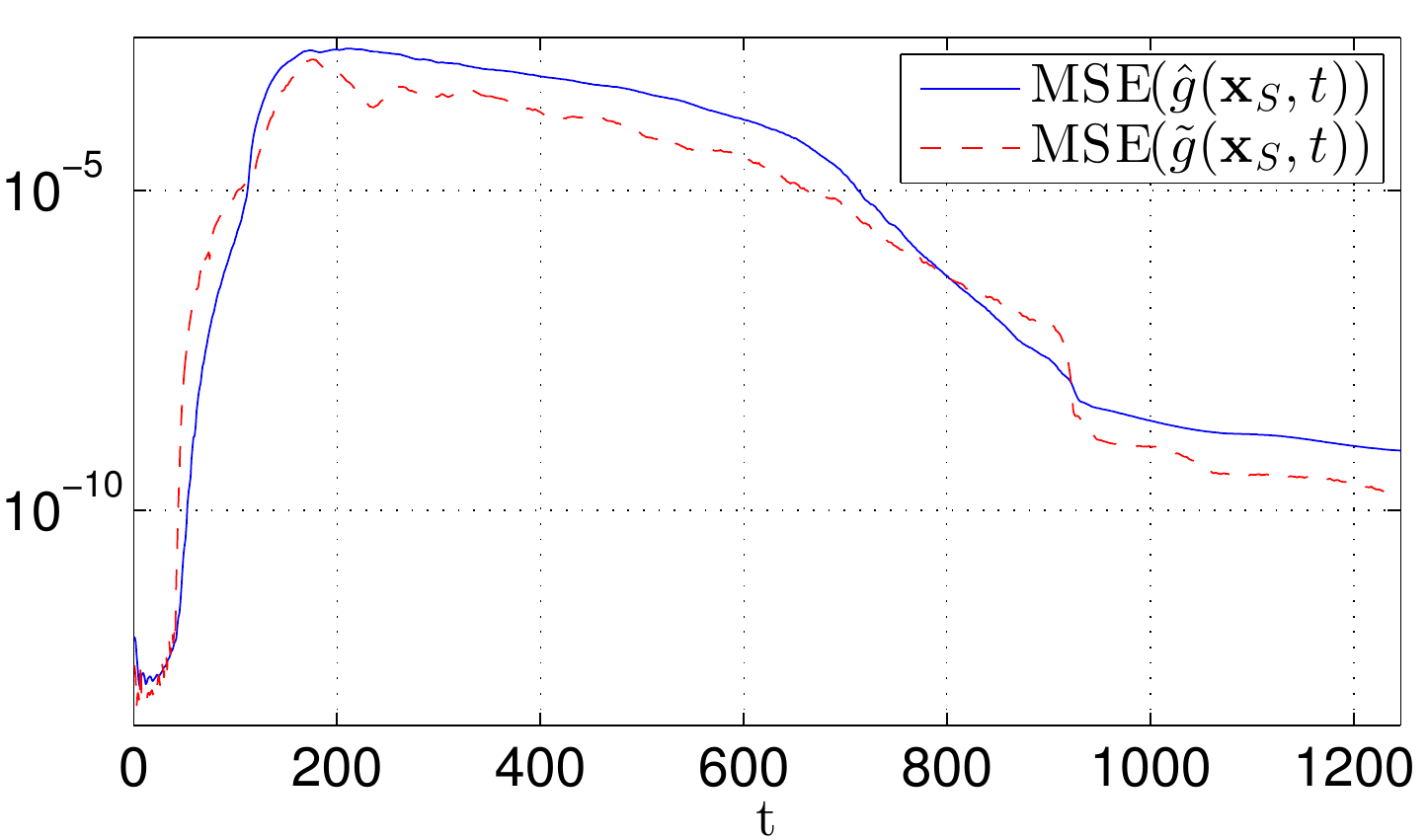}}\hfill
\subfloat[][$\Psi = \mc C_{3,192,192}^{256,256}$]{
\includegraphics[width=0.225\textwidth,trim={0 0 0 1em},clip]{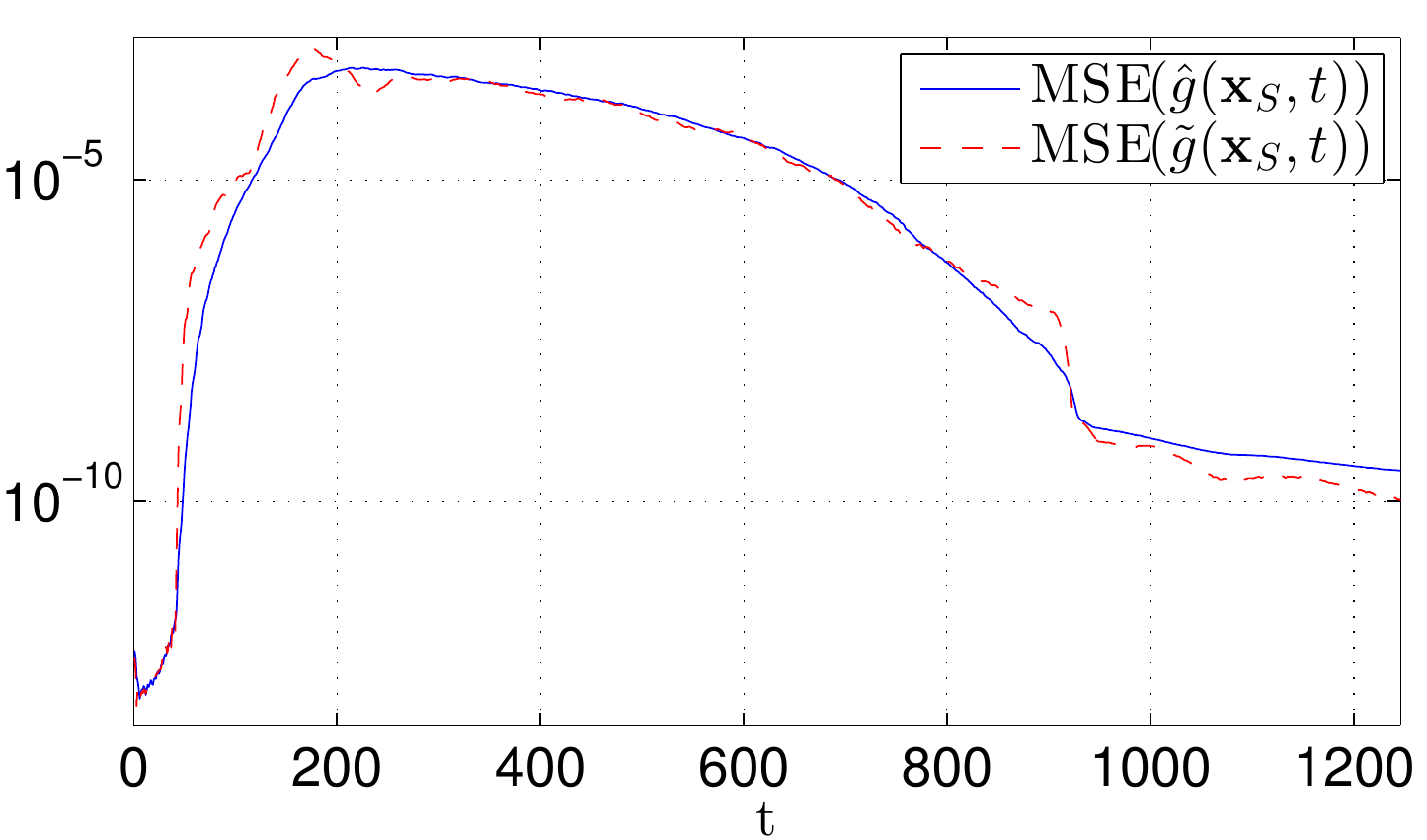}}
\caption{MSE of the compressed PAT data $\hat g(\bx_{\mc S},t)$ versus the reconstructed PAT data $\tilde g(\bx_{\mc S},t)$ over time, for (a) standard Curvelet transform, $\mc C_3^{256,256}$; (b) low-frequency Curvelet transform, $\mc C_{3,192,192}^{256,256}$.}
\label{fig:mse}
\end{figure}

\begin{figure}[ht]
\subfloat[][$g(\bx_{\mc S},t_{100})$]{ 
\includegraphics[width=0.160\textwidth,trim={0 0 1em 2.5em},clip]{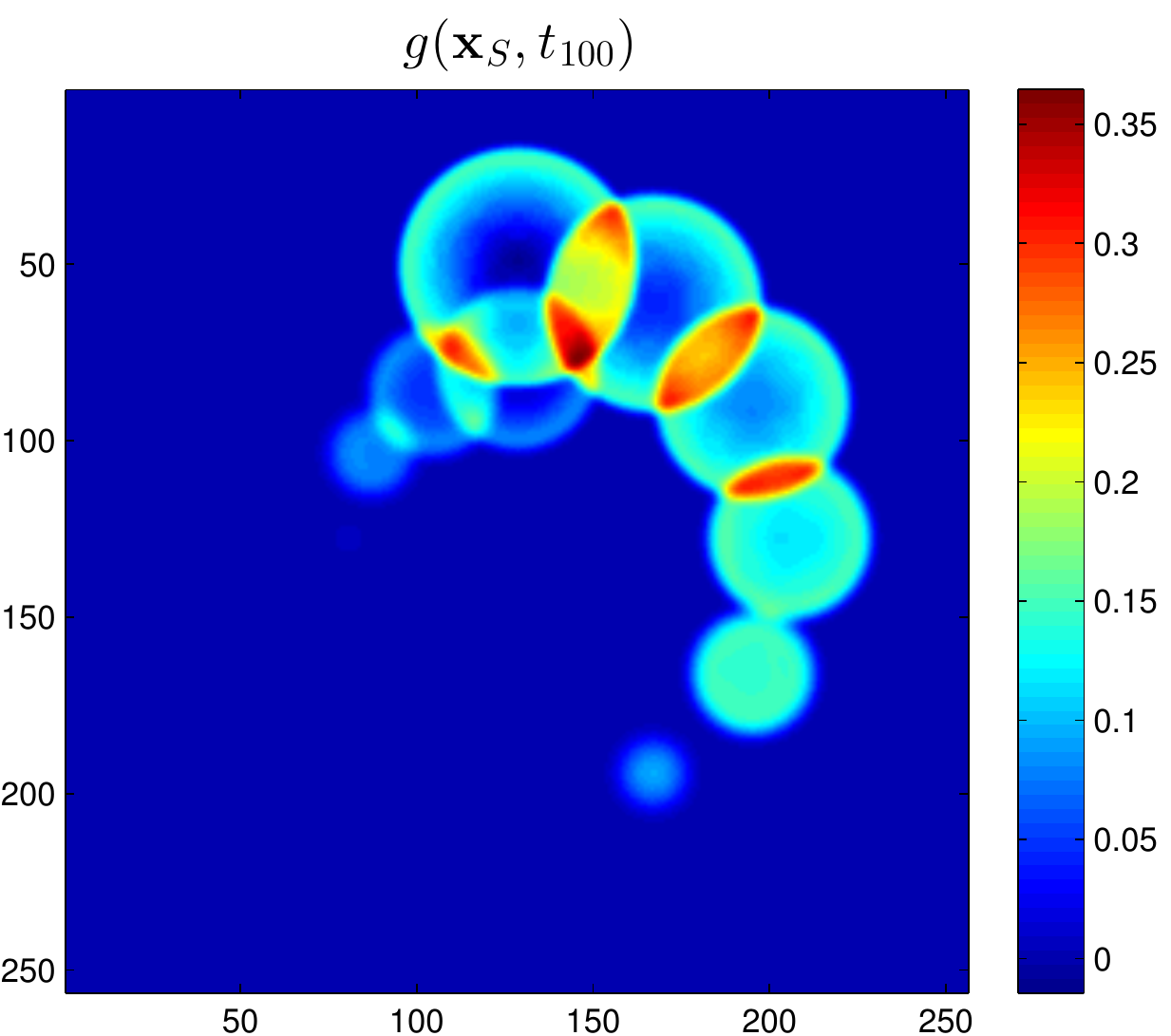}}
\subfloat[][$\hat g(\bx_{\mc S},t_{100})$]{ 
\includegraphics[width=0.137\textwidth,trim={0 0 0em 2.5em},clip]{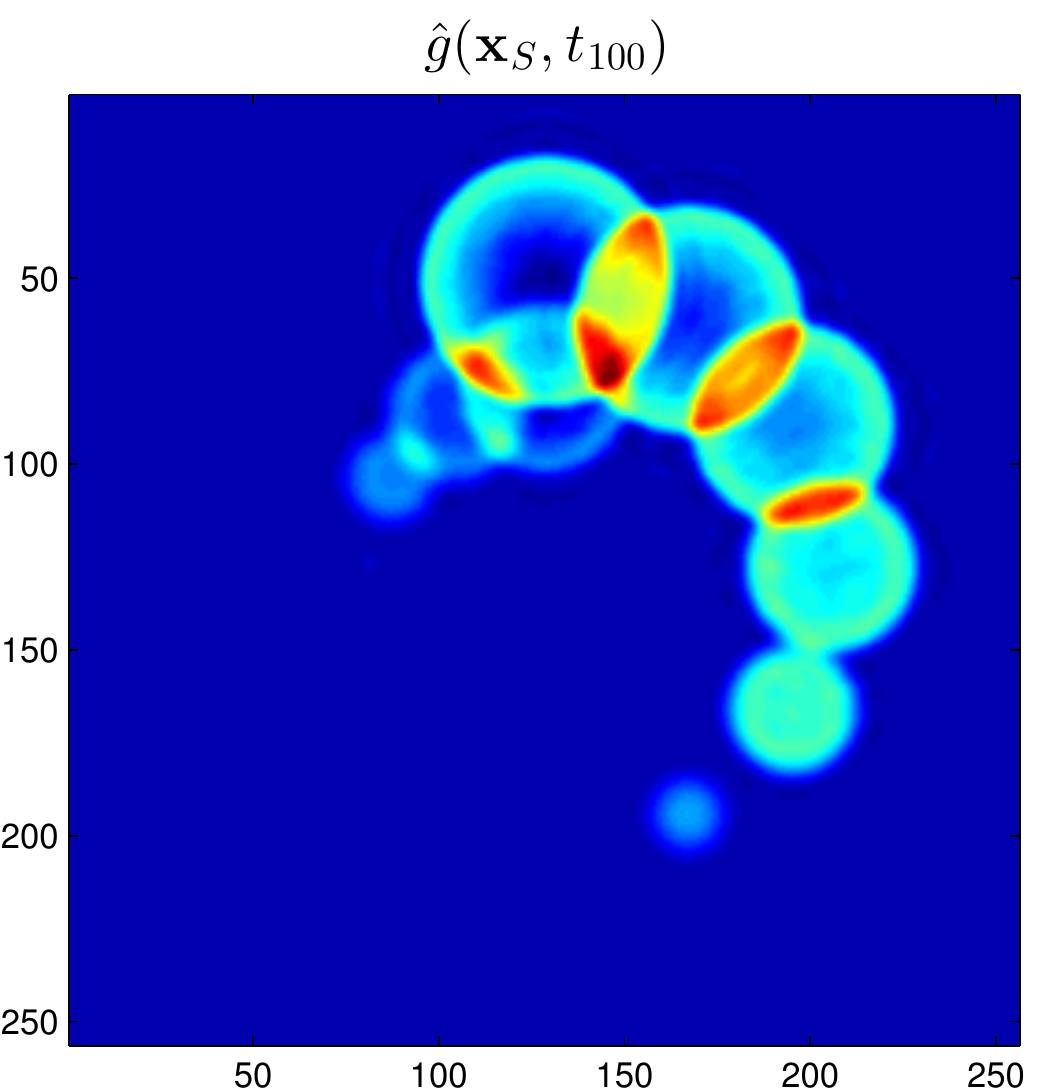}}
\hspace{0.008\textwidth}
\subfloat[][$\tilde g(\bx_{\mc S},t_{100})$]{ 
\includegraphics[width=0.137\textwidth,trim={0 0 0em 2.5em},clip]{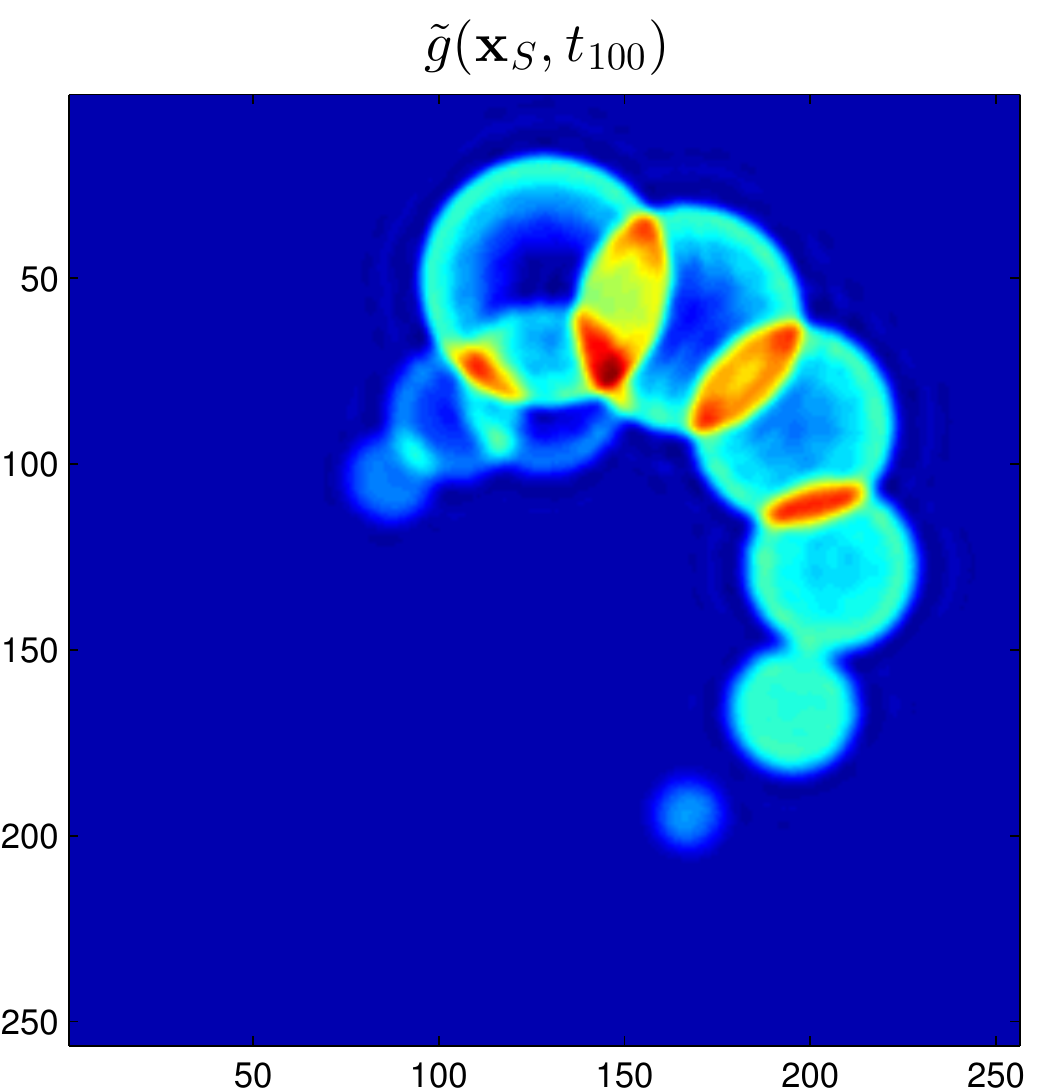}}
\\
\phantom{\hspace{0.001\textwidth}}
\subfloat[][$\Psi g(\bx_{\mc S},t_{100})$]{
\includegraphics[width=0.13\textwidth,trim={0 0 0em 2.5em},clip]{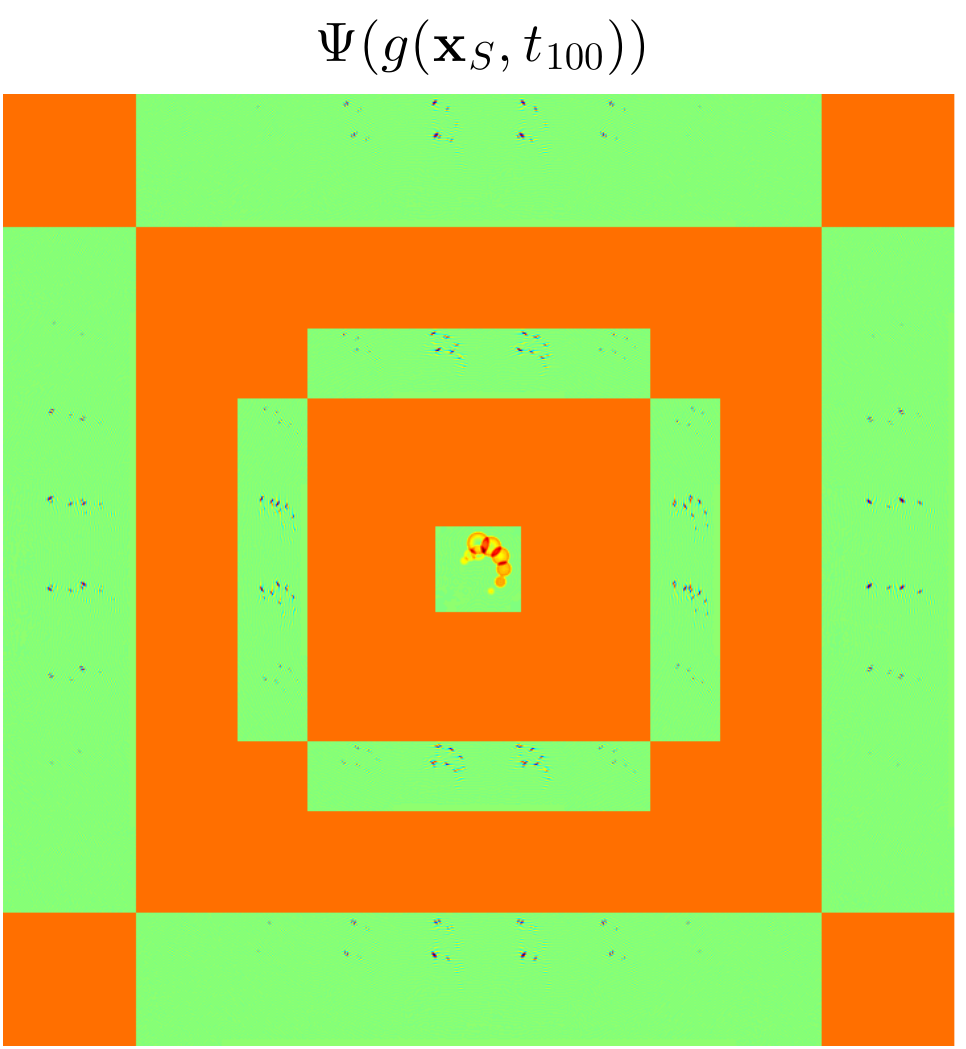}}
\hspace{0.02\textwidth}
\subfloat[][$\Psi \hat g(\bx_{\mc S},t_{100})$]{
\includegraphics[width=0.13\textwidth,trim={0 0 0em 2.5em},clip]{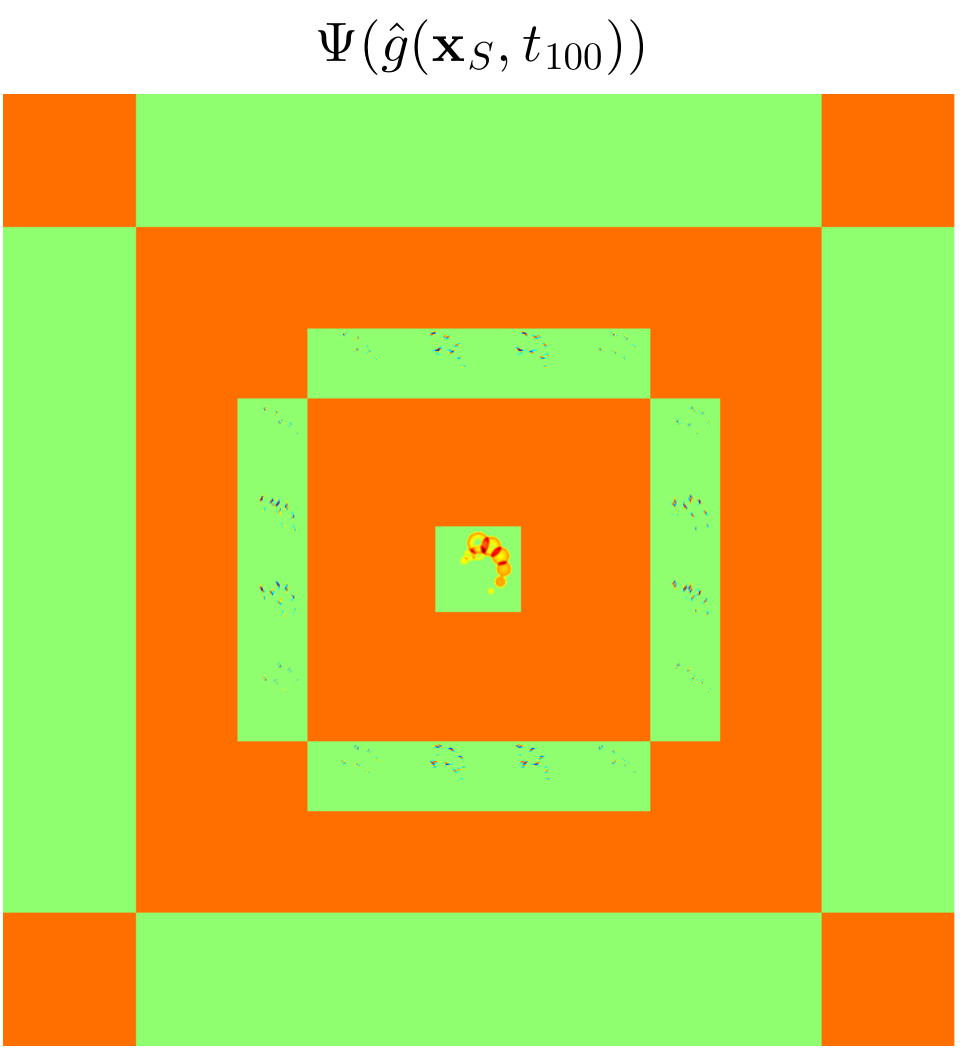}}
\hspace{0.018\textwidth}
\subfloat[][$\Psi \tilde g(\bx_{\mc S},t_{100})$]{
\includegraphics[width=0.13\textwidth,trim={0 0 0em 2.5em},clip]{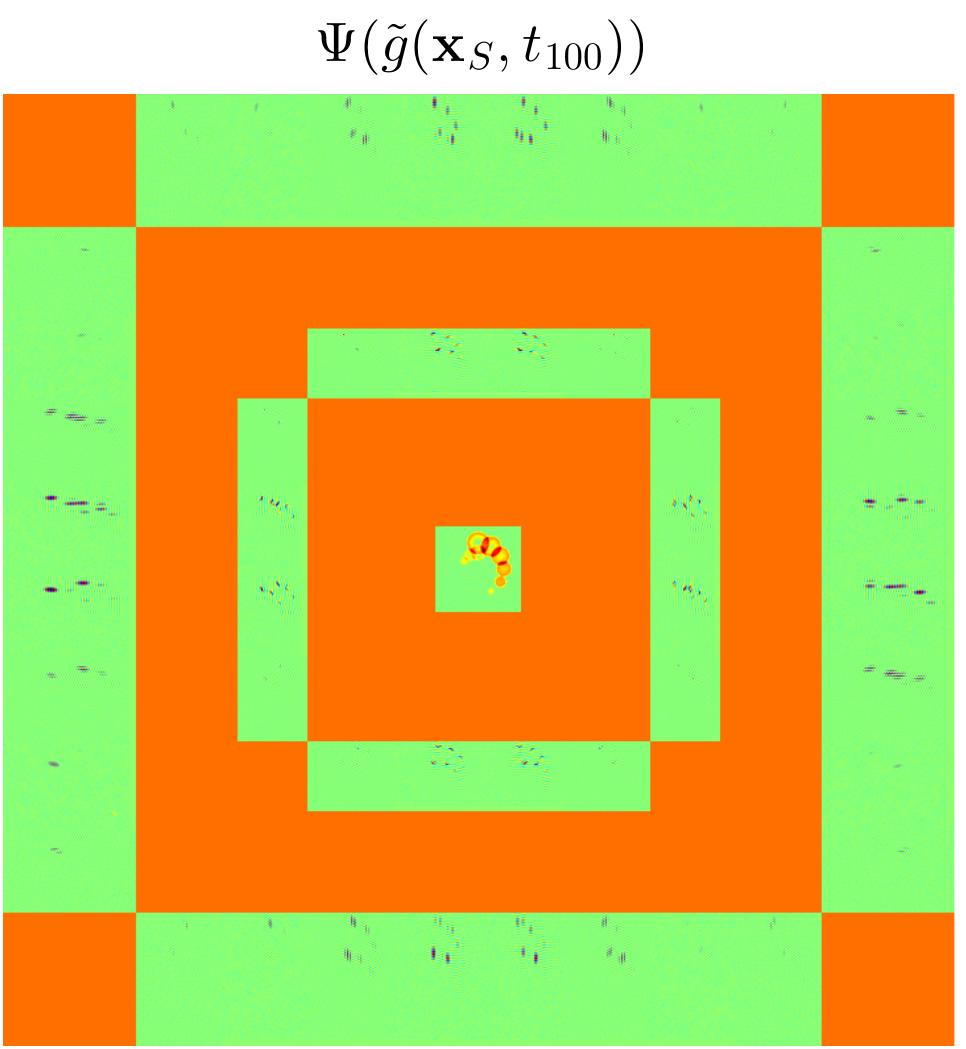}}
\\
\phantom{\hspace{0.16\textwidth}}
\subfloat[][Error of $\hat g(\bx_{\mc S},t_{100})$]{ 
\includegraphics[width=0.155\textwidth,trim={0.5em 0 1.5em 2.5em},clip]{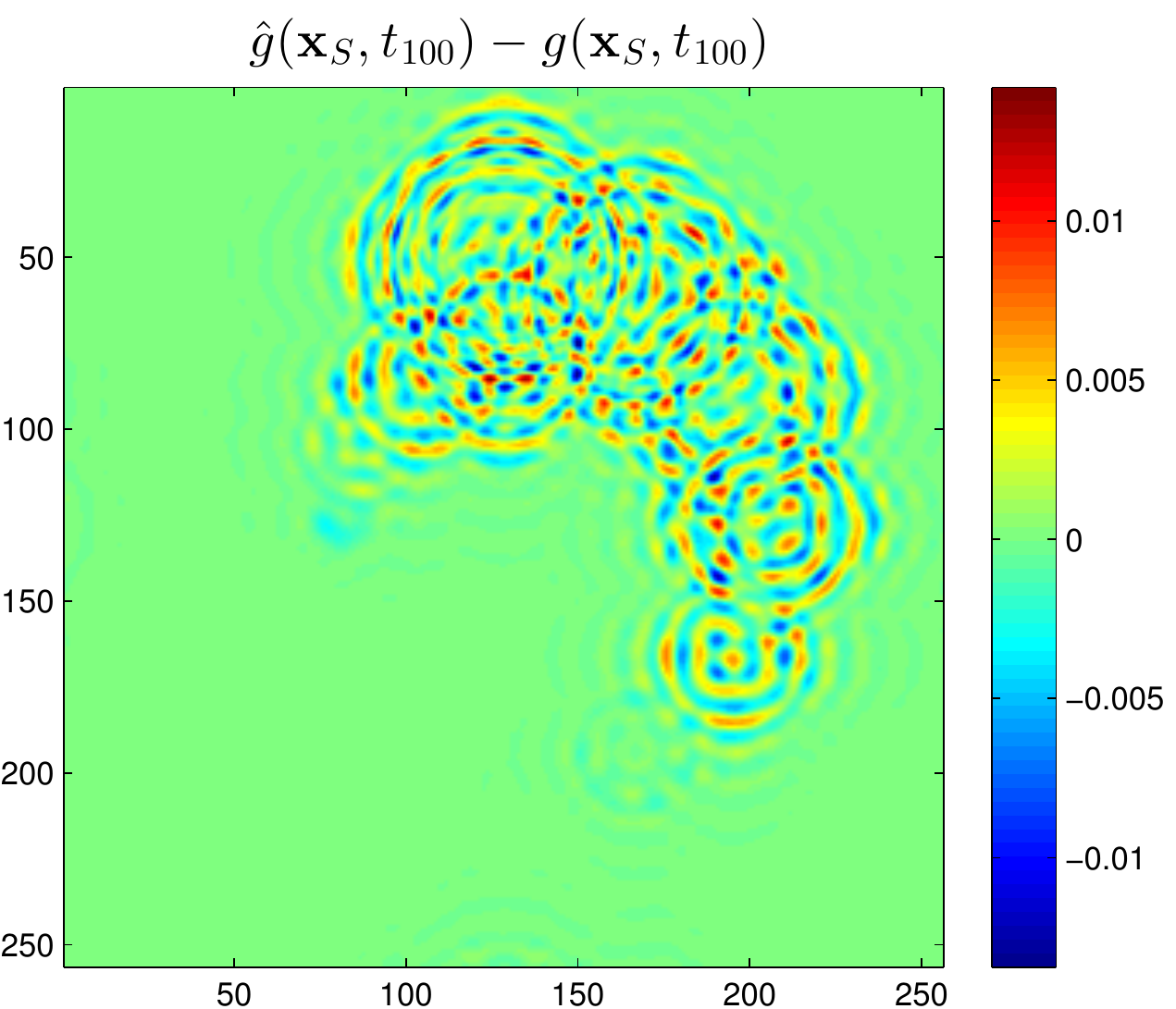}}
\subfloat[][Error of $\tilde g(\bx_{\mc S},t_{100})$]{ 
\includegraphics[width=0.157\textwidth,trim={0.5em 0 1em 2.5em},clip]{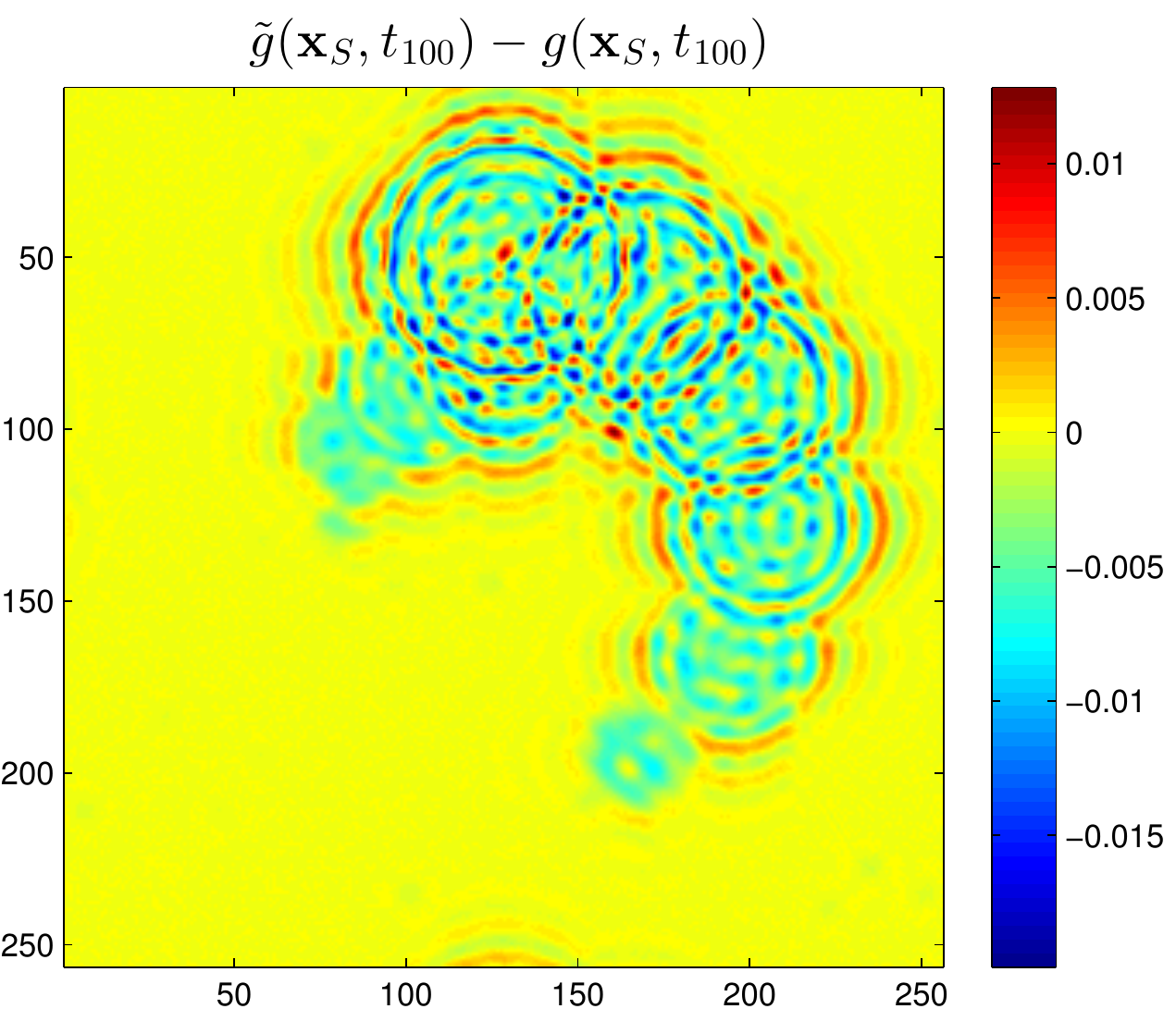}}
\caption{Clock phantom. PAT data at time step $t_{100}$, (a) simulated full data $g(\bx_{\mc S},t)$, (b) compressed $\hat g(\bx_{\mc S},t)$ and (c) reconstructed $\tilde g(\bx_{\mc S},t)$ data with low-frequency Curvelet transform, $\mc C_{3,192,192}^{256,256}$.
The corresponding Curvelet coefficients (d-f) and (g) compression error $\hat g(\bx_{\mc S},t) - g(\bx_{\mc S},t)$, (h) recovery error $\tilde g(\bx_{\mc S},t) - g(\bx_{\mc S},t)$.}
\label{fig:gt100CLP}
\end{figure}

\begin{figure}[ht]
\subfloat[][$g(\bx_{\mc S},t_{230})$]{ 
\includegraphics[width=0.160\textwidth,trim={0 0 1em 2.5em},clip]{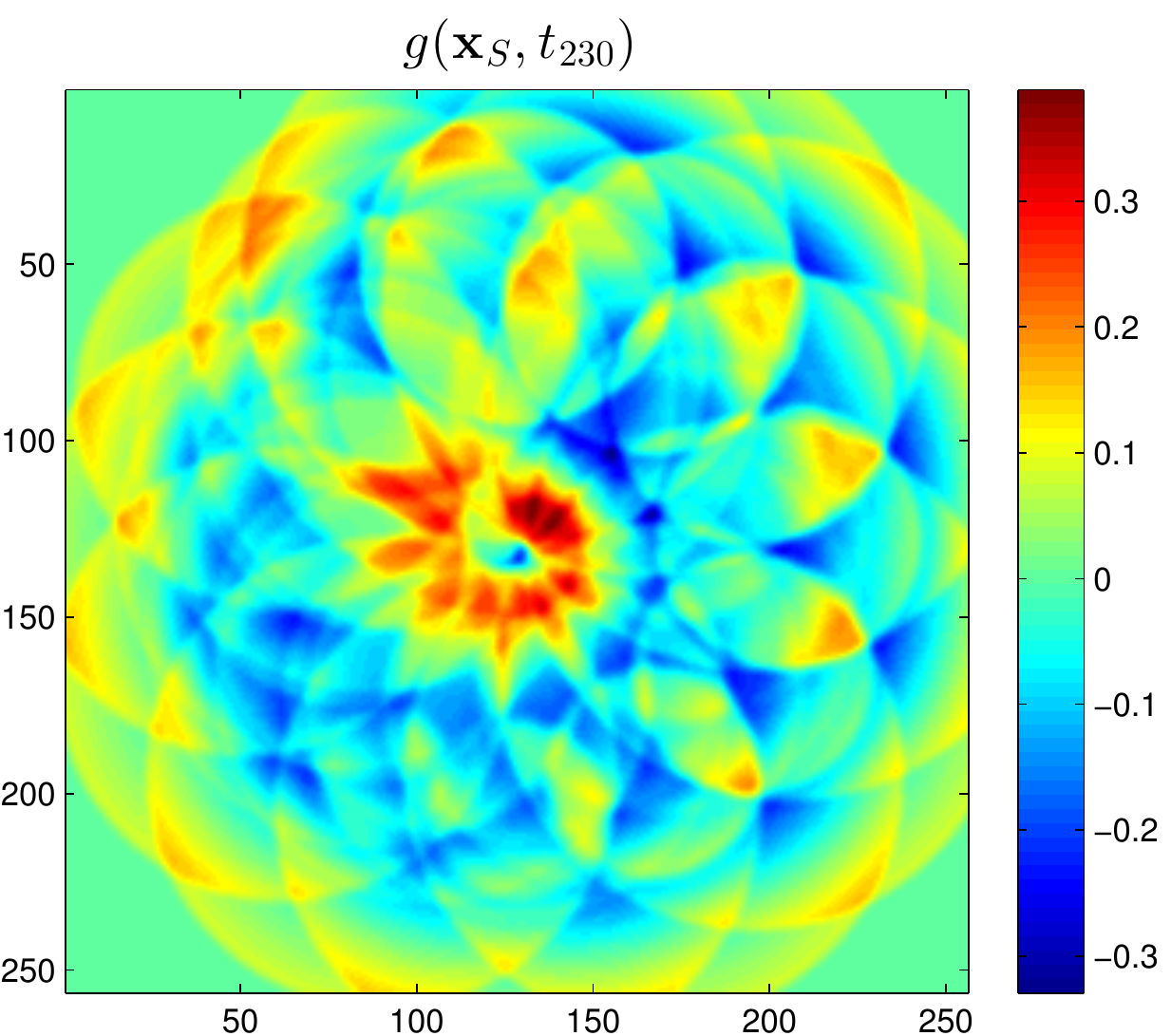}}
\subfloat[][$\hat g(\bx_{\mc S},t_{230})$]{ 
\includegraphics[width=0.137\textwidth,trim={0 0 0em 2.5em},clip]{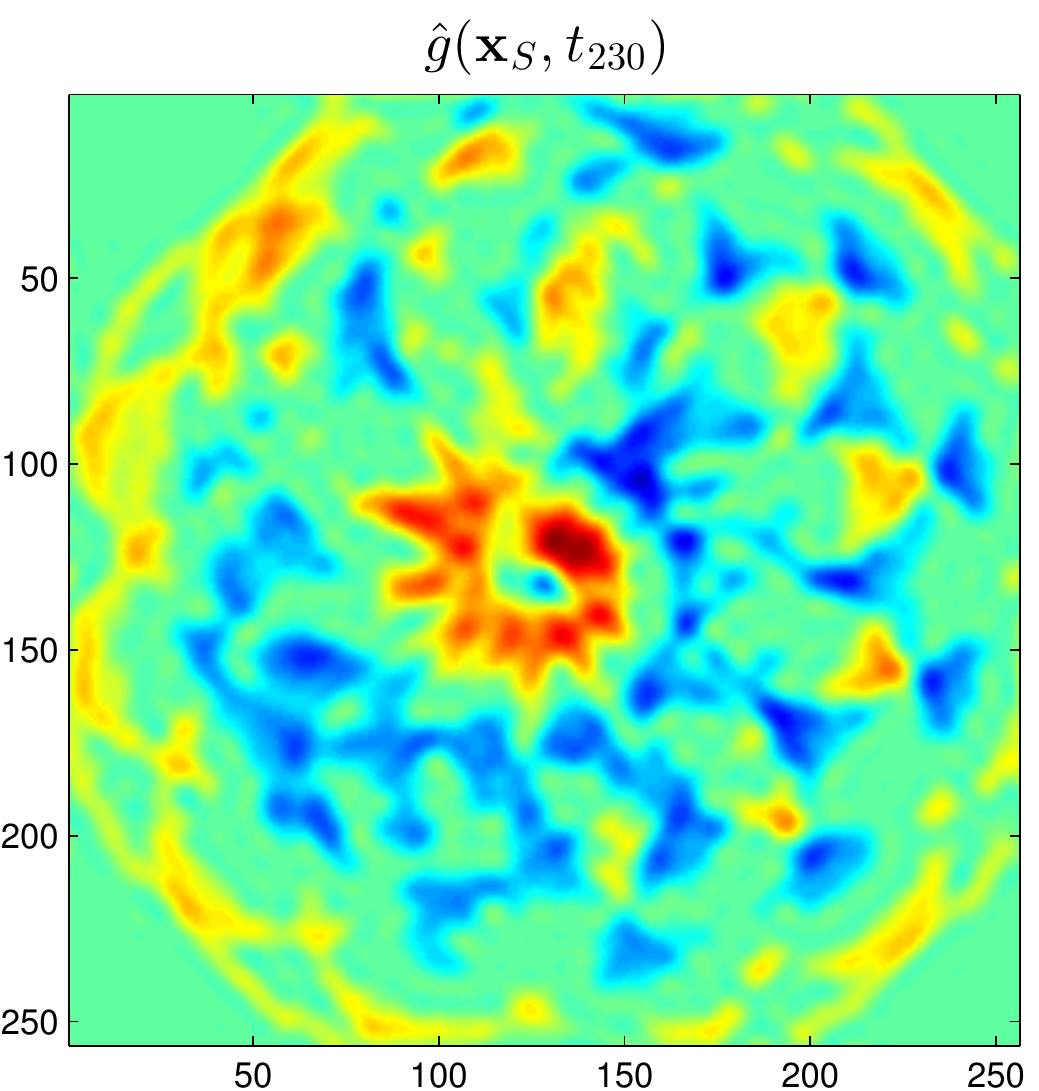}}
\hspace{0.008\textwidth}
\subfloat[][$\tilde g(\bx_{\mc S},t_{230})$]{ 
\includegraphics[width=0.137\textwidth,trim={0 0 0em 2.5em},clip]{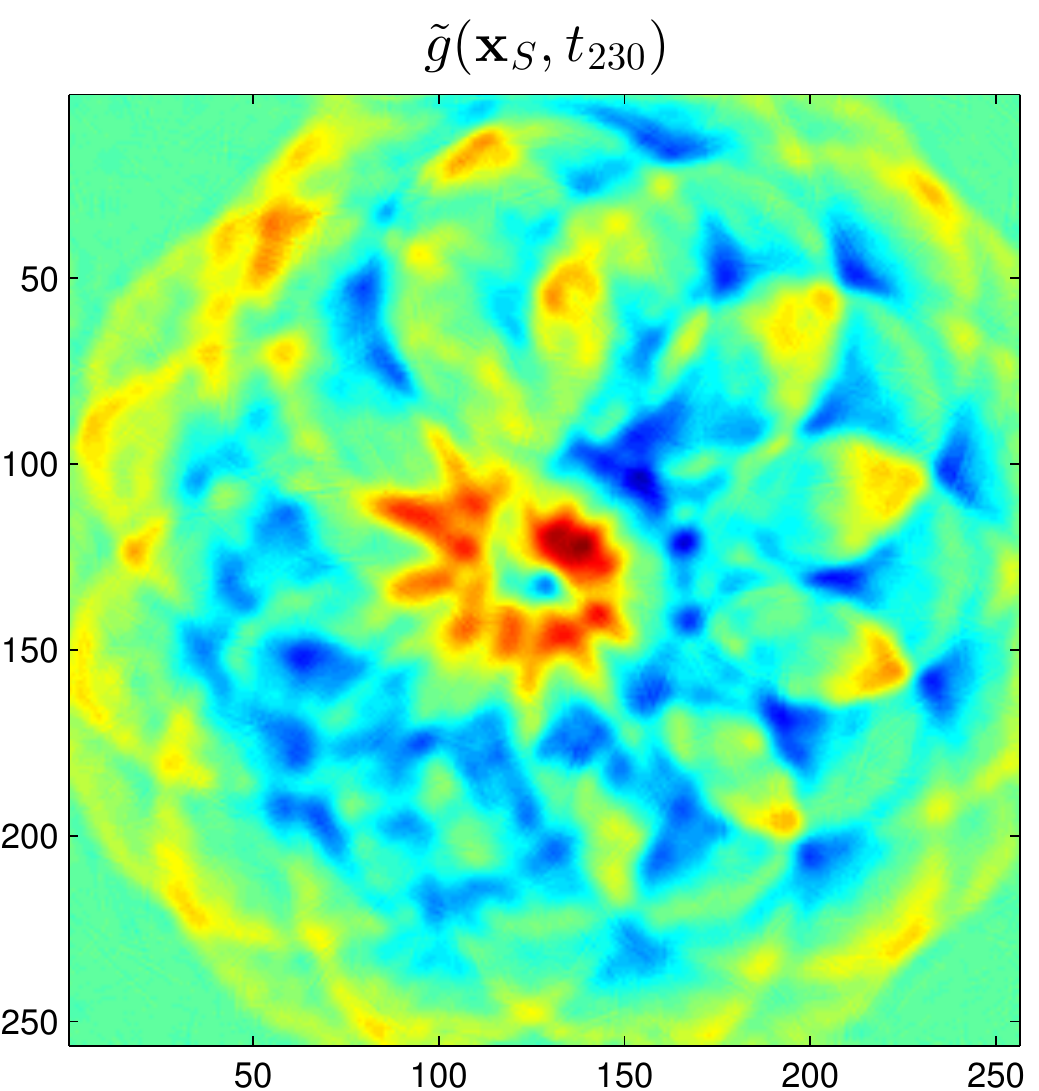}}
\\
\dontshow{
\phantom{\hspace{0.001\textwidth}}
\subfloat[][$\Psi g(\bx_{\mc S},t_{230})$]{
\includegraphics[width=0.13\textwidth,trim={0 0 0em 2.5em},clip]{clock_93x256x256_im_psi_sd_t230_fast_hadamard_11802x65536_curvelet_l3_real_flpcurv192x192_L1_tau0p01_SALSA_v2_stop1_tol0p0005_maxit100_mu500}}
\hspace{0.02\textwidth}
\subfloat[][$\Psi \hat g(\bx_{\mc S},t_{230})$]{
\includegraphics[width=0.13\textwidth,trim={0 0 0em 2.5em},clip]{clock_93x256x256_im_psi_csd_t230_fast_hadamard_11802x65536_curvelet_l3_real_flpcurv192x192_L1_tau0p01_SALSA_v2_stop1_tol0p0005_maxit100_mu500}}
\hspace{0.018\textwidth}
\subfloat[][$\Psi \tilde g(\bx_{\mc S},t_{230})$]{
\includegraphics[width=0.13\textwidth,trim={0 0 0em 2.5em},clip]{clock_93x256x256_im_psi_rsd_t230_fast_hadamard_11802x65536_curvelet_l3_real_flpcurv192x192_L1_tau0p01_SALSA_v2_stop1_tol0p0005_maxit100_mu500}}
\\
}
\phantom{\hspace{0.16\textwidth}}
\subfloat[][Error of $\hat g(\bx_{\mc S},t_{230})$]{ 
\includegraphics[width=0.155\textwidth,trim={0.5em 0 1.5em 2.5em},clip]{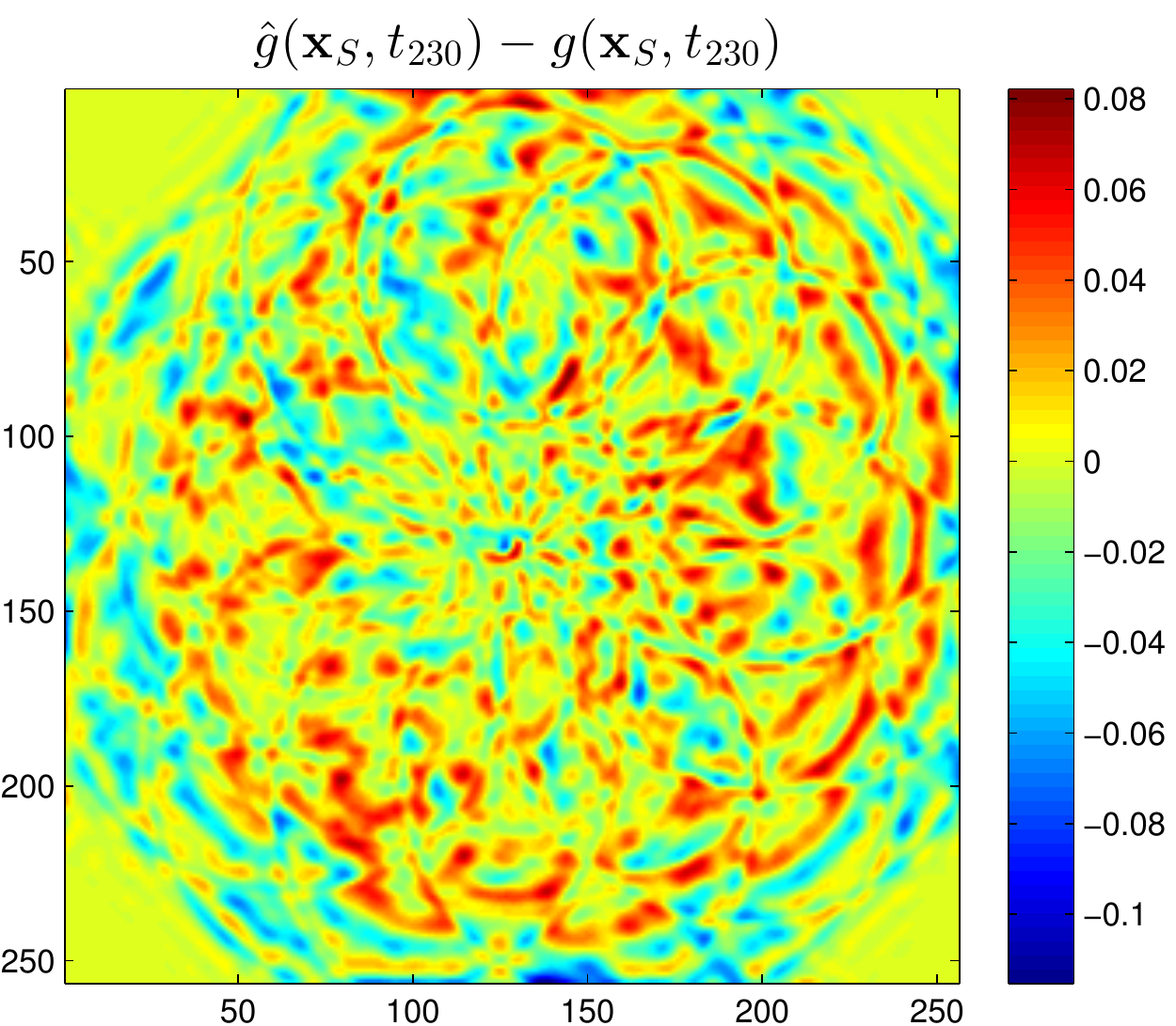}}
\subfloat[][Error of $\tilde g(\bx_{\mc S},t_{230})$]{ 
\includegraphics[width=0.157\textwidth,trim={0.5em 0 1em 2.5em},clip]{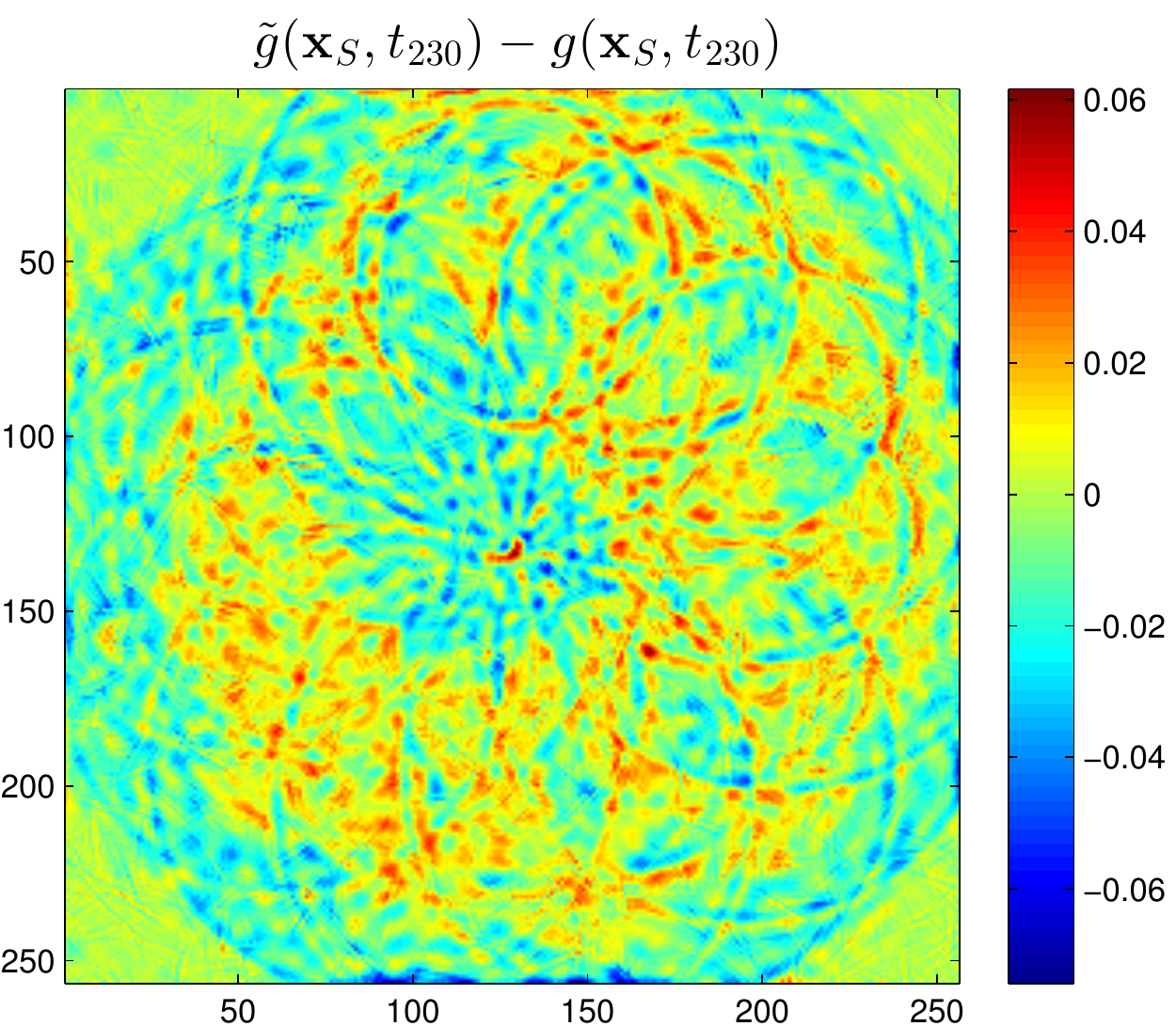}}
\caption{Clock phantom. PAT data at time step $t_{230}$, (a) simulated full data $g(\bx_{\mc S},t)$, (b) compressed $\hat g(\bx_{\mc S},t)$ and (c) reconstructed $\tilde g(\bx_{\mc S},t)$ data with low-frequency Curvelet transform, $\mc C_{3,192,192}^{256,256}$. (d) Compression error $\hat g(\bx_{\mc S},t) - g(\bx_{\mc S},t)$, (e) recovery error $\tilde g(\bx_{\mc S},t) - g(\bx_{\mc S},t)$.}
\label{fig:gt230CLP}
\end{figure}

Next, we compare the compression with $3\%$ of coefficient to recovery from $18\%$ of measurements which corresponds to $6$ times the assumed sparsity, an empirically chosen factor for binary scrambled Hadamard patterns of size $256^2$. 
Figure \ref{fig:mse} shows the mean square error (MSE) of the compressed versus recovered PAT data for both transforms. While the MSEs of the reconstructed data are almost identical, for the data compressed with the low-frequency Curvelets, after some initial time steps the MSE becomes lower than for the standard Curvelets, and most importantly it matches closer the MSE of the reconstructed data. This demonstrates that the low-frequency Curvelets are an adequate (while cheaper) representation of the PAT data. The series of Figures \ref{fig:gt100C},\ref{fig:gt100CLP},
\ref{fig:gt230C},\ref{fig:gt230CLP} shows the PAT data $g(\bx_{\mc S}, t)$, its compression $\hat g(\bx_{\mc S}, t)$ and reconstruction $\tilde g(\bx_{\mc S}, t)$ at different time steps. Consistently, we observe that the higher scale coefficients are eliminated by compression while they partially reappear in the reconstruction hinting that the factor 6 maybe somewhat pessimistic. This is also evident in the higher frequency appearance of the error of the reconstruction in comparison to the compression. 

\begin{figure}[ht]
\subfloat[][$TR(g(\bx_{\mc S},t))$]{ 
\includegraphics[width=0.160\textwidth,trim={0 0 1em 2.5em},clip]{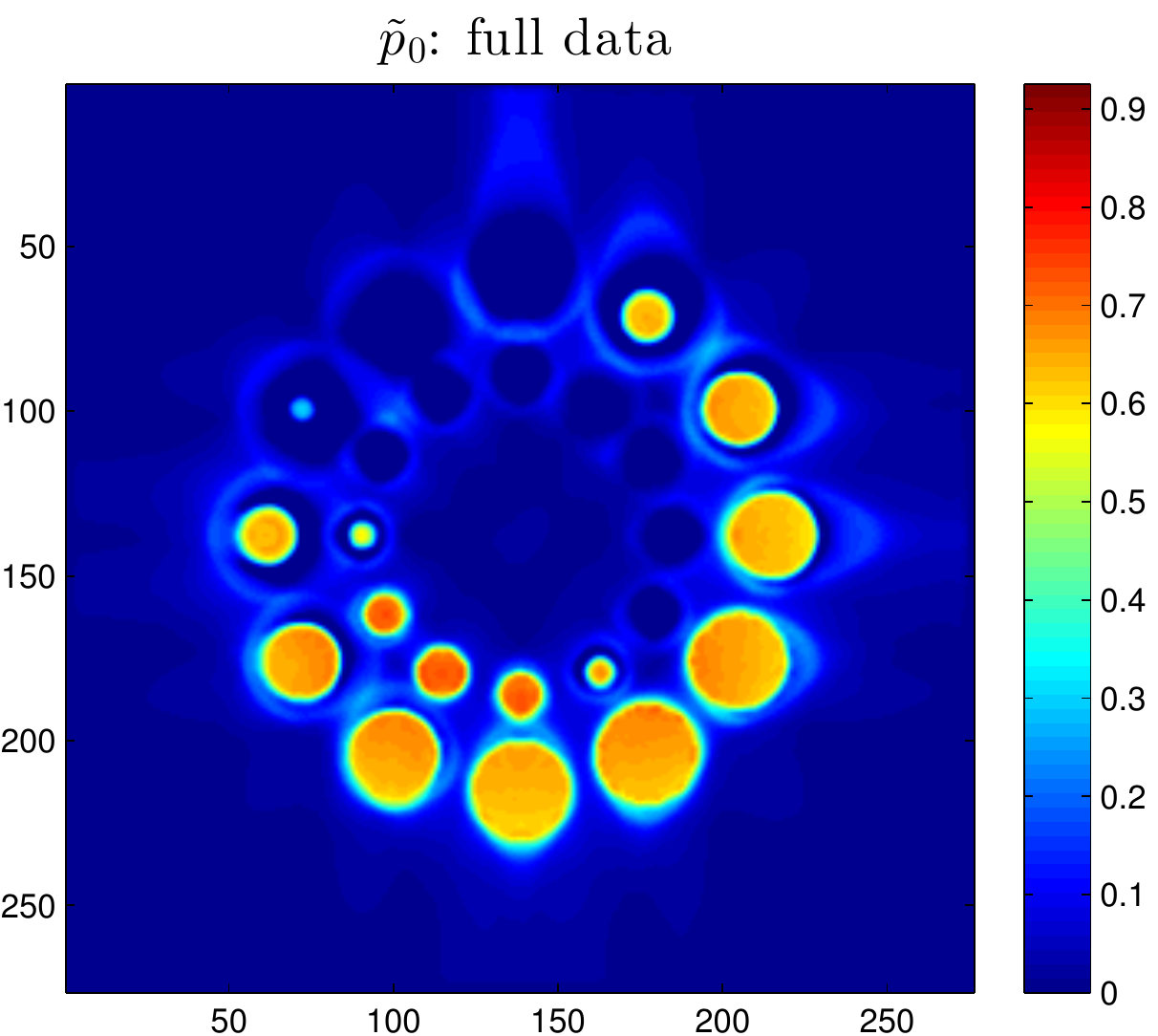}}
\subfloat[][$TR(\tilde g(\bx_{\mc S},t))$]{ 
\includegraphics[width=0.142\textwidth,trim={0 0 0em 2.5em},clip]{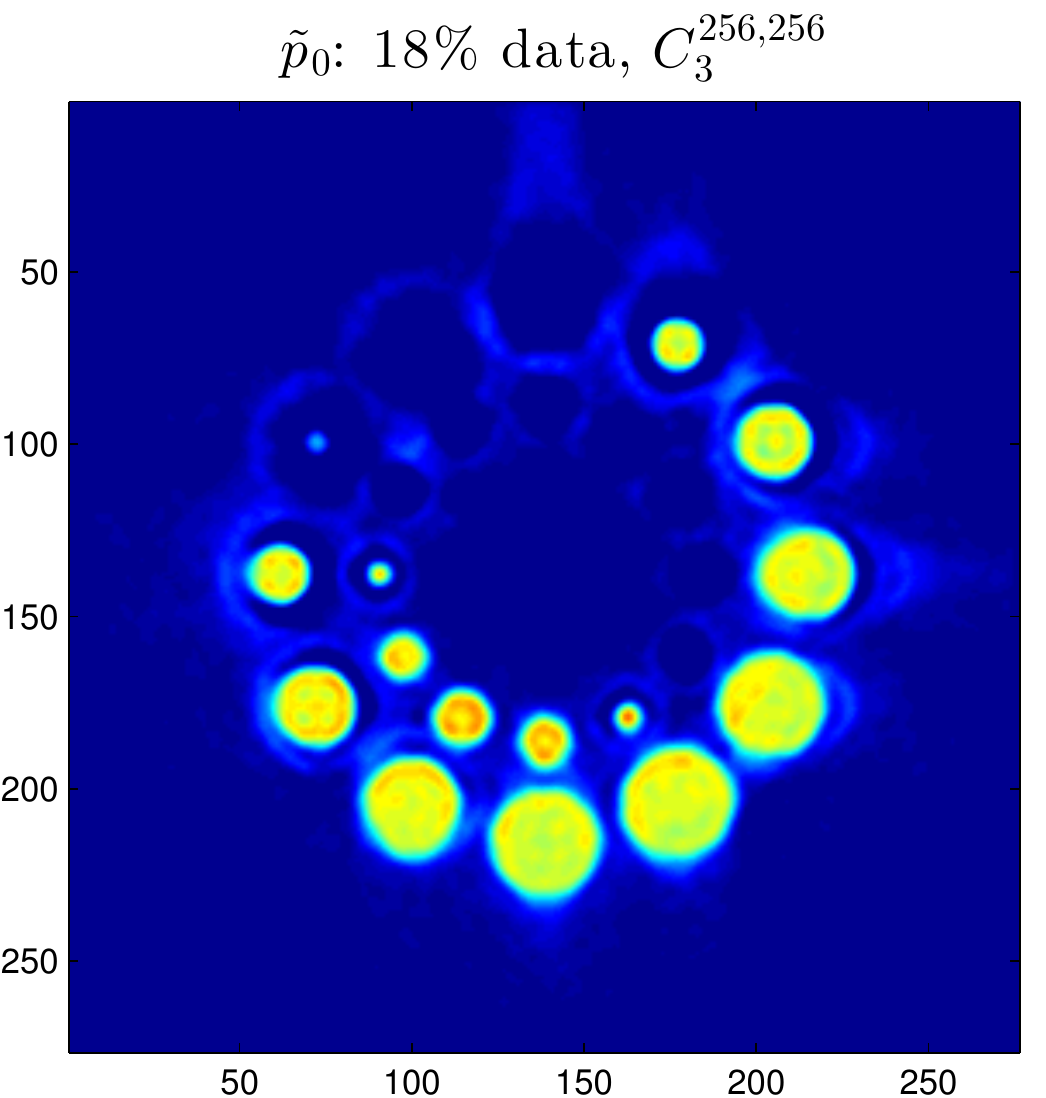}}
\hspace{0.005\textwidth}
\subfloat[][$TR(\tilde g^{LF}(\bx_{\mc S},t))$]{ 
\includegraphics[width=0.142\textwidth,trim={0 0 0em 2.5em},clip]{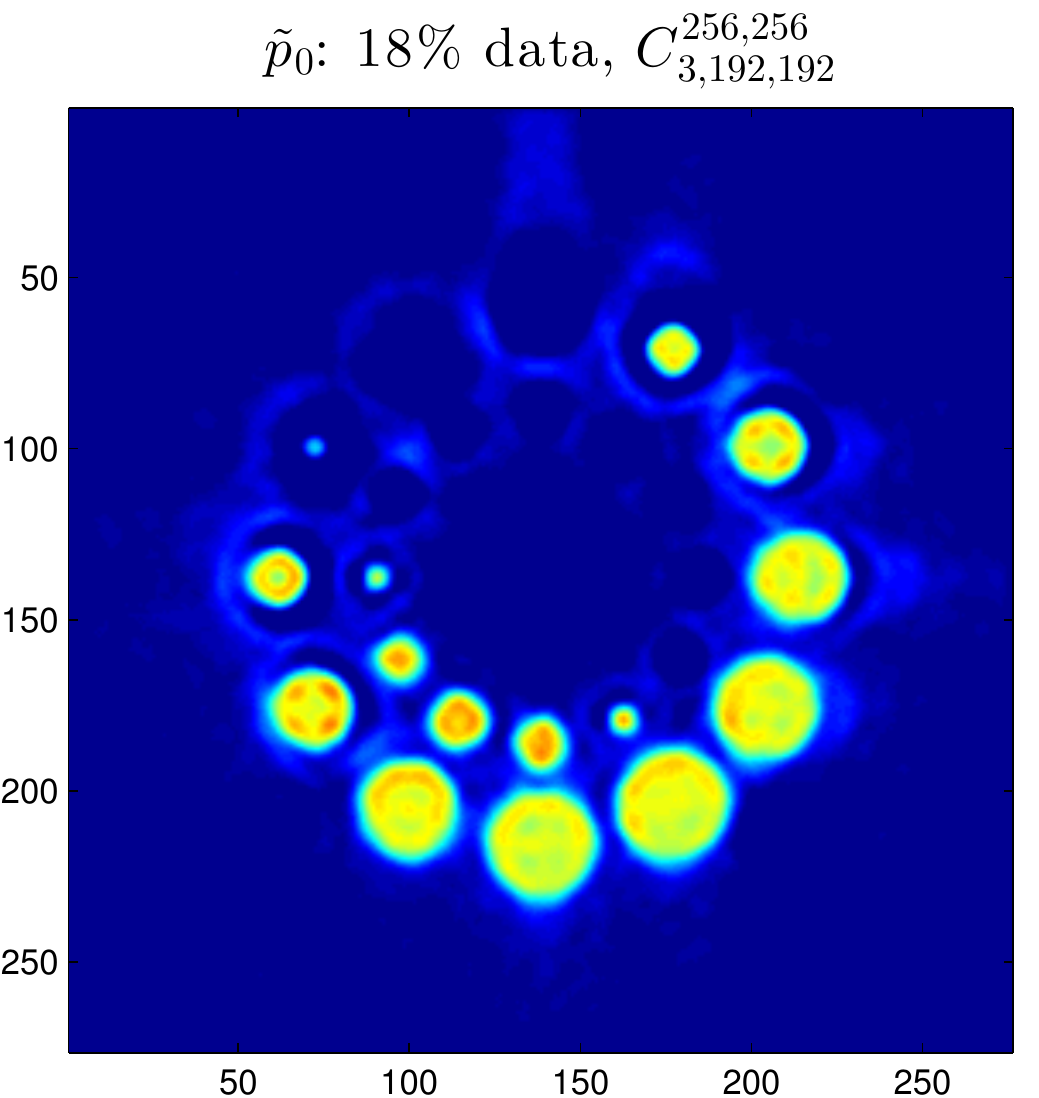}}
\\
\phantom{\hspace{0.16\textwidth}}
\subfloat[][Error (b)-(a)]{ 
\includegraphics[width=0.16\textwidth,trim={0.5em 0 1.5em 2.5em},clip]{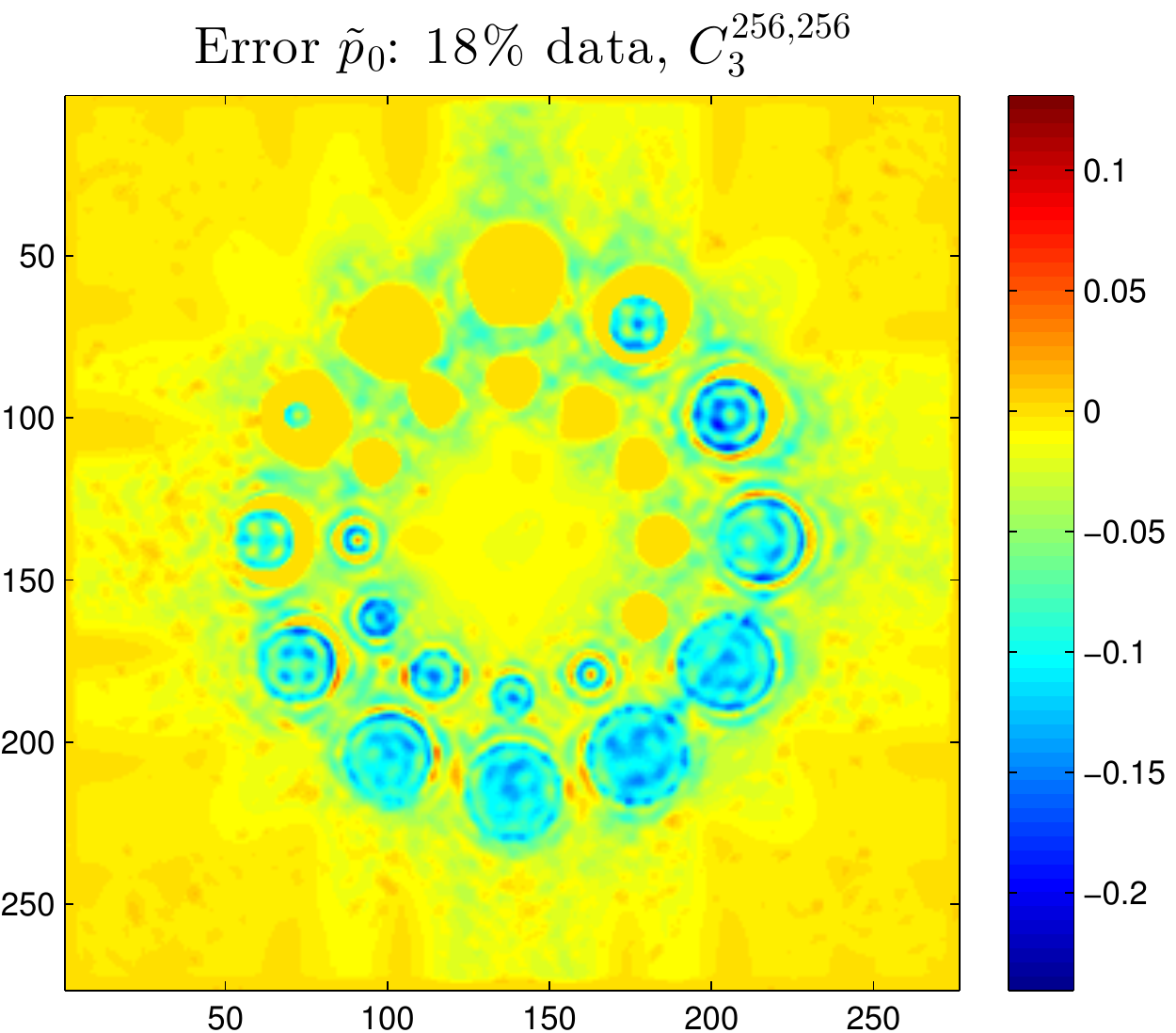}}
\hspace{-0.015\textwidth}
\subfloat[][Error (c)-(a)]{ 
\includegraphics[width=0.142\textwidth,trim={0 0 0em 2.5em},clip]{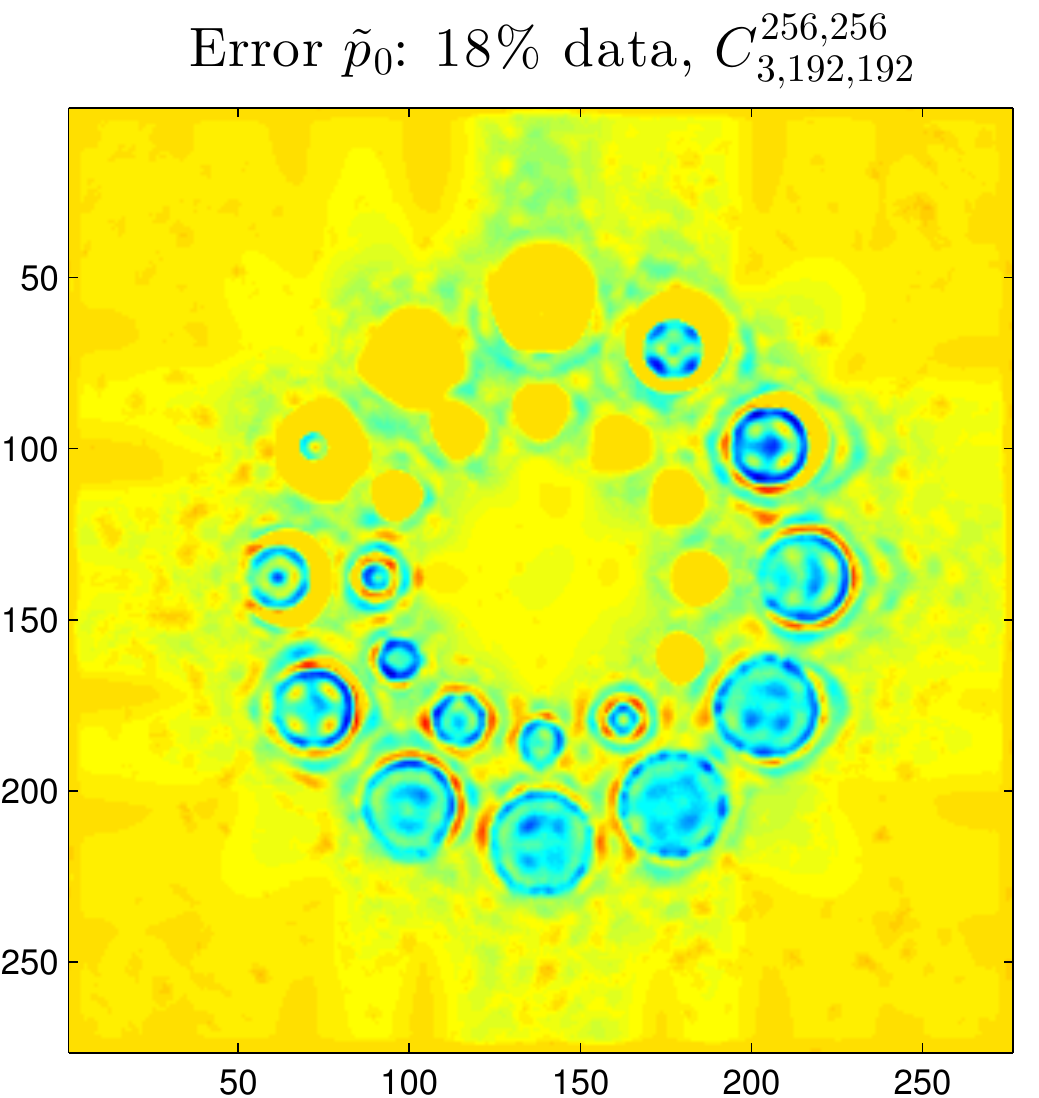}}
\caption{Clock phantom. Central slice through reconstructed PAT image $\tilde p_0$ from (a) full data $g(\bx_{\mc S},t)$, (b) data reconstructed with standard Curvelet transform $\mc C_{3}^{256,256}$, (c) data reconstructed with low-frequency Curvelet transform $\mc C_{3,192,192}^{256,256}$. 
The corresponding error (d),(e).}
\label{fig:p0recon}
\end{figure}

The PAT image reconstruction from the recovered PAT data for both transforms is depicted in Figure \ref{fig:p0recon}, which are visually very similar with MSE of $ 4.4243\cdot 10^{-4}$ for the standard Curvelets and $ 4.6885\cdot 10^{-4}$ for the low-frequency Curvelets, where the reconstruction from full data has been used as the ground truth in MSE calculations. 

\subsection{Synthesized pattern data: knotted tubes filled with ink}
Next, we present reconstruction from compressed measurements synthesized from the point-by-point FP sensor measurements. The purpose of such data is to demonstrate the reconstruction with realistic noise but good signal to noise ratio, which at present is a limiting factor for SPOC. 

\begin{figure}[ht]
\centering
\subfloat[][]{
\includegraphics[width=0.31\textwidth]{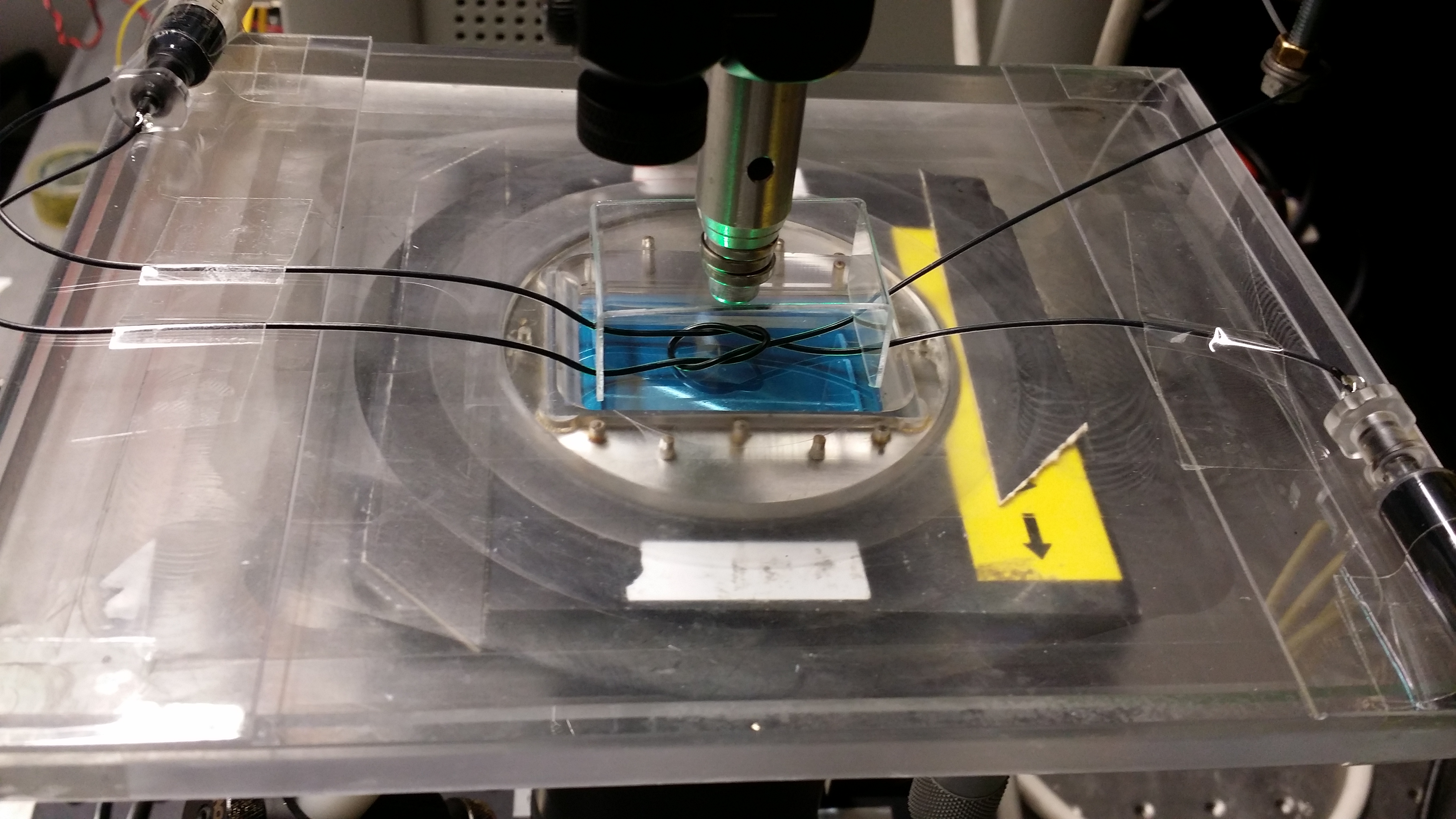}}
\subfloat[][]{
\includegraphics[width=0.18\textwidth]{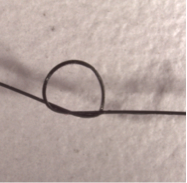}}
\caption{(a) Knotted ink tubes phantom on the FP sensor. (b) Photo of the artificial hair phantom.}\label{fig:Knots}
\end{figure}

Two polythene tubes were filled with 10\% and 100\% ink and tied into a knot, see Figure \ref{fig:Knots}(a). The tubes were immersed in a $1$\% intralipid solution. The wavelength of the excitation laser was $1064$\si{\nano\meter} delivering energy of approximately $20$\si{\milli\joule}. A full scan data consists of $128 \times 128$ locations corresponding to spatial resolution of $150$\si{\micro\meter} $\times$ $150$\si{\micro\meter}, sampled at $625$ time points corresponding to time resolution of $12$\si{\nano\second}. 

\begin{figure}[ht]
\subfloat[][$TR(g(\bx_{\mc S},t))$]{ 
\includegraphics[width=0.160\textwidth,trim={0 0 1em 2.5em},clip]{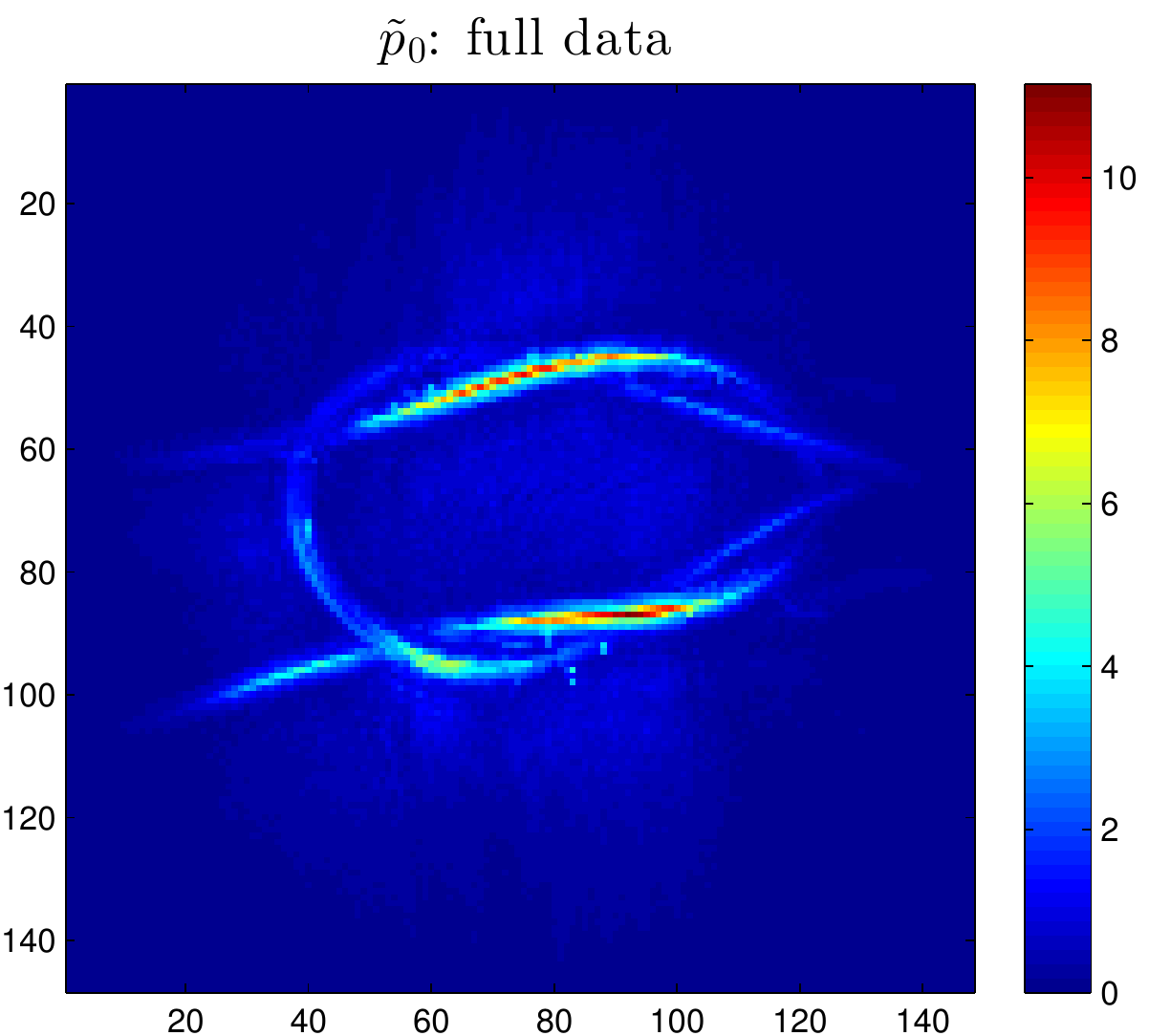}}
\subfloat[][$TR(\tilde g^{lin}(\bx_{\mc S},t))$]{ 
\includegraphics[width=0.142\textwidth,trim={0 0 0em 2.5em},clip]{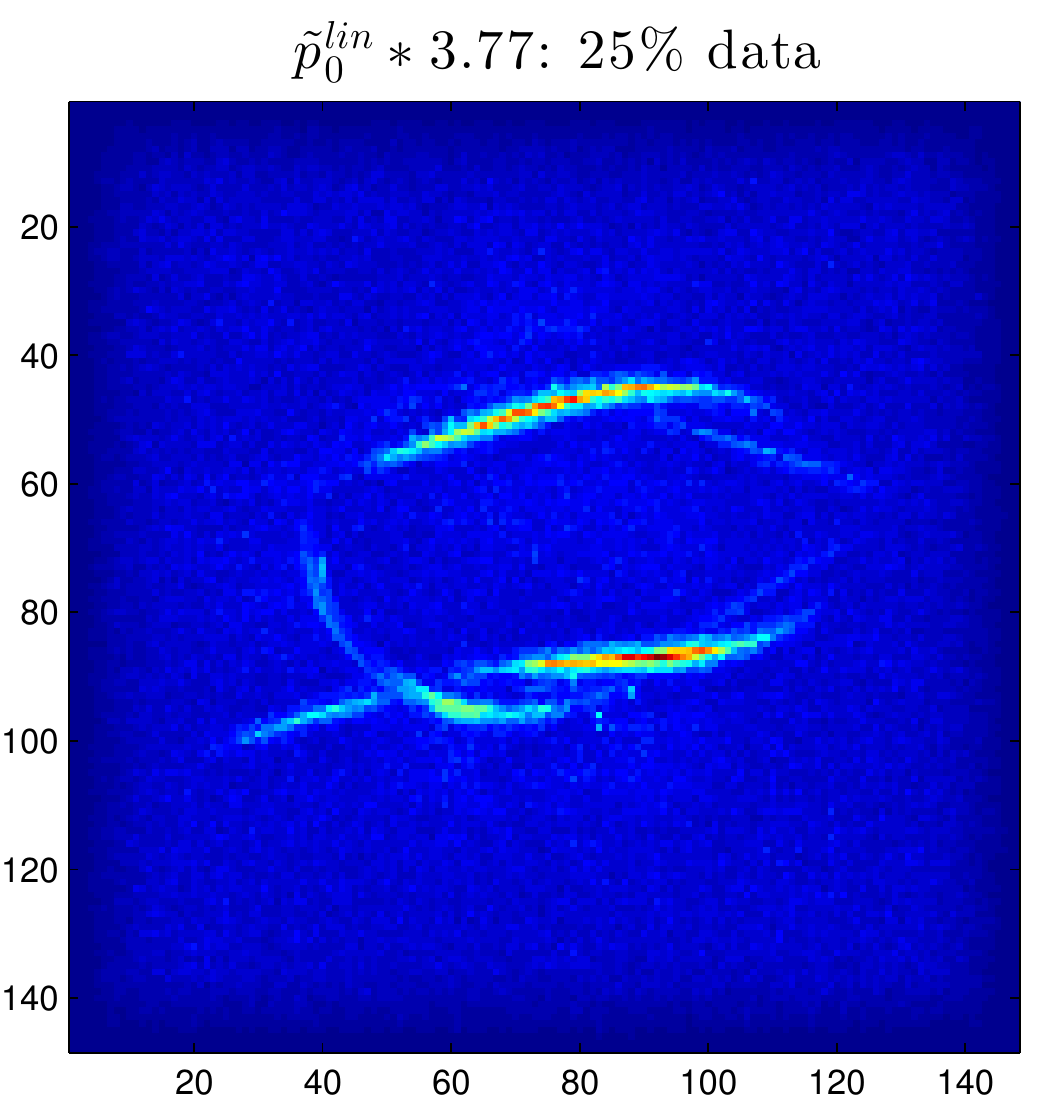}}
\hspace{0.008\textwidth}
\subfloat[][$TR(\tilde g(\bx_{\mc S},t))$]{ 
\includegraphics[width=0.142\textwidth,trim={0 0 0em 2.5em},clip]{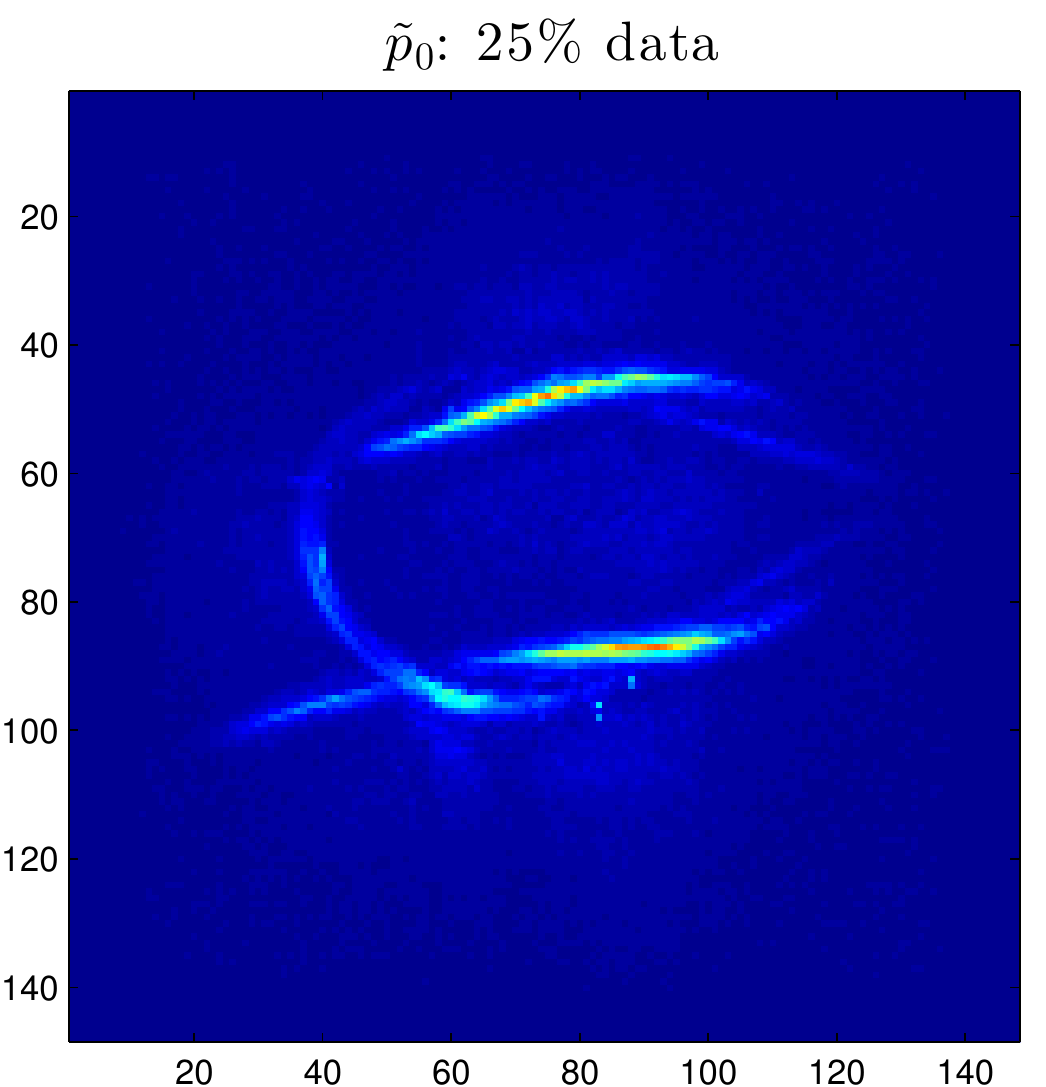}}
\\
\phantom{\hspace{0.16\textwidth}}
\subfloat[][Error (b) - (a)]{ 
\includegraphics[width=0.160\textwidth,trim={0 0 1em 2.5em},clip]{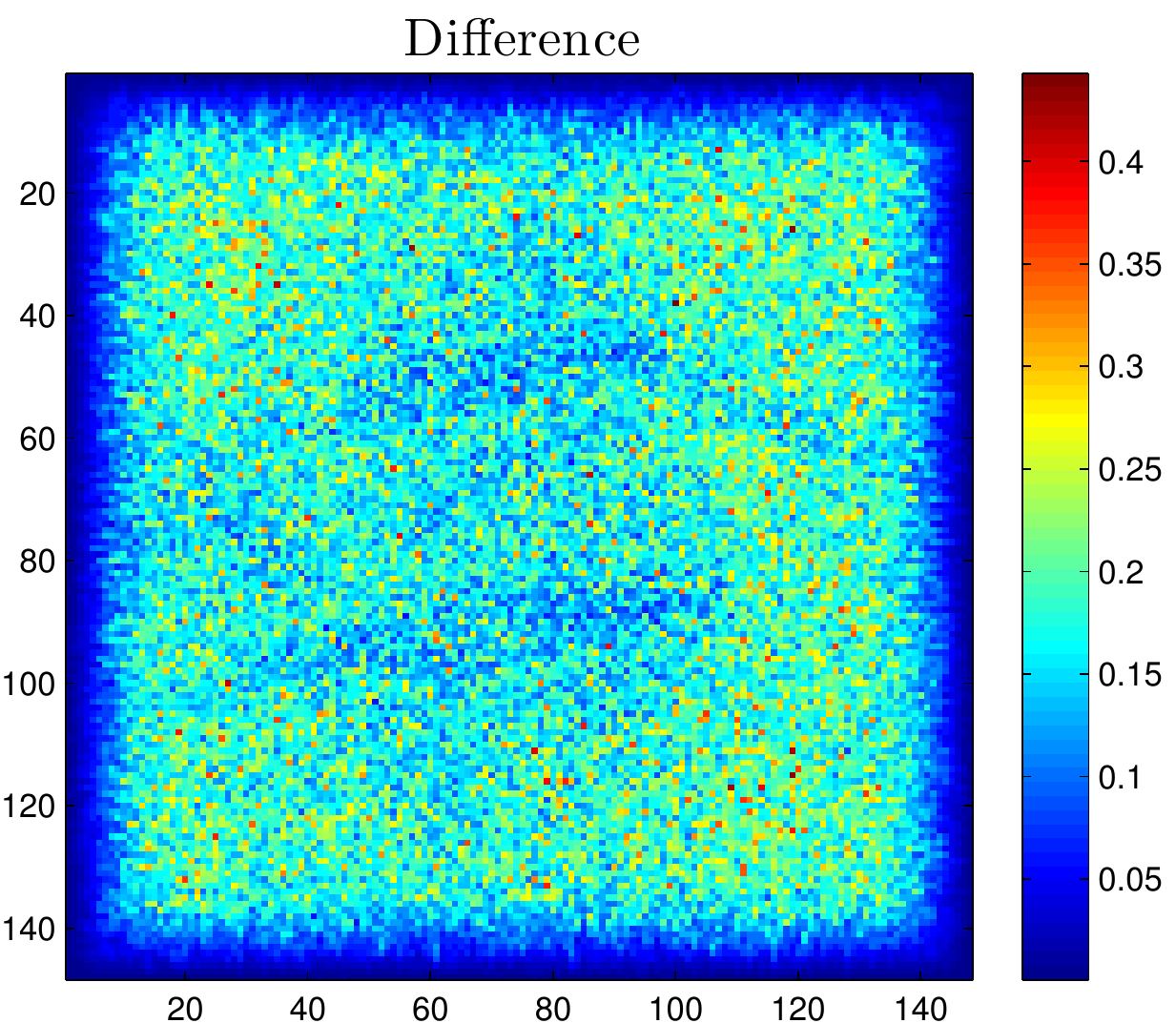}}
\subfloat[][Error (c) - (a)]{ 
\includegraphics[width=0.1575\textwidth,trim={0 0 1em 2.5em},clip]{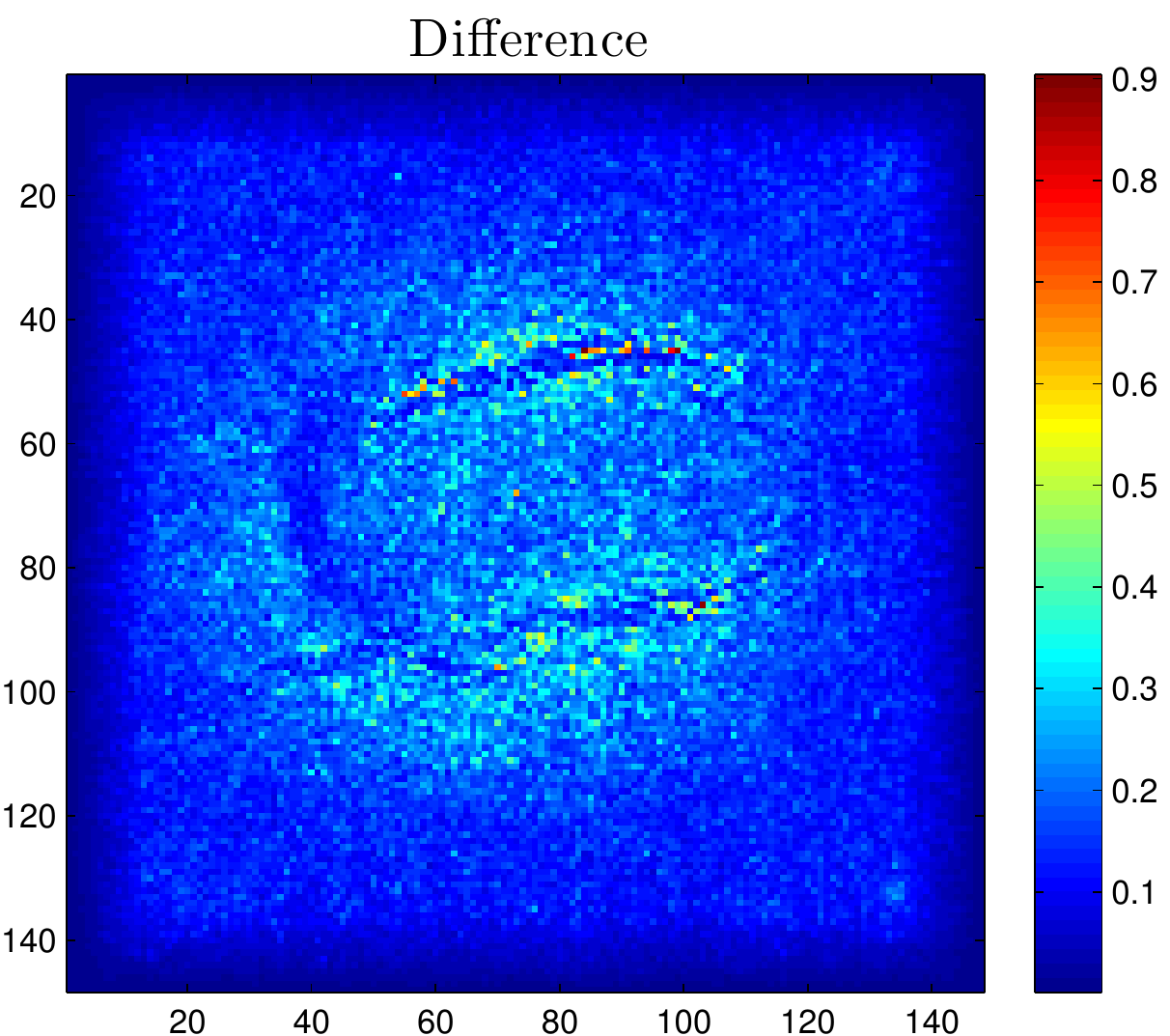}}
\caption{Knotted tubes. Maximum intensity projection of reconstructed PAT image $\tilde p_0$ from (a) full data $g(\bx_{\mc S},t)$, (b) linearly reconstructed data $\tilde g^{lin}(\bx_{\mc S},t)$ (here plotted scaled by a factor $3.77$), (c) nonlinearly reconstructed data $\tilde g(\bx_{\mc S},t)$ with low-frequency Curvelet transform $\mc C_{3,192,192}^{256,256}$, 
and (d),(e) the corresponding error.}
\label{fig:ThKn_p0recon}
\end{figure}

As a gold standard we take the reconstruction from a full point-by-point data set shown in Figure \ref{fig:ThKn_p0recon}(a). 
%
The linear reconstruction from 25\% of patterns obtained by setting the missing Hadamard pattern measurements to 0 is shown in Figure \ref{fig:ThKn_p0recon}(b) while the nonlinear reconstruction with SALSA (with parameters $\tau = 0.01 \max(|\Psi\Phi\tr b^t|)$, $\mu = 5 \max(|\Psi\Phi\tr b^t|)/\|b^t\|_2$ and stopping tolerance $5\cdot 10^{-4}$ or after $100$ iterations) in Figure \ref{fig:ThKn_p0recon}(c). The nonlinear reconstruction effectively restores the contrast lost in the linear 0-padded reconstruction.
The MSE with respect to the full data reconstruction is $6.7582\cdot 10^{-3}$ for the linear 0-padded reconstruction and $2.0824\cdot 10^{-3}$ for the nonlinear reconstruction.
\begin{figure}[ht]
\centering
\includegraphics[width=0.225\textwidth]{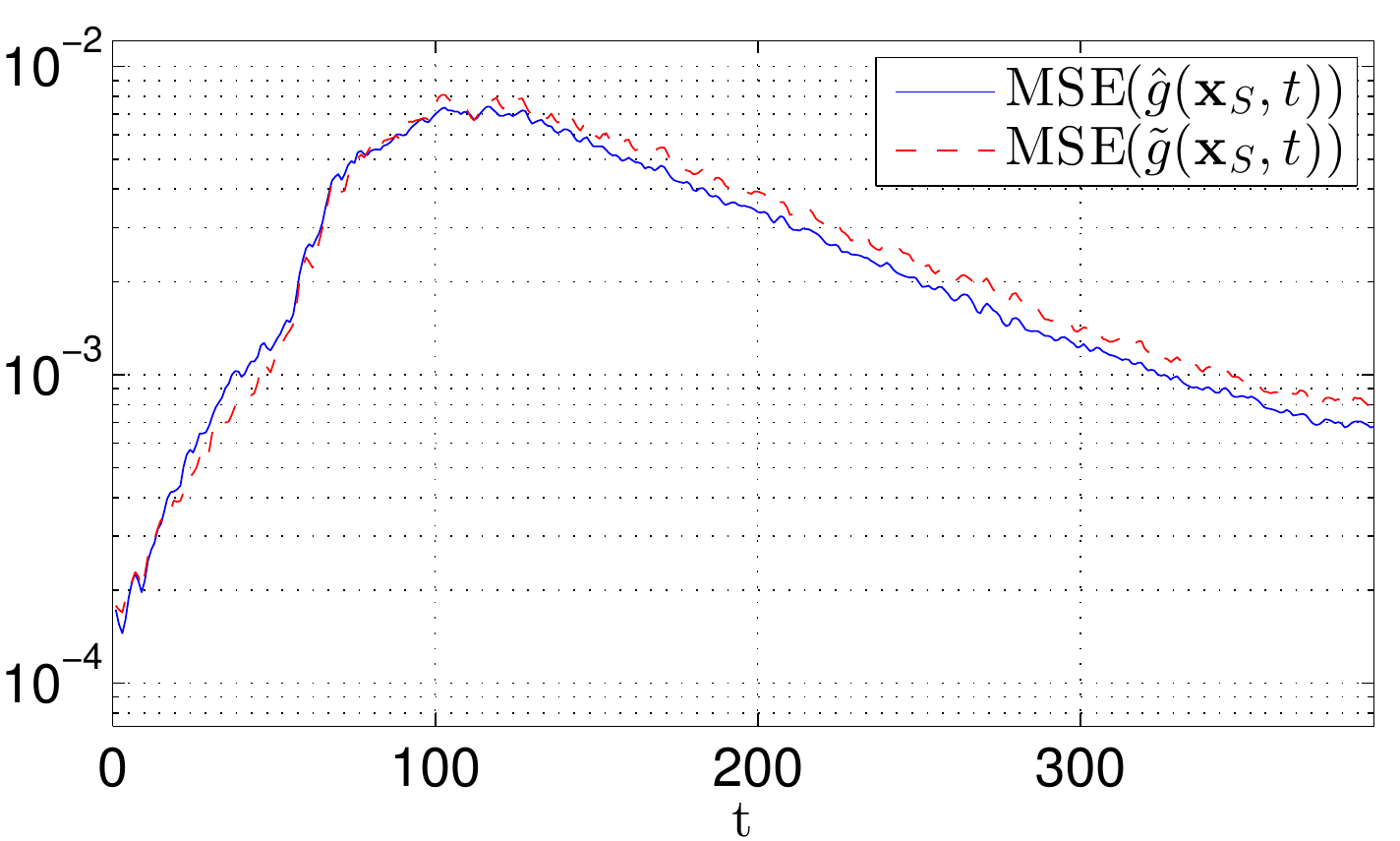}
\caption{Knotted tubes. MSE of the compressed PAT data $\hat g(\bx_{\mc S},t)$ versus the reconstructed PAT data $\tilde g(\bx_{\mc S},t)$ over time for low-frequency Curvelet transform $\mc C_{3,96,96}^{128,128}$.}
\label{fig:ThKn_mse}
\end{figure}

\subsection{Experimental pattern data: hair knot}
Our last example is a synthetic hair knot phantom of diameter $\sim$150 \si{\micro\meter}, immersed in 1\% intralipid solution and positioned approximately $2$\si{\milli\meter} above the sensor and $3$\si{\milli\meter} deep below the intralipid surface. The photo of the phantom is shown in Figure \ref{fig:Knots}(b). 
The area of DMD corresponding to $640 \times 640$ micromirrors grouped in $5 \times 5$ was used to form $128 \times 128$ pixels. Due to the angle of the optical path, each such pixel corresponds to an area of $62.12$ $\times$ $68$ \si{\micro\meter} on the FP sensor. The measurement was averaged over four excitations. 
The speed of sound used for time reversal was 1490 \si{\meter/\second}.

\begin{figure}[ht]
\subfloat[][$TR(g(\bx_{\mc S},t))$]{ 
\includegraphics[width=0.164\textwidth,trim={0 0 1em 2.5em},clip]{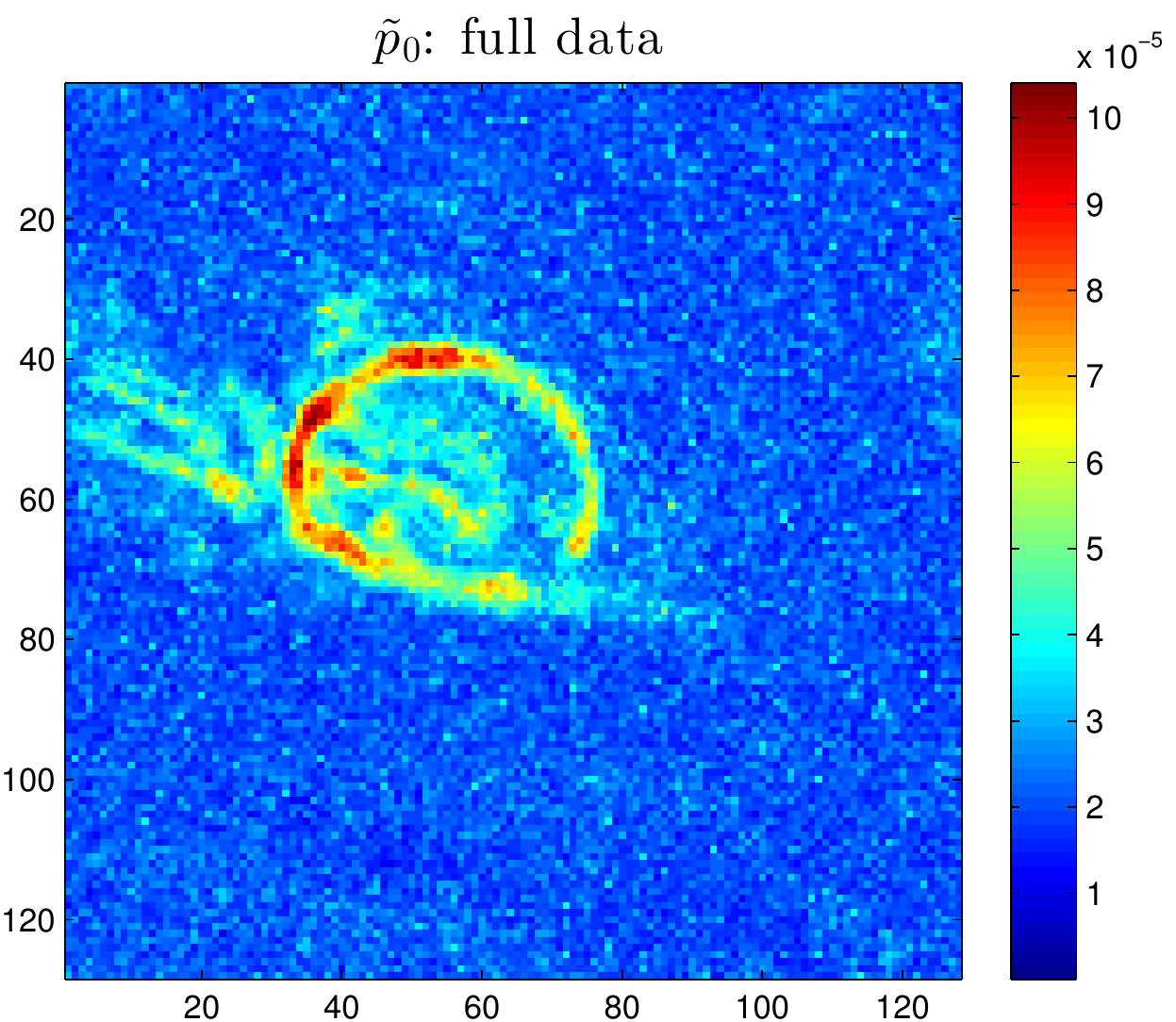}}
\subfloat[][$TR(\tilde g^{lin}(\bx_{\mc S},t))$]{ 
\includegraphics[width=0.142\textwidth,trim={0 0 0em 2.5em},clip]{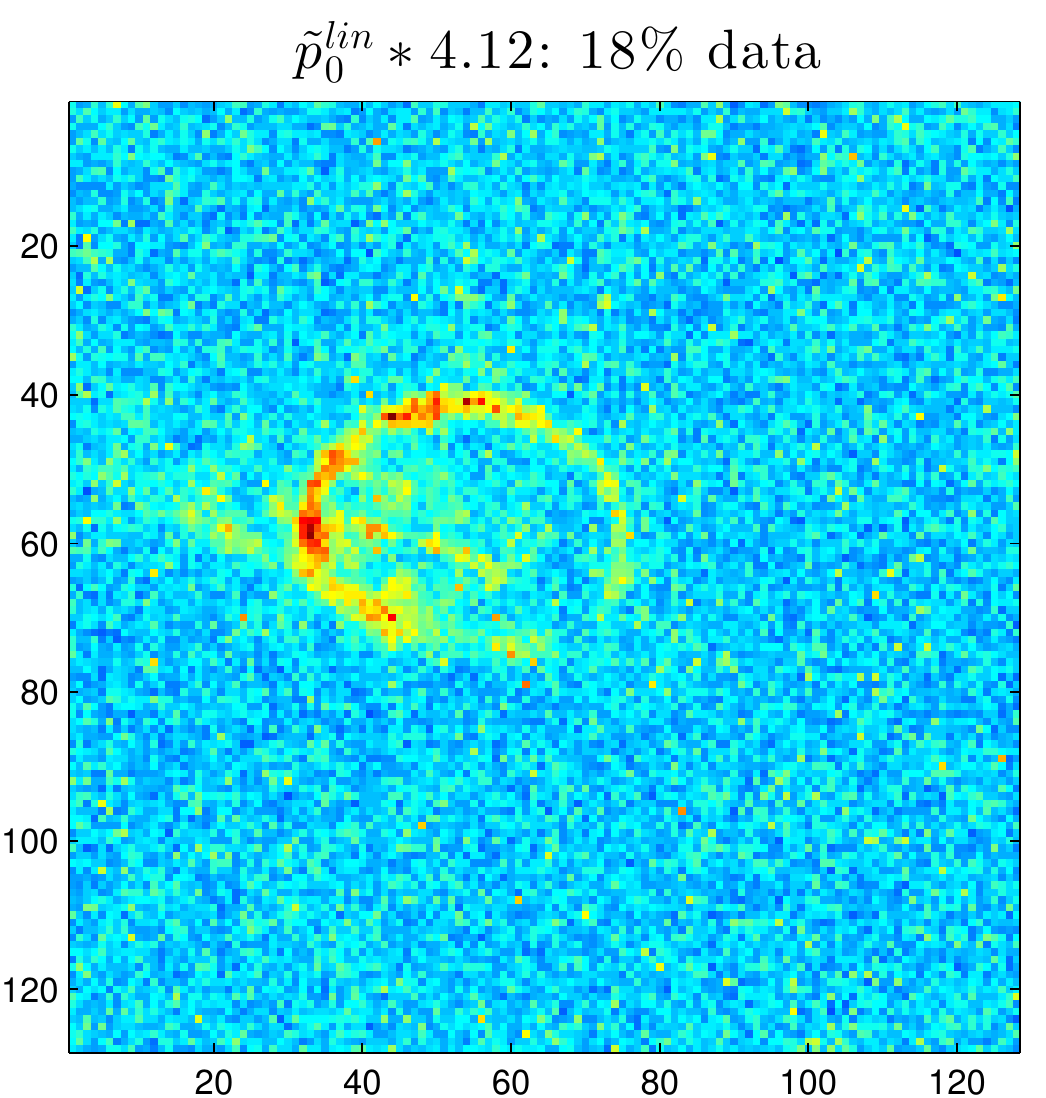}}
\hspace{0.008\textwidth}
\subfloat[][$TR(\tilde g(\bx_{\mc S},t))$]{ 
\includegraphics[width=0.142\textwidth,trim={0 0 0em 2.5em},clip]{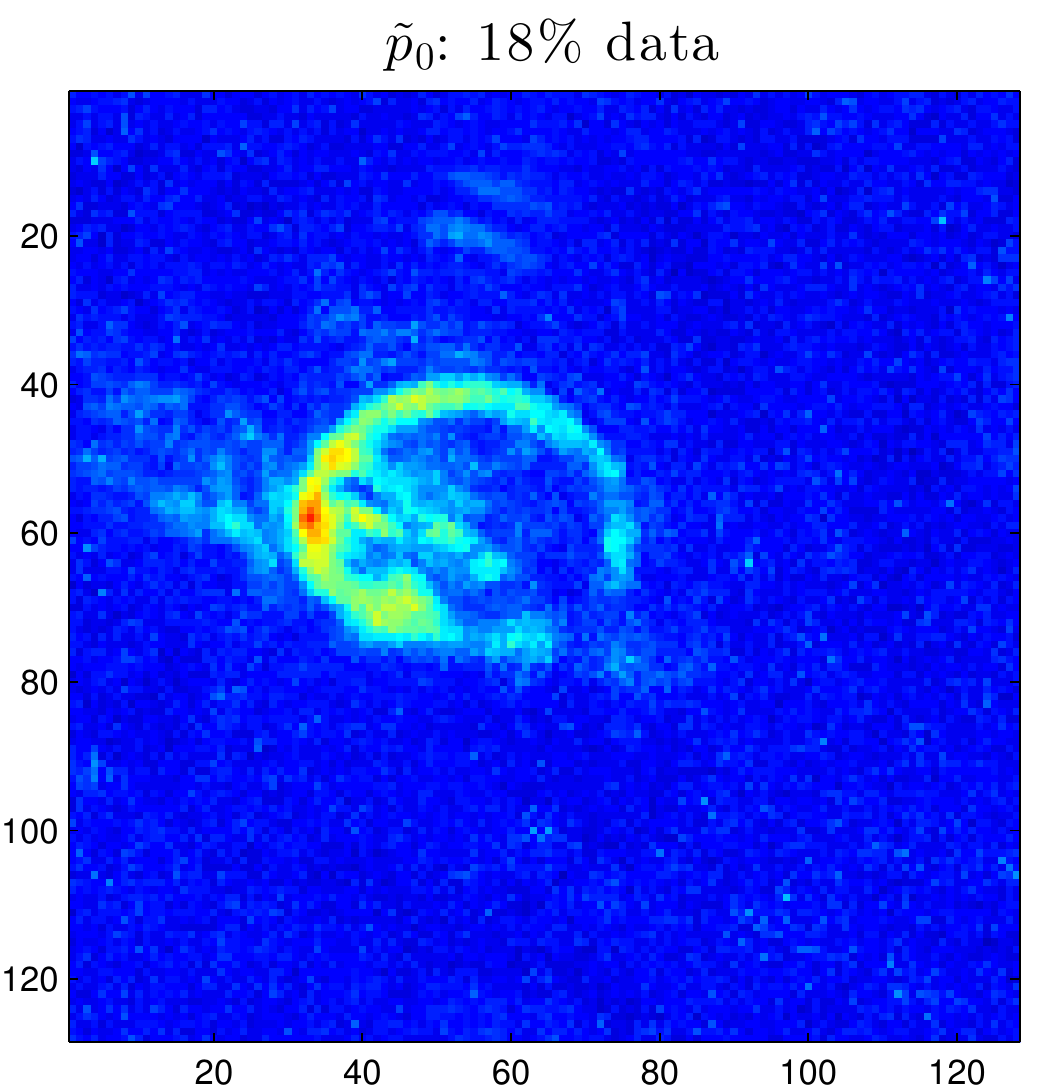}}
\caption{Hair knot. Maximum intensity projection of reconstructed PAT image $\tilde p_0$ from (a) full data $g(\bx_{\mc S},t)$, (b) linearly reconstructed data $\tilde g^{lin}(\bx_{\mc S},t)$ (here plotted scaled by a factor $4.12$), (c) nonlinearly reconstructed data $\tilde g(\bx_{\mc S},t)$ with low-frequency Curvelet transform $\mc C_{3,96,96}^{128,128}$.}
\label{fig:HrKn_p0recon}
\end{figure}

Figure \ref{fig:HrKn_p0recon}(a) shows the PAT image reconstruction obtained from the full set of scrambled binary Hadamard $128^2$ patterns (the full point data was computed inverting scrambled Hadamard transform). The $\sim 4.12$ amplified linear reconstruction from 0-padded 18\% measurements is shown in Figure \ref{fig:HrKn_p0recon}(b) and the nonlinear reconstruction obtained with SALSA ($\tau = 0.05\max(\Psi\Phi\tr b^t)$, $\mu = 0.75\max(\Psi\Phi\tr b^t)/\|b^t\|_2$ and the stopping tolerance $5\cdot 10^{-3}$ or maximum $100$ iterations) is illustrated in Figure \ref{fig:HrKn_p0recon}(c). 

\section{Conclusions and discussion}
We presented a method for acoustic field reconstruction from a limited number of patterned measurements from an optical ultrasound detector. When compressively sensing photoacoustic signals, our method 
recovered the pressure field on the sensor (the PAT data)
at each time step independently using the sparsity of the data in a low-frequency Curvelet frame, which is a modification of the standard Curvelets tailored to account for the smoothing of the wave front during optical acquisition. The major advantage of such a scheme is that 
the series of problems to solve are standard (2D) CS recovery problems and as they are independent they can be solved in parallel.  
Furthermore, decoupling the CS and acoustic reconstruction affords more flexibility in PAT modeling, for instance including absorption or nonlinear effects, and allows for the use of highly optimized readily available software for 
photoacoustic image reconstruction.

One of the major challenges for compressed sensing of photoacoustic signals is that the same number of interrogation patterns is applied in each time step. As the sparsity of the wave field on the detector varies, in the proposed scheme this inevitably results in different quality reconstruction in different time steps. In particular, we first lose the ``sharpness'' of the wave front, because this is reflected in the  coefficients at the highest scale which magnitude is generally smaller. This effect is partially counterweight by the proposed low-frequency Curvelet representation which is tailored to frequency range of photoacoustic data and hence boosts those coefficients. 
\dontshow{
One of the major challenges for compressed sensing of photoacoustic signals is the time varying sparsity of the wave field on the detector, which when sampled with the same number of interrogation patterns inevitably results in different quality reconstruction in different time steps. 
In particular, we first lose the definition (sharpness) of the wave fronts, because this is reflected in very small coefficients at the highest scale. To counterweight this effect we 
proposed a low-frequency Curvelet representation which is tailored to frequency range of photoacoustic data and hence allows for higher compressibility. 
}
Furthermore, in our experiments we observed that the data recovery errors were partially alleviated during the acoustic inversion. As the initial pressure $p_0$ is mapped to the entire time series, the entire time series carries the information about the wave front and the acoustic inversion acts to average out data errors. 
In \cite{Arridge:2016adj, Arridge:2016apat} this problem was tackled utilizing the sparsity directly in the photoacoustic image $p_0$, rather than in the data. There we solve one large (3D) CS recovery problem where the CS sensing operator is a composition of the acoustic propagation and pattern measurements. Such forward operator is significantly more expensive to apply and its incoherence properties have to be analyzed. Finally, nonlinear effects such as acoustic absorption are out of scope of the standard linear CS framework. 

In future work, we intend to extend the here proposed acoustic field reconstruction method to efficiently handle dynamic problems.


%



\section*{Acknowledgment}

This research was supported by Engineering and Physical Sciences Research Council, UK (EP/K009745/1).

\ifCLASSOPTIONcaptionsoff
  \newpage
\fi



\bibliographystyle{IEEEtran}
\bibliography{PATCS}
\end{document}